\documentclass[11pt]{article}

\usepackage{amssymb,amsmath,amsthm,amsopn,amsfonts}
\usepackage[T1]{fontenc}
\usepackage{ae,aecompl}
\usepackage{graphicx}

\usepackage{bm}
\usepackage{fancybox,layout,color}
\usepackage{fancyhdr}

\usepackage{scrtime}

\usepackage{setspace}
\usepackage{url}
\usepackage{float}
\usepackage{placeins}

\numberwithin{equation}{section}
\definecolor{purple}{RGB}{160,32,40}

\usepackage{ulem}  % for the command \sout{}
\usepackage[makeroom]{cancel}  % for the command \cancel{}

\newtheorem{teo}{Theorem}[section]
\newtheorem{nota}[teo]{Remark}
\newtheorem{ex}[teo]{Example}
\newtheorem{coro}[teo]{Corollary}
\newtheorem{defi}[teo]{Definition}
\newtheorem{prop}[teo]{Proposition}
\newtheorem{lema}[teo]{Lemma}

\newcommand{\R}{\ensuremath{\mathbb{R}} }

\newcommand{\dist}{\mathrm{dist}}

\newcommand{\epi}{\mathrm{epi}}
\newcommand{\Aff}{\mathrm{Aff}}

\DeclareMathOperator{\co}{\mathsf{co}}

\DeclareMathOperator{\diam}{\mathsf{diam}}
\DeclareMathOperator{\sprt}{\mathsf{spprt}}

\title{ Compensated Convexity, Multiscale Medial Axis Maps\\ and Sharp Regularity of the Squared Distance Function} 
 
\author{\normalsize  Kewei Zhang\thanks{School  of Mathematical Sciences, 
				University of Nottingham, University Park, Nottingham, NG7 2RD, UK, (Kewei.Zhang@nottingham.ac.uk)},  
	Elaine Crooks\thanks{Department of Mathematics, Swansea University, 
				Singleton Park, Swansea, SA2 8PP, UK (e.c.m.crooks@swansea.ac.uk)}
		and
	Antonio Orlando\thanks{CONICET, Inst. de Estructuras \& Dept. de Mec\'{a}nica, Universidad Nacional de Tucum\'an,
Argentina (aorlando@herrera.unt.edu.ar)}
}

\date{ }

\begin{document}

%%%%%%%%%%%%%%%%%%%%%%%%%%%%%%%

\maketitle

%%%%%%%%%%%%%%%%%%%%%%%%%%%%%%%

\singlespacing 
\pagestyle{fancy}
\fancyhead{}
% \fancyfoot[OR,ER]{\tiny \today\,\, \thistime}
\cfoot{\thepage}

%%%%%%%%%%%%%%%%%%%%%%%%%%%%%%%%

\begin{abstract}  
In this paper we introduce a new stable mathematical model for locating and 
measuring the medial axis of geometric 
objects, called the {\it quadratic multiscale medial axis map of scale $\lambda$},
and provide a sharp regularity result for the squared-distance function 
to any closed non-empty subset $K$ of $\R^n$. 
Our results exploit  properties of the function $C^l_{\lambda}(\dist^2(\cdot;\, K))$
obtained by applying the quadratic lower compensated convex transform of parameter $\lambda$ \cite{Zha08a} 
to $\dist^2(\cdot;\, K)$, the Euclidean squared-distance function to $K$. 
Using a quantitative estimate 
for the tight approximation of $\dist^2(\cdot;\, K)$  by $C^l_{\lambda}(\dist^2(\cdot;\, K))$, we prove the 
$C^{1,1}$-regularity of $\dist^2(\cdot;\, K)$ outside a neighbourhood of the closure of the medial axis $M_K$ of $K$, 
which can be viewed as a weak Lusin-type theorem for $\dist^2(\cdot;\, K)$, and give an asymptotic expansion formula 
for $C^l_{\lambda}(\dist^2(\cdot;\, K))$  in terms of 
the scaled squared distance transform to the set and to the convex hull of the 
set of points that realize the minimum distance to $K$. 
The multiscale medial axis map, denoted by $M_{\lambda}(\cdot;\, K)$, is a family  of non-negative 
 functions, parametrized by  $\lambda>0$, whose limit as $\lambda \to \infty$ exists and is called 
 the multiscale medial axis landscape map, $M_{\infty}(\cdot;\, K)$. We show that $M_{\infty}(\cdot;\, K)$
 is   strictly positive on the medial axis $M_K$
and zero elsewhere. We give conditions that ensure $M_{\lambda}(\cdot;\, K)$  keeps a constant height
along the parts of $M_K$ generated by two-point subsets with the value of the height dependent on
the scale of the distance between the generating points, thus providing a
hierarchy of heights (hence, the word 'multiscale') between different parts of $M_K$ that 
enables subsets of $M_K$ to be selected by simple thresholding.
Asymptotically, further understanding of the multiscale effect is provided by our exact representation of $M_{\infty}(\cdot;\, K)$.
Moreover, given a compact subset $K$ of $\R^n$, while it is well known that 
$M_K$ is not Hausdorff stable, we prove that in contrast, 
$M_{\lambda}(\cdot;\, K)$  is stable under the
 Hausdorff distance, and deduce implications  for the localization of the stable parts of $M_K$. 
Explicitly calculated prototype examples of  medial axis maps are also presented and used to illustrate the 
theoretical findings. 

\end{abstract}

\footnotesize{

\medskip
{\bf 2000 Mathematics Subjects Classification number}:
 53A05,  26B25, 52B55,  52A41, 65D17,  65D18

\medskip
{\bf Email}:  kewei.zhang@nottingham.ac.uk,    e.c.m.crooks@swansea.ac.uk, aorlando@herrera.unt.edu.ar
}

\normalsize

%%%%%%%%%%%%%%%%%%%%%%%%%%%%%%%%%%%%%%%%%%%%%%%%%%%%%%%%%%%%%%%%%%%%%%%%%%%%%%%%%%%%%%%%%%%%%%%%%%%%%
%%%%%%%%%%%%%%%%%%%%%%%%%%%%%%%%%%%%%%%%%%%%%%%%%%%%%%%%%%%%%%%%%%%%%%%%%%%%%%%%%%%%%%%%%%%%%%%%%%%%%
%%%%%%%%%%%%%%%%%%%%%%%%%%%%%%%%%%%%%%%%%%%%%%%%%%%%%%%%%%%%%%%%%%%%%%%%%%%%%%%%%%%%%%%%%%%%%%%%%%%%%

\section{Introduction}
The medial axis of an object is a geometric structure that was introduced by Blum \cite{Blu67} as a
means of providing a compact representation of a shape. 
Initially defined as the set of the shock points of a grass fire lit on the boundary and allowed to propagate uniformly
inside the object, closely related definitions of skeleton \cite{CH68} and cut-locus \cite{Wol93}
have since been proposed, and have 
served for the study of its topological properties \cite{ACNS13,CCM91,Lie04,Mat88,SPW96}, 
its stability \cite{ChoiS04,CS04} and for the development of fast and efficient algorithms for its computation 
\cite{AAAHJR09,AM97,AL13,KSKB95,OK95}. 
Applications of the medial axis are ample in scope and nature, ranging from computer vision 
to image analysis, from mesh generation to computer aided design. We refer to \cite{SP08} and the 
references therein for applications and accounts of some recent theoretical developments. 

An inherent drawback of the medial axis is, however, its sensitivity to boundary details, in the sense that 
small perturbations of the object (with respect to the Hausdorff distance) can produce huge variations 
of the corresponding medial axis.   
This observation has prompted a large body of research that has roughly followed two lines, both aimed 
at the definition of some stable modification of the medial axis: 
one consists of reducing the complexity of the medial axis by pruning the less important parts of the domain
\cite{SB98}, the other considers the definition of filter conditions that identify subsets of 
the medial axis which are stable to perturbations of the sets and retain some of its topological
properties, for instance,  homotopy equivalence with the object. 
Within this second line of research, we mention, among others, the $\lambda-$medial axis 
introduced in \cite{CL05}, the $\theta-$homotopy preserving medial axis introduced initially in \cite{FLM03} and 
subsequently modified in \cite{SFM07}  to ensure the homotopy equivalence, and the power crust method 
 \cite{ACK01b}.  
The $\lambda$-medial axis and the $\theta$-homotopy preserving medial axis are explicitly defined as subsets 
of the medial axis, being collections of those points of the medial axis that meet 
some geometrical criteria. Such criteria are expressed in terms of a bound either on the distance to the boundary of the 
object or on the separation angle $\theta$ (see Definition \ref{Sec3.Def.OpnAng} below), respectively. 
The power crust, in contrast, and also the algorithm discussed in \cite{Dey06,DZ04}, provide a
continuous approximation of the medial axis constructed using a subset of the vertices, called poles, of the Voronoi diagram
of a finite point sample of the object boundary. In all such works, 
the stable modifications are sought  by
identifying directly points of the medial axis or of an approximation of it. 
The excellent survey paper \cite{ABE09}
contains a thorough discussion of such approaches and of  the related stability issues.

We adopt in this paper a fundamentally distinct strategy which, if compared with the  works mentioned above,  
represents an indirect approach relying on the use of the compensated convex transforms \cite{Zha08a}.
The theory of compensated convex transforms
has been introduced and applied  in the calculus of variations
for finding the quasiconvex envelope of a function \cite{Zha99,Zha01a, Zha02,Zha03b}
and for finding tight smooth approximations of the maximum function and the squared-distance 
function  \cite{Zha08b}. Compensated convex transforms, however, 
also provide a natural and stable global 
method to extract geometric singularities, such as ridges, valleys and edges, from a given function 
by manipulating its `landscape' \cite{ZOC14,ZOC14b}, and it is in this way that the transforms, 
in particular the lower compensated convex transform (hereafter, called also the lower transform), 
will be used in this paper. Whether one applies  the lower compensated convex transform 
or the upper compensated convex transform depends on the type of geometric singularities
to be extracted. The works \cite{ZOC14,ZOC14b} present a systematic study on the use 
of these basic transforms to extract singularities from the graph of functions in general,
or from the characteristic functions of compact sets, 
whereas the patent application \cite{ZOC09} contains various applications 
including our method for extracting the multiscale medial axis map.
The key properties 
that are exploited to highlight and/or to design a specific singularity are: the tight 
approximation of the compensated transforms, their regularity and the manner in which they 
respond to the type of curvature. More specifically, 
\cite{ZOC14} focuses on the basic use of these 
transforms to detect ridges, valleys and saddle points of graph of functions,  whereas 
\cite{ZOC14b} presents the design of a transform which is capable of filtering out the `regular points'
and the `regular directions' on manifolds.

The application of the lower transform to study the medial axis of a set is motivated fundamentally by 
the identification of the medial axis with 
the singularity set of the distance function \cite[Lemma 8.5.12]{Hor83} and by the geometric structure 
of this set \cite{ACNS13,CP01,MM03}.
On the other hand, the distance function, its regularity and its geometric structure, are well-studied 
both in geometric measure theory \cite{Fed59} and in the theory of partial differential equations \cite{CSi04,Eva10,GT98,Hor83}.
If the set $K$ is a smooth compact submanifold of $\R^n$, there are many local regularity results of the distance
function near $K$ \cite{KP81,Fo84,DZo94,DZo98}, whereas, 
for a general bounded open set $\Omega$, some results by Albano \cite{Alb14} imply
that the distance function $\dist(\cdot;\, \Omega^c)$ is locally  $C^{1,1}$ in 
$\Omega\setminus\overline{M_{\Omega^c}}$ 
in the sense that if $x\in \Omega\setminus\overline{M_{\Omega^c}}$, there is a $\delta>0$ such that
$\dist(\cdot;\,\Omega^c)\in C^{1,1}(\overline{B}(x;\,\delta))$. 

In the following, however, it is more convenient to refer to the squared-distance function 
and use the identification of the singular set of the distance function with the set of points
where the squared distance function fails to be locally $C^{1,1}$. Here, we just note that 
the advantage of referring to the squared distance function rather than to the distance 
function has also been realized in other contexts, such as, in the study of the motion of surfaces by its mean 
curvature represented by manifolds with codimension greater 
than one \cite{Deg96,AS96}. We refer to \cite{AM98} for a detailed study on the properties of the squared distance function 
and on its applications in the geometric evolution problems. 

Using properties of the lower transform, we 
apply the lower compensated convex transform to the Euclidean squared-distance function which gives 
a smooth $(C^{1,1}$) tight approximation outside a neighbourhood of the closure of the medial axis 
(see Theorem \ref{Sec3.Thm.TghtApprx.Eq1}), and define our multiscale medial axis map as a scaled difference 
between the squared-distance function and its lower transform. 
From the property of the tight approximation of the lower transform of 
the squared-distance function, we also deduce a sharp
$C^{1,1}$-regularity result (see Corollary \ref{Sec3.Cor.Reg} and Example \ref{Sec3.Ex.ShrpEst})  
of the squared-distance function outside a neighbourhood of the closure of the medial axis of $K$, 
which can be viewed as a weak Lusin type theorem for the squared-distance function and extend regularity results
of the squared-distance function to any closed non-empty subset of $\R^n$. This result also 
offers an instance of application of the compensated convex transform to obtain a fine result of geometric measure theory
and is, somehow,  related to the behaviour of semiconcave functions 
(see \cite{CSi04} and Remark \ref{Sec2.Rem.Sem}$(c)$ below). 
We observe that, in general, the regularity of 
$\dist^2(\cdot;\, K)$ cannot be better than $C^{1,1}$ even for a compact convex set $K$ where $M_K=\varnothing$. 
A simple example is given by the square $K=[0,1]\times [0,1]\subset \mathbb{R}^2$. 
In this case, it can be easily verified that $\dist^2(x;\, K)$  is 
globally $C^{1,1}$ but not $C^2$.

The application of the lower compensated convex transform of scale $\lambda$ to the squared-distance function 
produces a continuous function in $\mathbb{R}^n$ 
that remains strictly positive on the medial axis and tends to zero outside of it as a positive parameter $\lambda$
becomes very large (see Proposition \ref{Sec3.LimInfSup}). We will, in fact, characterize the limit 
of the multiscale medial axis map of 
scale $\lambda$ as $\lambda$ approaches to infinity (see Theorem \ref{Sec3.Teo.LimMMA}) and refer to this geometric structure
as the quadratic multiscale medial axis landscape map of $K$. The values of this map are well
separated, in the sense that they are zero outside the medial axis and remain strictly positive on it. Furthermore, 
we will give conditions (see Proposition \ref{Sec3.Pro.OpAng} and Section \ref{SecEx})
that ensure that the multiscale medial axis map of 
scale $\lambda$ actually keeps a constant height
along the parts of the medial axis generated by two-point subsets, with the value of the height dependent on
the distance between the two generating points. Such values can, therefore, be used to define a 
hierarchy between different parts of the medial axis and we can thus select the relevant parts through 
simple thresholding, that is, by taking suplevel sets of the multiscale medial axis map. 
To reflect this property, we use the word "multiscale". 
For each branch of the medial axis, the multiscale medial axis map automatically defines a scale associated with it. In other words, 
a given branch has a strength which depends on some geometric features of the part of the set that generates that branch.

Given a closed non-empty subset $K$ of $\R^n$, we will also prove that, despite the medial axis of $K$ not being Hausdorff stable, 
the quadratic multiscale medial axis map, is indeed
Hausdorff stable (see Theorem \ref{Sec4.Teo.HAS}). It follows  that
the graph of the medial
axis map carries more information than the medial axis itself, which allows the definition of a 
hierarchy between the parts of the medial axis and the selection of the relevant ones through 
simple thresholding, that is, by 
taking suplevel sets of the medial axis map. In this manner, it is possible to choose the main parts 
that reflect genuine geometric features of the object and remove minor ones 
generated by noise. 

In conclusion, we  observe that
while our method seems to share similarities with those based 
on the extraction of ridges of the distance transform \cite{AS92,BRS91,KSKB95,Mon68,WCG95}, (that require, 
however, an a-priori definition of ridge,
based usually on an approximation of the derivative of the distance transform),  
the method we propose is, in fact, substantially different from such approaches, given that we obtain a neighbourhood of 
the singularities as the difference between the squared-distance transform and its smooth tight approximation. In this manner, 
as mentioned above, we
provide an indirect definition of the singularity, which does not require any derivative 
approximation or any differentiability assumption.

After this brief introduction, the next section will introduce the relevant notation and recall basic results
in convex analysis and  lower compensated convex transforms. 
Section \ref{SecMR} contains the definition of the multiscale medial axis map, and some of 
its principal properties, such as the tight approximation of the lower compensated transform 
to the squared-distance transform (see Theorem \ref{Sec3.Thm.TghtApprx.Eq1}) and as an application, we deduce a
sharp regularity result of the squared-distance function to any non-empty closed subset of $\R^n$ (see Corollary \ref{Sec3.Cor.Reg}). 
Section \ref{SecHAU} presents the Hausdorff stability of the multiscale medial axis map whereas 
Section \ref{SecEx} discusses some mathematical prototype models of explicitly calculated medial axis maps 
for a simple four point set to some more complicated three dimensional objects. Finally, 
Section \ref{SecPrf} concludes the paper with the proofs of the main results. 

%%%%%%%%%%%%%%%%%%%%%%%%%%%%%%%%%%%%%%%%%%%%%%%%%%%%%%%%%%%%%%%%%%%%%%%%%%%%%%%%%%%%%%%%%%%%%%%%%%%%%
%%%%%%%%%%%%%%%%%%%%%%%%%%%%%%%%%%%%%%%%%%%%%%%%%%%%%%%%%%%%%%%%%%%%%%%%%%%%%%%%%%%%%%%%%%%%%%%%%%%%%
%%%%%%%%%%%%%%%%%%%%%%%%%%%%%%%%%%%%%%%%%%%%%%%%%%%%%%%%%%%%%%%%%%%%%%%%%%%%%%%%%%%%%%%%%%%%%%%%%%%%%

\section{Notation, Basic Definitions and Preliminary Results}\label{SecNot}
Throughout the paper $\mathbb{R}^n$ denotes the $n$-dimensional Euclidean space, and $|x|$ and 
$x\cdot y$ the standard Euclidean norm and inner product, respectively, for  $x,\, y\in \mathbb{R}^n$. 
In some cases, we will also make use of the  notation 
$(x,y)$ to denote the point of $\R^n$ given by $xe_1+y_2e_2+\ldots+y_ne_n$, 
where $\{e_1,\ldots,e_n\}$ is an orthornormal basis of $\R^n$, 
$(x,\,y_2,\ldots,y_n)\in \R^n$ and $y=y_2e_2+\ldots+y_ne_n$.
Given a non-empty subset $K$ of $\R^n$,
$K^c$  denotes the complement of $K$ in $\R^n$, i.e. $K^c=\R^n\setminus K$,  $\overline{K}$  its closure
and $\co[K]$  the convex hull of $K$, that is, the smallest (with respect to inclusion) 
convex set that contains the set $K$. 
For $x\in \R^n$
and $r>0$, $B(x;\,r)$ indicates the open ball with center $x$ and radius $r$ whereas $S(x;\,r)$ denotes the 
sphere with center $x$ and radius $r$ and is the boundary of $B(x;\,r)$.
The distance transform of  a non-empty set $K\subset \mathbb{R}^n$ 
is the function that, at any point $x\in \mathbb{R}^n$, associates the distance of 
$x$ to $K$, which is defined as $\inf\{|x-y|,\; y\in K\}$ and is
denoted as $\dist(x;\,K)$. We use the notation $Df$ to denote the derivative of $f$.

Across the current literature, there is no uniform definition of the medial axis,
with its meaning changing from one author to another. What the medial axis is for one, 
becomes the skeleton for another, and in some cases subtle 
differences are present, especially in the continuum case, where 
the closure of such sets is considered. In this paper we adopt the definition given 
by Lieutier in \cite{Lie04}, but it is reformulated here to include a non-empty closed set 
$K\subset \R^n$ with $K\neq \R^n$ as well as a non-empty bounded open set $\Omega$.

\begin{defi} \label{Def.MA}
	For a given non-empty closed set $K\subset \mathbb{R}^n$, with $K \neq \R^n$, we define the 
	{\bf medial axis} $M_K$ of $K$ as  
	 the set of points  $x\in\mathbb{R}^n\setminus K$ such that
	$x\in M_K$ if  and only if there are at least two different points $y_1,\, y_2\in K$,
	satisfying $\dist(x;\,K)=|x-y_1|=|x-y_2|$.
	For a non-empty bounded open set $\Omega\subset \mathbb{R}^n$, the medial axis  of $\Omega$
	is defined by $M_\Omega:=\Omega\cap M_{\partial\Omega}$. 
\end{defi}

%%%%%%%%%%%%%%%%%%
\medskip
%%%%%%%%%%%%%%%%%%

\begin{nota}
\begin{itemize}
	\item[$(a)$] The definition of medial axis for a bounded open set $\Omega$ is equivalent to that of the closed set
		$\Omega^c$, since the definition of $M_K$ implies that $M_{\Omega^c} = M_{\partial \Omega}$, and hence $M_{\Omega} = M_{\Omega^c}$. 
	\item[$(b)$] 
	Other frequently used notions are those of the skeleton of $K$, denoted as $skl(K)$, and 
		the cut locus of a manifold,  denoted as $cl(K)$, which applies to 
		the more general case of Riemannian geometry. 
		Here, for a non-empty closed subset $K$ of the Euclidean space $\R^n$, we define the skeleton of $K$ to be the set 
		of the centers of the maximal (with respect to inclusion) open balls contained in $K^c$,
		whereas the cut locus of $K$ is taken to be the set of the cut locus of the points of the boundary of $K$
		in $\R^n$, where the cut locus of $p$ in $\partial K$ is 
		the set of points in the manifold where the geodesics starting at $p$ stop being minimizing. 
		It can then be shown that $ cl(K)=\overline{M_K}$.
		As a result, the notions of medial axis, skeleton and cut locus 
		are related but are not the same;  for example, in general,   \cite{Mat88,Dey06}
	\begin{equation}\label{Eq.SkMA}
		M_{K}\subset skl(K)\subset  \overline{M_K}\,.
	\end{equation}
	\item[$(c)$] Our definition of the medial axis $M_K$ is, in particular, consistent with our main 
		convergence result (Theorem \ref{Sec3.Teo.LimMMA})  which recovers $M_K$ as the set on which the limit of the 
		medial axis map $M_{\lambda}(\cdot;\,K)$ (Definition \ref{Sec3.Def.MMAM}) as $\lambda$ tends to infinity is
		strictly positive.
\end{itemize}
\end{nota}

%%%%%%%%%%%%%%%%%%
\medskip
%%%%%%%%%%%%%%%%%%

Next we collect definitions and results from convex analysis for functions $f$ taking finite values, i.e. for
$f:\mathbb{R}^n\mapsto \mathbb{R}$, which will be used in this paper, and refer to 
\cite{HUL01,Roc66} for details and proofs.\\

Given a function $f:\mathbb{R}^n\mapsto \mathbb{R}$ bounded below, the convex envelope $\co[f]$ is 
the largest convex function not greater than $f$. We will often make use of the following 
characterization.

\begin{prop}\label{Prp.GlCnvx}
	Let $f:\mathbb{R}^n\mapsto \mathbb{R}$ be coercive in the sense that
	$f(x)/|x|\to\infty$ as $|x|\to \infty$, and $x_0\in\R^n$. Then
   \begin{itemize}
	\item[(i)]
		The value $\co\left[f\right](x_0)$ of the convex envelope of $f$ at
		$x_0\in \mathbb{R}^n$ is given by
		\begin{equation}\label{Eq.Prp.GlCnvx.1}
			\co\left[f\right](x_0)=\underset{i=1,\ldots,n+1}{\inf}\,\Bigg\{ \sum^{n+1}_{i=1}\lambda_if(x_i):
				\; \sum^{n+1}_{i=1}\lambda_i=1,\; \sum^{n+1}_{i=1}\lambda_ix_i=x_0,\;
				\lambda_i\geq 0,\; x_i\in \R^n
			\Bigg\}\,.
		\end{equation}
		If, in addition, $f$ is lower semicontinuous, the infimum is reached by some
		$(\lambda_i^\ast,\, x_i^\ast)$ for $i=1,2,\dots,n+1$ with
		$\left(x_i^\ast, f(x_i^\ast)\right)$'s lying in the intersection
		of a supporting plane of
		the epigraph of $f$, $\epi(f)$, and
		$\epi(f)$.
	\item[(ii)]
		The value $\co\left[f\right](x_0)$, for $f$ taking only finite values, can also be obtained as follows:
		\begin{equation}\label{Eq.Prp.GlCnvx.2}
			\co\left[f\right](x_0)=\sup\left\{
				\ell(x_0):\; \ell\;\;\text{ affine}\quad\text{ and}\quad
				\ell(y)\leq f(y)\;\;\text{ for all }y\in\R^n
				\right \}
		\end{equation}
	with the $\sup$ attained by an affine function $\ell^\ast\in\Aff(\R^n)$.
     \end{itemize}
\end{prop}
 
%%%%%%%%%%%%%%%%%%
\medskip
%%%%%%%%%%%%%%%%%%

We then introduce the following definition \cite{BKK00}, which is needed in Lemma \ref{Sec2.Lem.UpDif}.

\begin{defi} \label{Def.UpSemiDif}
Assume $x_0\in \R^n$. We say that a function $f:\R^n\mapsto \R$ is  upper-semidifferentiable at $x_0$ if there exists $a\in\R^n$ such that
\[
	\limsup_{y\to 0}\frac{f(x+y)-f(x)-a\cdot y}{|y|}\leq 0\,.
\]
\end{defi}

%%%%%%%%%%%%%%%%%%
\medskip
%%%%%%%%%%%%%%%%%%

The following lemma, concerning the existence and properties of an optimal affine function, will be needed  for the proofs of 
Proposition \ref{Sec2.Th.LipGradLwTr} and Theorem \ref{Sec3.Thm.TghtApprx}. 

\begin{lema}\label{Sec2.Lem.UpDif} 
	Suppose $f:\R^n\mapsto \R$ is continuous, upper-semidifferentiable, coercive in the sense that 
	$\lim_{|x|\to\infty}f(x)/|x|=+\infty$, and $\co[f]\in C^{1,1}(\R^n)$. 
	If $\co[f](0)<f(0)$, then
	there is an affine function $\ell(x)=2a\cdot x+b$ and distinct points 
	$x_1,x_2,\dots,x_k\in\R^n$ and $\lambda_1>0,\lambda_2>0,\dots,\lambda_k>0$ 
	satisfying $2\leq k\leq n+1$, $\sum^k_{i=1}\lambda_i=1$, $\sum^k_{i=1}\lambda_ix_i=0$ and
	$x_i\neq 0$, $x_i\neq x_j$ if $1\leq i\neq j\leq k,$ such that
\begin{equation}\label{Sec2.Lem.UpDif.Eq} 
	\begin{array}{ll}
		(i)   &\displaystyle \ell(x)\leq f(x)	\text{ for all }x\in \R^n\,;\\
		(ii)  &\displaystyle \ell(x_i)=f(x_i)	\text{ for }i=1,2,\dots, k\,;\\
		(iii) &\displaystyle 2a=D\ell(x_i)=Df(x_i) \text{ for }i=1,2,\dots, k;\\
		(iv)  &\displaystyle  b=\ell(0)=\co[f](0)\,;\\
		(v)  &\displaystyle 2a=D\ell(0)=D\co[f](0)\,.
	\end{array}
\end{equation}
\end{lema}

%%%%%%%%%%%%%%%%%%
\medskip
%%%%%%%%%%%%%%%%%%

The quadratic lower compensated convex transform, introduced in \cite{Zha08a}, will play a pivotal role in the definition
of our multiscale medial axis map. We next recall  its definition and some of its properties, and refer to 
\cite{Zha08a,ZOC14, ZOC14b} for details and proofs.

\begin{defi}\label{Sec2.Def.LwTr}
Let $f:\mathbb{R}^n\mapsto \mathbb{R}$ be a lower semicontinuous function \cite{Roc66,HUL01} satisfying
\begin{equation} \label{Eq.Sec2.Bndf}
	f(x)\geq -A_1|x|^2-A_2
\end{equation}
for some constants $A_1,\, A_2\geq 0$. 
The (quadratic) lower compensated convex transform (lower transform for short) for $f$ with scale $\lambda>A_1$ 
is defined for $x\in \R^n$ by
\begin{equation}\label{Eq.Sec2.LwTr}
	C^l_\lambda(f)(x)=\co[f+\lambda|\cdot|^2](x) -\lambda|x|^2
\end{equation}
If $f$ is bounded below, we may set $\lambda\geq 0$.
\end{defi}

%%%%%%%%%%%%%%%%%%
\medskip
%%%%%%%%%%%%%%%%%%

\begin{nota}\label{Sec2.Rem.Sem}
\begin{itemize}
	\item[$(a)$]	The requirement of the lower semicontinuity of $f$ is to guarantee  that 
		$C^l_\lambda(f)(x)\to f(x)$ as $\lambda\to\infty$ for all $x\in \mathbb{R}^n$, since
		otherwise, the lower transform will converge to the lower semicontinuous
		envelope of $f$. 
	\item[$(b)$]	From \eqref{Eq.Sec2.LwTr} it also follows that
		$C^l_\lambda(f)(x)$ is the envelope of all the quadratic functions with fixed quadratic term $\lambda |x|^2$ that
		are less than or equal to $f$, that is, for $x\in \R^n$ \cite[Eq. (1.4)]{ZOC14}
		\begin{equation} \label{LwTrEnvPrb}
			C^l_\lambda(f)(x)=\sup\left\{
				-\lambda|x|^2+\ell(x):\; -\lambda|y|^2+\ell(y)\leq f(y)\;\;\text{\rm for all }y\in\R^n
				\;\;\text{\rm and }\ell\;\;\text{\rm affine}\right \}\,.
		\end{equation}
	\item[$(c)$] Recalling from \cite{CSi04} that a function $f:\mathbb{R}^n\mapsto \mathbb{R}$ is called $c-$semiconvex if, for 
		some constant $c>0$, the function $f(x)+c/2|x|^2$ is convex, we observe that 
		the lower compensated convex transform for $f$ with scale $\lambda$, $C^l_\lambda(f)$ is a $2\lambda$-semiconvex
		function. In fact, $C^l_\lambda(f)$ represents the $2\lambda$-semiconvex envelope of $f$. We sometimes use such a property
		to extend some properties of semiconvex functions to $C^l_\lambda(f)$. 
	\item[$(d)$] To gain further geometric insight {into} the lower compensated convex transform defined by \eqref{Eq.Sec2.LwTr}, 
		in Figure \ref{Sec2.FigCnstLwTr}  
		we display the steps of the construction of 
		$C_{\lambda}^l(f)$ for $f(x)=\dist^2(x;\,K)$ with $K=\{-1,\,1\}$ and $\lambda=2$. 
		The graph of the augmented function $f+\lambda |x|^2$ is displayed in Figure \ref{Sec2.FigCnstLwTr}$(b)$ along with $f$, 
		whereas  Figure \ref{Sec2.FigCnstLwTr}$(c)$ shows the convex envelope of the augmented function. Figure \ref{Sec2.FigCnstLwTr}$(d)$ 
		displays finally the graph of $C_{\lambda}^l(f)$ which is compared with that of $f$. Note that 
		the convex envelope of the augmented function is different from $f+\lambda |x|^2$ only in a 
		neighbourhood of the singular point $0$ of $f$, so that, when then we subtract the weight, 
		the final effect is a smoothing of $f$ only in such neighbourhood.
		This simple example, along with the ones discussed in Section \ref{SecEx}, 
		enables one also to understand the role of the parameter $\lambda$ and our meaning of scale.
		The parameter $\lambda$ acts as a scale parameter in the sense that it controls the curvature of the 
		lower compensated convex transform in the neighbourhood of the singularity of the function and allows the extraction
		of the singularity with a value which gives somehow a measure of its strength. Also one may observe the so-called 'tightness' of
		the lower compensated convex transform approximation of the original function from below (see Proposition \ref{Sec2.Pro.Tight}), 
		which agrees with the original function
		except in the neighbourhood near the singular point.
		\begin{figure}[htbp]
			\centerline{
			$\begin{array}{cc}
				\includegraphics[height=0.35\textwidth]{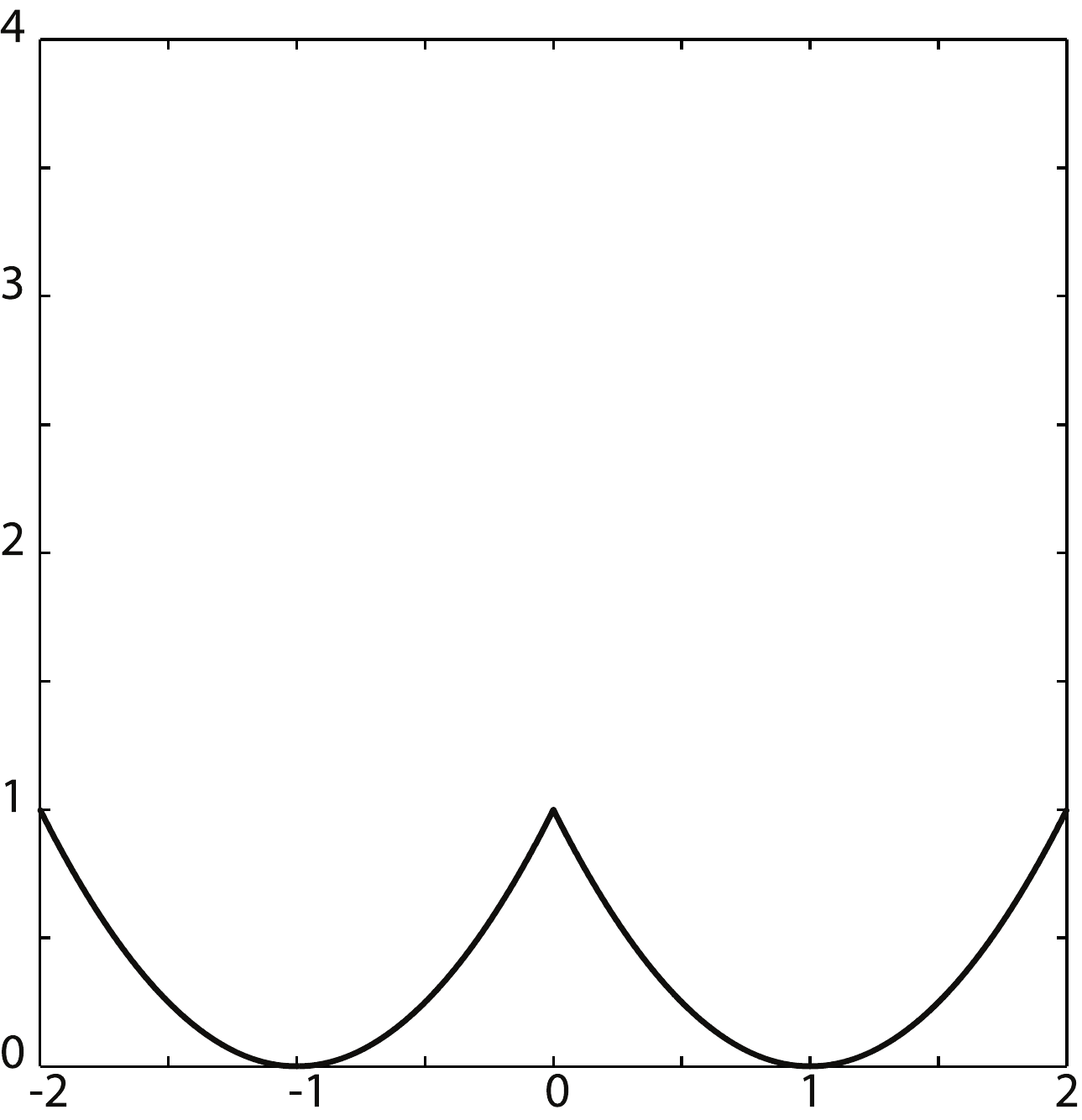}&
				\includegraphics[height=0.35\textwidth]{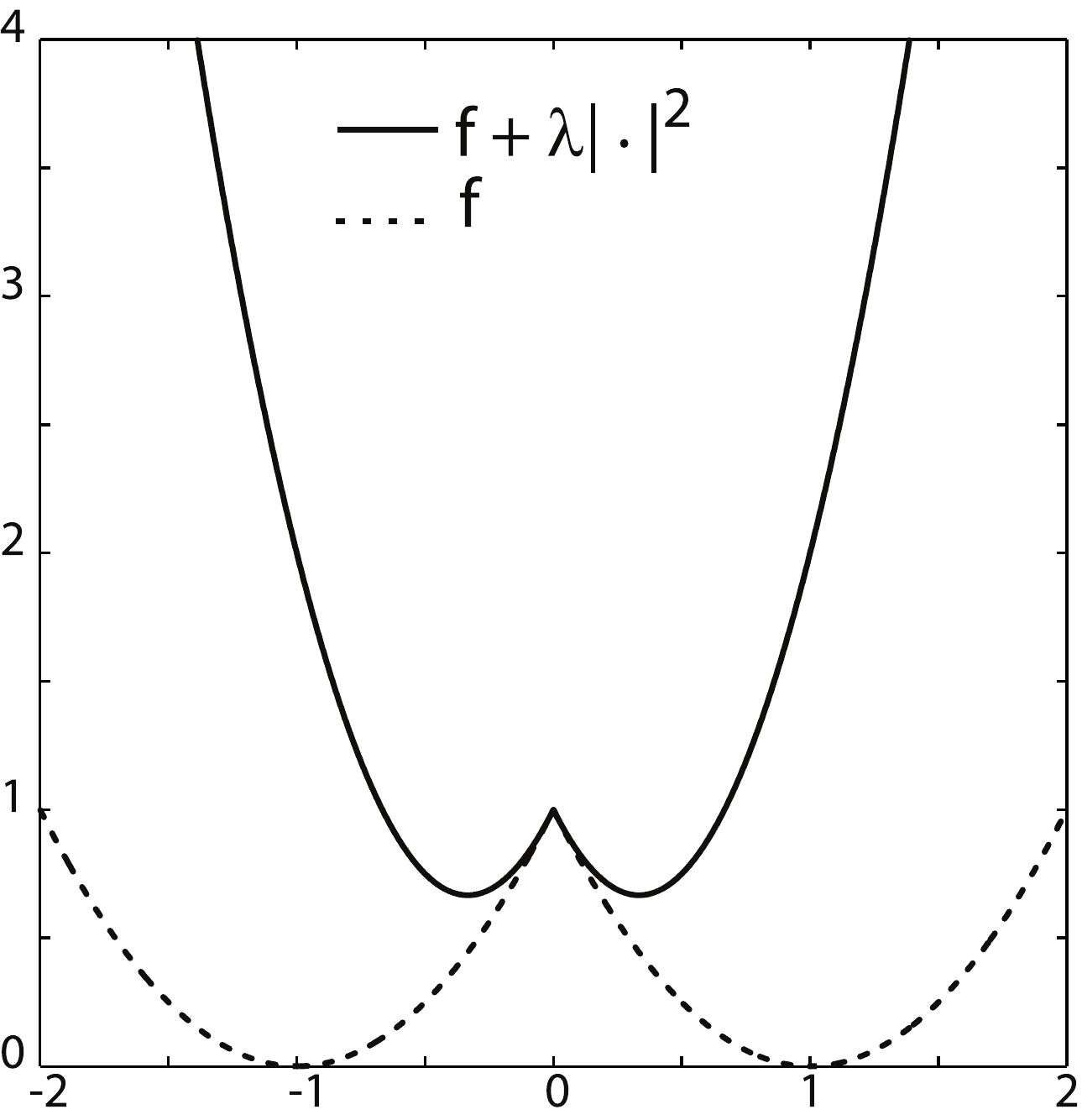}\\
				(a) & (b)\\
				\includegraphics[height=0.35\textwidth]{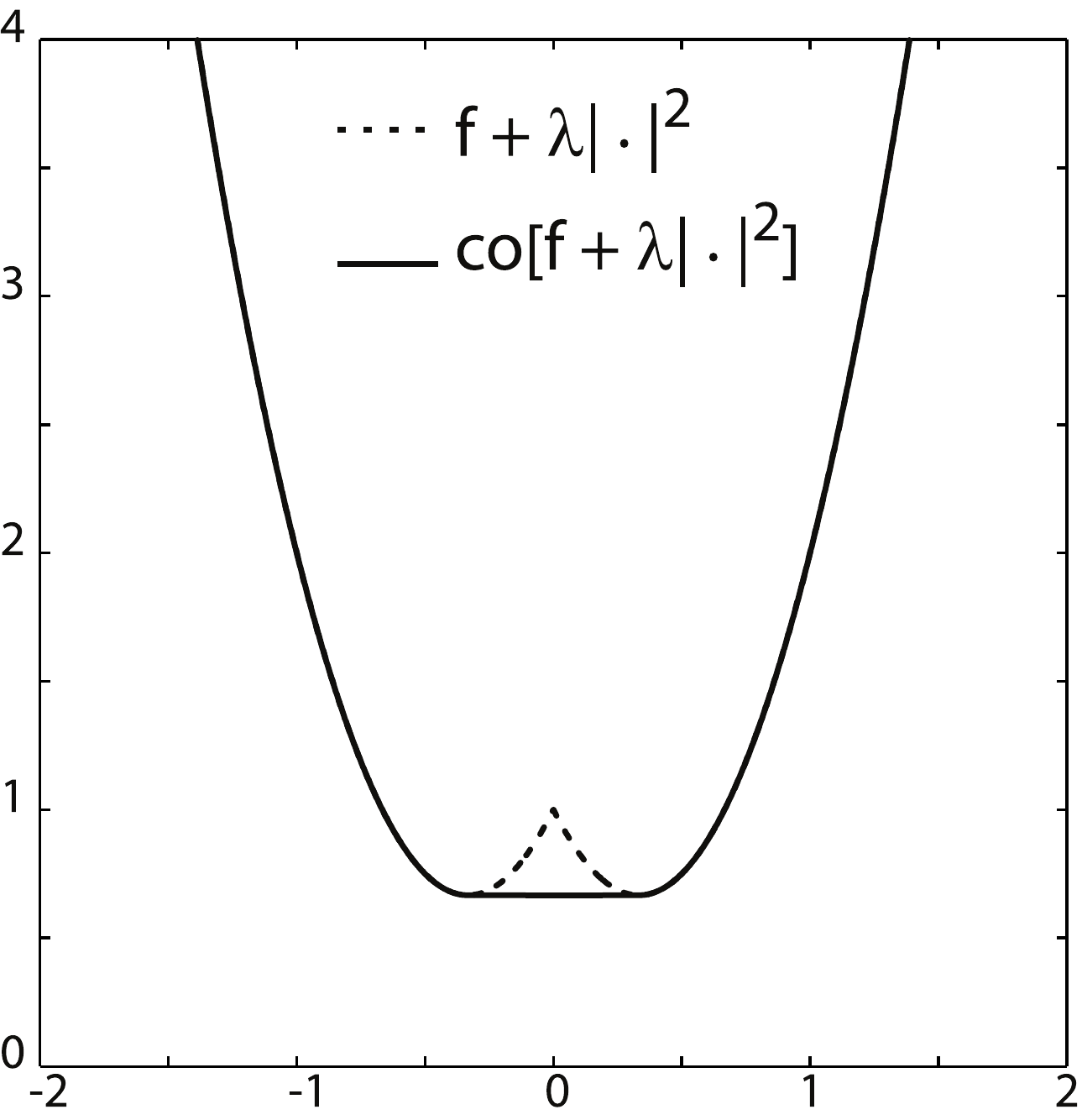}&
				\includegraphics[height=0.35\textwidth]{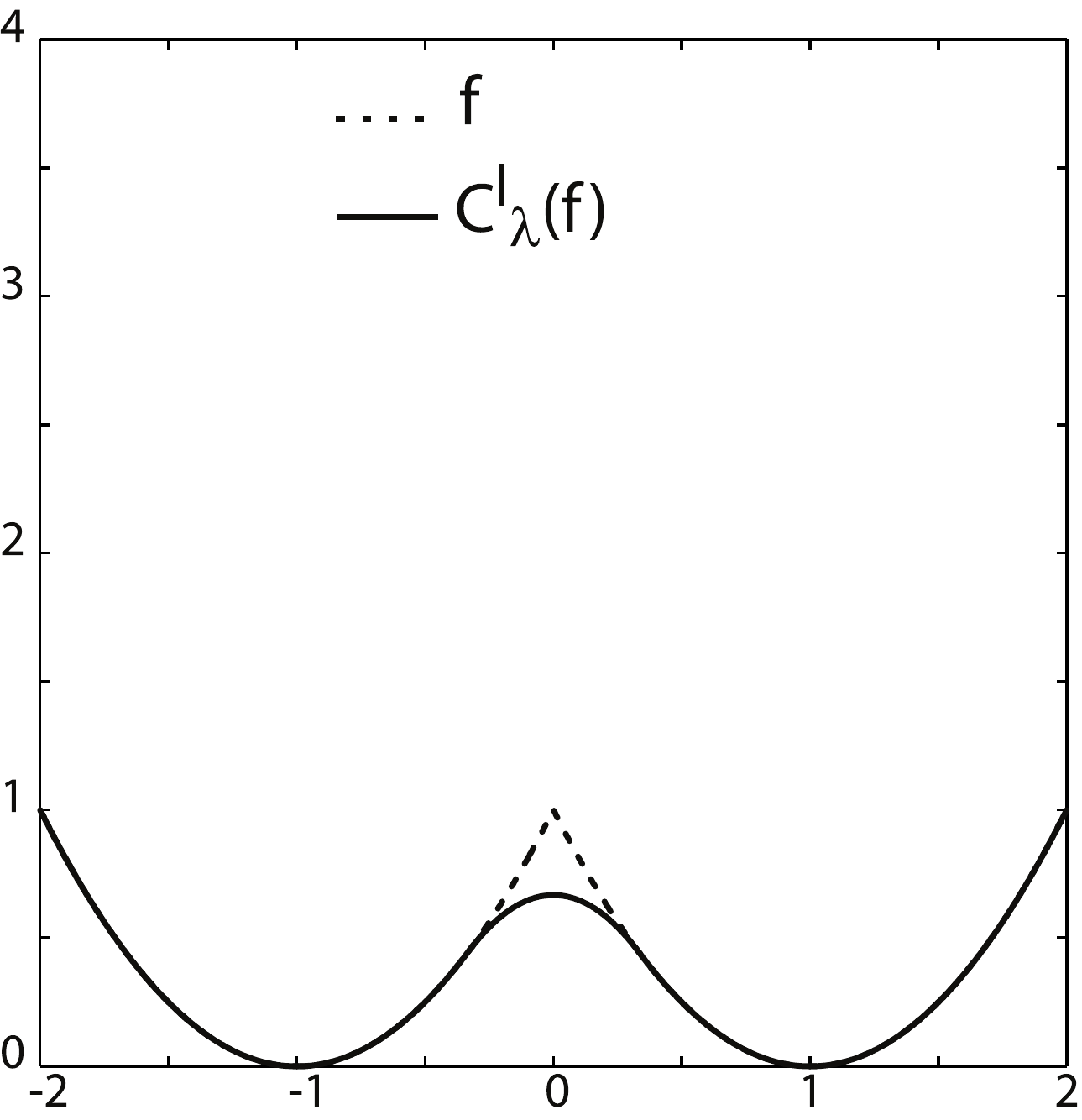}\\
				(c) & (d)
			\end{array}$
			}
			\caption{\label{Sec2.FigCnstLwTr} 
				Steps illustrating the construction of the lower compensated convex transform 
				of $f(x)=\dist^2(x;\,K)$ with $K={-1,\,1}$. 
				$(a)$ Graph of the function $f$;  
				$(b)$ Graph of the augmented function $f+\lambda|\cdot|^2$ with $\lambda=2$; 
				$(c)$ Graph of the convex envelope of $f+\lambda |\cdot|^2$ compared to that of $f+\lambda|\cdot|^2$; 
				$(d)$ Graph of $C_{\lambda}^l(f)$ compared to that of $f$.
				}
			\end{figure}
\end{itemize}
\end{nota}

%%%%%%%%%%%%%%%%%%
\medskip
%%%%%%%%%%%%%%%%%%

The following properties of $C^l_\lambda(f)$ will also be used. 

\begin{prop}\label{Sec2.Pro.MonLwTr}
Given $f:\mathbb{R}^n\mapsto \mathbb{R}$ that satisfies \eqref{Eq.Sec2.Bndf}, then for all $A_1<\lambda<\tau<\infty$, 
we have
	\begin{equation}
		C^l_\lambda(f)(x)\leq C^l_\tau(f)(x)\qquad\text{for }x\in\R^n\,,
	\end{equation}
and, for $\lambda>A_1$
	\begin{equation}
		C^l_\lambda(f)(x)\leq f(x)\qquad\text{for }x\in\R^n\,.
	\end{equation}
\end{prop}

\begin{prop}\label{Sec2.Pro.OrdLwTr}
If $f\leq g$ in $\R^n$ and satisfy \eqref{Eq.Sec2.Bndf}, then 
\begin{equation}
	C^l_\lambda(f)(x)\leq C^l_{\lambda}(g)(x)\qquad\text{for }x\in\R^n\text{ and }
	\lambda\geq \max\{A_{1,f}, \,A_{1,g}\}\,.
\end{equation}
\end{prop}

%%%%%%%%%%%%%%%%%%
\medskip
%%%%%%%%%%%%%%%%%%

The transform  $C^l_\lambda(f)$ realizes a `tight' approximation
of the function $f$, in the following sense (see \cite[Theorem 2.3$(iv)$]{Zha08a}).

\begin{prop}\label{Sec2.Pro.Tight}
	Let $f\in C^{1,1}(\overline B(x_0;\,r))$, with $B(x_0;\,r)$ the open ball of center $x_0$ and radius $r>0$.
	Then for sufficiently large $\lambda>0$, we have that $f(x_0)=C^l_\lambda(f)(x_0)$.
\end{prop}

%%%%%%%%%%%%%%%%%%
\medskip
%%%%%%%%%%%%%%%%%%

Such a property motivates  the definition of the multiscale ridge transform which was introduced 
in \cite{ZOC14} to extract ridges of general functions and  shown to be invariant with 
respect to translation. 
This multiscale ridge transform will  be used in Section \ref{SecMR} to define the multiscale medial
axis map (see Definition \ref{Sec3.Def.MMAM}).  

\begin{defi}
Given $\lambda >0$, the ridge transform of scale $\lambda$, for a given function $f:\mathbb{R}^n\mapsto \mathbb{R}$ satisfying 
\eqref{Eq.Sec2.Bndf}, is defined as:
\begin{equation}\label{Eq.Sec2.Rdg} 
	R_\lambda(f(x)):=f(x)-C^l_\lambda(f)(x),\quad x\in\mathbb{R}^n\,.
\end{equation}
\end{defi}

%%%%%%%%%%%%%%%%%%
\medskip
%%%%%%%%%%%%%%%%%%

We now present some regularity properties of $C^l_\lambda(\dist^2(\cdot;\, K))$, which will be exploited to analyze 
the behaviour of the multiscale medial axis map. We recall first the following result given in
\cite[Lemma 4.3]{Zha08a}.

\begin{lema}\label{Sec2.Lem.KZLem43} 
	Suppose $f:\R^n\to \R$ is convex and such that $|f(x)|\leq c|x|^2+c_1$ for $x\in \R^n$ with $c,\,c_1>0$. 
	Assume $\lambda>0$ and define $f_{\lambda}=\lambda|\cdot|^2-f$. Then for $\lambda >c$, 
	\begin{equation}\label{Sec2.Eq.KZLem43}
		\co[f_{\lambda}](x+y)- \co[f_{\lambda}](x)-D\co[f_{\lambda}](x)\cdot y \leq \lambda |y|^2
	\end{equation}
for $x,y\in \R^n$.
\end{lema}

The next proposition improves a result in 
\cite[Theorem 3.1]{Zha08a}.

\begin{prop}\label{Sec2.Th.RegLwTr} 
	Suppose $K\subset \mathbb{R}^n$ is a non-empty closed set. 
	Then for $\lambda>0$, 
	\[
		C^l_\lambda(\dist^2(\cdot;\, K))\in C^{1,1}(\mathbb{R}^n)\,.
	\]	
	Furthermore, the Lipschitz constant of the gradient 
	$DC^l_\lambda(\dist^2(\cdot,\, K))$ is at most $2\max\{1,\,\lambda\}$.
\end{prop}

%%%%%%%%%%%%%%%%%%
\medskip
%%%%%%%%%%%%%%%%%%

The next property is a useful inequality for the derivative of the lower transform,  $DC^l_\lambda(\dist^2(\cdot;\, K))$.

\begin{prop} \label{Sec2.Th.LipGradLwTr}
Suppose $K\subset \mathbb{R}^n$ is a non-empty closed set. Then
\begin{equation} 
\label{Sec2.ineq.deriv}
	|DC^l_\lambda(\dist^2(\cdot;\, K))(x)|^2\leq 4C^l_\lambda(\dist^2(\cdot;\, K))(x)
\end{equation}
and  equality holds in (\ref{Sec2.ineq.deriv}) if and only if $C^l_\lambda(\dist^2(\cdot;\, K))(x)=\dist^2(x;\, K)$.
\end{prop}

%%%%%%%%%%%%%%%%%%
\medskip
%%%%%%%%%%%%%%%%%%

\begin{nota}
	Given that $u=\dist^2(\cdot;\, K)$ is known to be a viscosity solution of the Hamilton-Jacobi equation \cite{MM03,CSi04} 
\begin{equation}\label{Sec2.Eq.HJ}
	|Du|^2=4u\,,
\end{equation}
	 Proposition \ref{Sec2.Th.RegLwTr} and Proposition \ref{Sec2.Th.LipGradLwTr} together
	imply that the lower transform $C^l_\lambda(\dist^2(\cdot;\, K))$ is a $C^{1,1}$ 
	subsolution of the Hamilton-Jacobi equation \eqref{Sec2.Eq.HJ}.
\end{nota}

%%%%%%%%%%%%%%%%%%
\medskip
%%%%%%%%%%%%%%%%%%

We next introduce  the sets $K(x)$ and $K_{2,\lambda}(x)$, which will be used to gain insight into  the geometric 
structure of $C^l_\lambda(\dist^2(\cdot;\, K))$.

\begin{defi}\label{Def.SetK}
Let $K\subset \R^n$ be a non-empty closed set. 
For any $x\in \mathbb{R}^n$, let $r(x)=\dist(x;\,K)$. 
We then define the following sets:
\begin{equation}\label{Eq.Def.SetK}
	K(x)=\Big\{x+y\in \R^n:\, x+y\in\partial K\text{ and }|y|=r(x)\Big\}
\end{equation}
and for $\lambda>0$,
\begin{equation}\label{Eq.Def.SetK2}
	K_{2,\lambda}(x):=\Big\{x+\frac{y}{1+\lambda}\in \R^n:\,x+y\in \partial K\text{ and }|y|=r(x)\Big\}\,.
\end{equation}
\end{defi}

%%%%%%%%%%%%%%%%%%
\medskip
%%%%%%%%%%%%%%%%%%

\begin{nota}
\label{Rk.cpt}
	If $x \not\in K$, the set $K(x)$ is the set of points  of $\partial K$ that realize the distance of $x$ to $K$. 
	Note also that it follows from   
	\eqref{Eq.Def.SetK} that $K(x)=\overline B(x;\,r(x))\cap K$, so in particular, $K(x)$ is compact, and if $x \not\in K$, 
	$B(x;\,r(x))\subset K^c$.
\end{nota}

%%%%%%%%%%%%%%%%%%
\medskip
%%%%%%%%%%%%%%%%%%

The following result, obtained in the proof  of \cite[Theorem 3.7]{Zha08a}, gives an explicit expression of the 
lower transform of $\dist^2(\cdot ;\,K(x))$, the squared distance to the set $K(x)$, and
will be used to produce a bound on the multiscale medial axis map (see Theorem \ref{Sec3.Teo.BndMMA}$(i)$). 
\begin{prop}\label{Sec2.Pro.LwTrSetK}
Let $K\subset \R^n$ be a non-empty closed set and $M_K$ the medial axis of $K$. 
Assume $x\in M_K$ and denote by $K(x)$ and $K_{2,\lambda}(x)$ the sets defined by 
\eqref{Eq.Def.SetK} and by \eqref{Eq.Def.SetK2}, respectively. Then, for all $y\in\R^n$, 
\begin{equation}\label{Eq.LwTrSetK}
	C^l_\lambda(\dist^2(\cdot;\,K(x)))(y)=(1+\lambda)\dist^2(y;\,\co[K_{2,\lambda}(x)])+
	\frac{\lambda}{1+\lambda}r^2(x)-\lambda|y-x|^2
\end{equation}
where $\co[K_{2,\lambda}(x)]$ is the convex hull of $K_{2,\lambda}(x)$.
\end{prop}

We will also need, for the proof of  Theorem \ref{Sec3.Teo.LimMMA}, 
the following  explicitly calculated formula of the lower transform for compact sets contained in a 
sphere $S(0;\,r)=\{x\in\R^n,\; |x|=r\}$ centred at $0\in\R^n$ with radius $r>0$. 
The formula is easy to derive following similar calculations to those in the proof of  \cite[Theorem 1]{Zha08b},
or of  \cite[Theorem 5.1]{Zha08a}. 

\begin{prop}\label{Sec2.formula-set-on-sphere}
Let $K\subset S(0;\,r)$ be a non-empty compact set. Then for every $x\in\R^n$,
\begin{equation}
	C^l_\lambda(\dist^2(\cdot;\,K))(x)=\frac{\lambda r^2}{1+\lambda}+(1+\lambda)\dist^2
	\left(x;\, \co\left[\frac{K}{1+\lambda}\right]\right)-
\lambda|x|^2,
\end{equation}
where $\co[K/(1+\lambda)]$ is the convex hull of $K/(1+\lambda)=\{x/(1+\lambda),\; x\in K\}$.
\end{prop}

%%%%%%%%%%%%%%%%%%
\medskip
%%%%%%%%%%%%%%%%%%

We will  invoke  the following technical lemma several times (see  Lemma 3.2 in \cite{Zha08a}).
\begin{lema}\label{Sec2.Lem.LwDst} 
	Assume $\rho>0$. Let $K=B^c(0;\,\rho)=\{x\in\mathbb{ R}^n,\quad |x|\geq \rho\}$ be the complement of the open
	ball  $B(0;\,\rho)$ with center the point $0$ and radius $\rho$, then
\begin{equation} 
	C^l_{\lambda}(\dist^2(\cdot;\,K))(x)=\left\{\aligned
	&\frac{\lambda}{1+\lambda}\rho^2-\lambda |x|^2,\quad |x|\leq \frac{\rho}{1+\lambda},\\
	&\dist^2(x;\,K),\quad\qquad |x|\geq \frac{\rho}{1+\lambda}\,.
	\endaligned \right.
\end{equation}
\end{lema}

%%%%%%%%%%%%%%%%%%
\medskip
%%%%%%%%%%%%%%%%%%

In the next lemma, which generalizes slightly \cite[Lemma 3.3]{Zha08a},
we give the expression of the lower transform of the squared distance 
to a set of two points. The two points, without loss of generality, are assumed to lie 
along a basis vector of $\R^n$, specifically, along the basis vector $e_1\in \R^n$. 
This lemma will be used extensively when we investigate
the behaviour of the multiscale medial axis map with respect to 
perturbations of the boundary of $K$.  

\begin{lema}\label{Sec3.Lem.2Pnt} 
	Assume $n\geq 2$ and let $\{e_1,\,\dots,e_n\}$ be an orthonormal basis of the  Euclidean space $\mathbb{ R}^n$. Let
	$K=\{-\alpha e_1,\,\alpha e_1\}$, where $\alpha>0$. We write $y=e_2y_2+\cdots+e_ny_n\in\mathbb{ R}^{n-1}$ and
	represent the point $xe_1+y\in \R^n$ as the pair $(x,y)$, which therefore denotes the point
	$(x,y_2,\dots, y_n)\in\mathbb{ R}^n$. Then for every $\lambda>0$, we have
\begin{equation}
	C^l_{\lambda}(\dist^2(\cdot;\,K))(x,y)=\left\{\begin{array}{ll}
		\displaystyle \frac{\lambda}{1+\lambda}\alpha^2-\lambda x^2+|y|^2\,,	
					& \displaystyle  |x|\leq \frac{\alpha}{1+\lambda}\,,\\[1.5ex]
		\displaystyle \dist^2((x,y);\,K)\,	    				
					& \displaystyle  |x|\geq \frac{\alpha}{1+\lambda}\,.
		\end{array}\right.
\end{equation}
In particular, 
\begin{equation}
	C^l_{\lambda}(\dist^2(\cdot;\,K))(0,y)=\frac{\lambda}{1+\lambda}\alpha^2+|y|^2<\alpha^2+|y|^2=\dist^2((0,y);\,K)\,.
\end{equation}
\end{lema}

We conclude this section with the definition of $\delta-$neighbourhood of a set,
of Hausdorff distance between two sets \cite{AT04}, and of $\epsilon-$sample of a set \cite{ABE09}.

\begin{defi}
	Given a non-empty subset $E$ of $\mathbb{R}^n$ and $\delta>0$, we define the 
	$\delta$-neighbourhood $E^\delta$ of $E$ by
	\begin{equation*}
		E^\delta=\{x\in \mathbb{R}^n:\; \dist(x;\, E)<\delta\}.
	\end{equation*}
\end{defi}
Note that $E^\delta$ is an open subset of $\R^n$.
\begin{defi} \label{Sec2.Def.HausDist}
	Let $E,\, F$ be non-empty subsets of $\R^n$. The Hausdorff distance
	between $E$ and $F$ is defined in \cite{AT04} by
	\begin{equation}\label{Sec2.Eq.HausDist}
		\dist_{\mathcal{H}}(E,F)=\inf\left\{\delta>0: F\subset E^\delta\; \text{and}\;\, E\subset F^\delta\right\}.
	\end{equation}
\end{defi}
This definition is also equivalent to saying that 
\[
	\dist_{\mathcal{H}}(E,F)=\max\Big\{\sup_{x\in E}\dist(x;\,F),\; \sup_{x\in F}	\dist(x;\,E) \Big\}\,.
\]

\begin{defi}
Let $K$ be a compact subset of $\R^n$. A sample $S$ of the boundary of $K$ is a finite 
set of points of the boundary of $K$,  i.e.  $S\subset \partial K$ and  $\#(S)\in \mathbb{N}$ where $\#(S)$ denotes the 
cardinality of $S$. An $\epsilon-$sample of $\partial K$ is a sample whose Hausdorff
distance to $\partial K$ is less than $\epsilon$, that is, $\dist_{\mathcal{H}}(S,\partial K)<\epsilon$.
\end{defi}

A uniform $\epsilon-$sample $S$ of $\partial K$ is an $\epsilon-$sample of $\partial K$ such that
\begin{equation}
	\dist_{\mathcal{H}}(S,\partial K)<\epsilon\diam(K),
\end{equation}
where the diameter of $K$, $\diam(K)$, is defined as
\[
	\diam(K)=\sup_{x,y\in K}|x-y|\,.
\]

%%%%%%%%%%%%%%%%%%%%%%%%%%%%%%%%%%%%%%%%%%%%%%%%%%%%%%%%%%%%%%%%%%%%%%%%%%%%%%%%%%%%%%%%%%%%%%%%%%%%%
%%%%%%%%%%%%%%%%%%%%%%%%%%%%%%%%%%%%%%%%%%%%%%%%%%%%%%%%%%%%%%%%%%%%%%%%%%%%%%%%%%%%%%%%%%%%%%%%%%%%%
%%%%%%%%%%%%%%%%%%%%%%%%%%%%%%%%%%%%%%%%%%%%%%%%%%%%%%%%%%%%%%%%%%%%%%%%%%%%%%%%%%%%%%%%%%%%%%%%%%%%%

\section{The Multiscale Medial Axis Map} \label{SecMR}
In this section, we define the quadratic multiscale medial axis map $M_{\lambda}( \cdot; K)$, characterize some of its properties,
and establish its relation to the medial axis $M_{K}$. As a by-product, we also infer sharp regularity results 
for the squared distance function $\dist^2(\cdot; K)$, which are of independent interest.

\begin{defi}\label{Sec3.Def.MMAM} 
Let $K\subset \mathbb{R}^n$ be a non-empty closed set.
The {\bf quadratic multiscale medial axis map} of $K$ (medial axis map for short)
with scale $\lambda>0$ is defined for $x\in\mathbb{R}^n$ by
\begin{equation}\label{Eq.Def.MMA}
	M_{\lambda}(x;\,K):=(1+\lambda)R_\lambda(\dist^2(\cdot;\,K))(x)
			=(1+\lambda)\Big(\dist^2(x;\,K))-C_{\lambda}^l(\dist^2(\cdot;\,K)))(x)\Big)\,.
\end{equation}
For a bounded open set $\Omega\subset \mathbb{R}^n$
with boundary $\partial\Omega$, we define the {\bf quadratic multiscale medial axis map} of $\Omega$ 
with scale $\lambda>0$ as
\begin{equation}
	M_{\lambda}(x;\,\Omega):=M_{\lambda}(x;\,\partial\Omega)\qquad x\in\Omega.
\end{equation}
\end{defi}

%%%%%%%%%%%%%%%%%%
\medskip
%%%%%%%%%%%%%%%%%%

\begin{nota}
\begin{itemize}
	\item[$(a)$] The convergence of the lower transform to the original function as
		$\lambda\to \infty$ yields that $\lim_{\lambda\to\infty}R_\lambda(\dist^2(\cdot;\,K))(x)=0$, implying that the values of the ridge transform can be very small when $\lambda>0$ is large. 
		To make the height of our medial axis map
		on the medial axis bounded away from zero, we thus need to scale
		the ridge transform.  The factor $(1+\lambda)$ turns out to be the ``right" scaling factor, as  will  be justified  in Theorem \ref{Sec3.Teo.BndMMA} below, 
		where it will be shown that on the medial axis $M_K$, the medial axis map  
		$M_{\lambda}(x;\,K)$ is bounded both above and below by quantities independent of $\lambda$. 
	\item[$(b)$] The quadratic multiscale medial axis map can also be seen as a morphological operator \cite{Ser82},  equal to 
		the scaled top-hat transform of the squared distance transform with quadratic structuring function. Letting 
		$f(x)=\dist^2(x;\,K)$ and $b_{\lambda}(x)=-\lambda|x|^2$,  it can be shown that the lower transform 
		corresponds to 
		the grayscale opening operator with quadratic structuring function \cite{ZOC14}; i.e.,
		\[
			C_{\lambda}^l(f)=\big(f\ominus b_{\lambda}\big)\oplus b_{\lambda}\,,
		\]
		and thus,
		\[
			M_{\lambda}(\cdot;\,K)=(1+\lambda)\bigg[f-\big(f\ominus b_{\lambda}\big)\oplus b_{\lambda}\bigg]\,.
		\]
		Notwithstanding such an interpretation, it is convenient to view Definition \ref{Sec3.Def.MMAM} in terms of 
		the lower compensated convex transform.
		The exploitation of properties of such transforms permits  a relatively easy evaluation of the geometrical
		properties of  $M_{\lambda}(\cdot;\,K)$ and also permits an easy numerical realization of 
		$M_{\lambda}(\cdot;\,K)$. This relies 
		 on the availability of numerical schemes for computing the lower transform of a given function, which entails
		the availability of schemes to compute the convex envelope of a function. We refer to \cite{ZOC14a} for the algorithmic and 
		implementation details of the schemes for realizing the lower transform of a function.
\end{itemize}
\end{nota}

%%%%%%%%%%%%%%%%%%
\medskip
%%%%%%%%%%%%%%%%%%

We begin with a key quantitative estimate of the tight approximation of the squared distance function $\dist^2(\cdot; K)$ by its lower transform $C^l_\lambda(\dist^2(\cdot;\, K))$.
This result not only underpins our study of the r\^{o}le of $M_{\lambda}$ in characterizing the medial axis $M_K$, but also yields improved locality and regularity properties  of $C^l_\lambda(\dist^2(\cdot;\, K))$ and $\dist^2(\cdot; K)$ respectively (see Corollaries \ref{Sec3.Cor.LclProp}, \ref{Sec3.Cor.Reg} and \ref{Sec3.Cor.Reg.OpenDmn}) which are of interest in their own right.

\begin{teo}\label{Sec3.Thm.TghtApprx} 
Let $K\subset \mathbb{R}^n$ be a non-empty closed set and
denote by $M_K$ the medial axis of $K$. 
Suppose $M_K\not=\varnothing$, $\lambda>0$, assume $x\in \R^n\setminus \overline{M_K}$,  and
let $M_{\lambda}(x;\,K)$ be the multiscale medial axis map of $K$ with scale $\lambda$.
If 
\begin{equation}\label{Sec3.Thm.TghtApprx.Eq1} 
	 \lambda \geq\frac{\dist(x;\,K)}{\dist(x;\,M_K)}\,,
\end{equation}	 
then
\begin{equation}\label{Sec3.Thm.TghtApprx.Eq2}
		\dist^2(x;\,K)=C_{\lambda}^l(\dist^2(\cdot;\,K))(x)\,,
\end{equation}
and consequently,  \begin{equation}
\label{Sec3.Thm.TghtApprx.Eq3}
M_{\lambda}(x;\,K)=0.
\end{equation}
\end{teo}   

%%%%%%%%%%%%%%%%%%
\medskip
%%%%%%%%%%%%%%%%%%

\begin{nota}
Note that  $M_K=\varnothing$ if and only if $K$ is convex (see, for example, \cite[Theorem 2.21]{GiaMod11}), in which case $\dist^2(\cdot;\,K)$ is convex, and 
therefore equals $C_{\lambda}^l(\dist^2(\cdot;\,K))$ in $\R^n$ for all $\lambda > 0$.
\end{nota}

%%%%%%%%%%%%%%%%%%
\medskip
%%%%%%%%%%%%%%%%%%

Now assume $\lambda>0$ and introduce the  set
\begin{equation} \label{Sec3.Def.Neigh}
	V_{\lambda,\,K}= \left\{x\in \mathbb{R}^n: \, \lambda\dist(x;\, M_K)\leq \dist(x;\,K)\right\}\,.
\end{equation}
Clearly $\overline{M_K} \subset V_{\lambda,\,K}$, so this defines a ``neighbourhood" of the medial axis $M_K$ of $K$ (note that it is possible that $\overline{M_K} \cap K \neq \emptyset$, so $V_{\lambda,\,K}$ is not necessarily a neighbourhood in the strict sense), and $V_{\lambda,\,K}$ is a closed set. Moreover, as $\lambda>0$ increases, $V_{\lambda,K}$ describes a family of shrinking sets such that
\begin{equation}
\label{Sec3.prop.vlk}
	\bigcap_{\lambda>0}V_{\lambda,K}=\overline{M_K}\,,
\end{equation}
and if we take the support
of the multiscale medial axis map,  Theorem \ref{Sec3.Thm.TghtApprx}  yields that
\begin{equation}
\label{Sec3.subset}
	\sprt(M_\lambda(\cdot;\,K))\subset V_{\lambda,\,K}\,.
\end{equation}

With the help of (\ref{Sec3.subset}), we can show the following result that characterizes $\overline{M_K}$ in terms of $\sprt(M_\lambda(\cdot;\,K))$.
\begin{coro}\label{Sec3.Cor.Sprt} 
	Suppose $K\subset \mathbb{R}^n$ is a non-empty closed set and $M_K\neq\varnothing$. Then
\begin{equation}\label{Sec3.Cor.Sprt.Eq}
	\bigcap_{0<\lambda<+\infty}\sprt(M_\lambda(\cdot;\,K))=\overline{M_K}\,. 
\end{equation}
\end{coro}

%%%%%%%%%%%%%%%%%%
\medskip
%%%%%%%%%%%%%%%%%%

An important consequence of Theorem \ref{Sec3.Thm.TghtApprx} is  the following 
locality property of the lower transform of the squared distance function. 
This result is also of independent interest, in particular because it quantifies the size of neighbourhood 
needed to evaluate $C^l_{\lambda}(\dist^2(\cdot; K))$, and also because it will be exploited in the proofs 
of Theorems \ref{Sec3.Teo.BndMMA} and \ref{Sec3.Teo.BndMMA} to establish results characterizing the properties of $M_{\lambda}(\cdot; K)$.

\begin{coro}\label{Sec3.Cor.LclProp} (Locality Property) 
Suppose $K\subset\mathbb{R}^n$ is a non-empty closed set. Then for every $x_0\in \mathbb{R}^n$,
\begin{equation}
	C^l_\lambda(\dist^2(\cdot;\, K))(x_0)=\co_{\overline B(x_0;\,r(x_0))}[\dist^2(\cdot;\, K)+\lambda|(\cdot)\,-\,x_0|^2](x_0),
\end{equation}
where
\begin{equation}
	r(x_0):=\frac{2}{\lambda}\dist(x_0;\,K)\,.
\end{equation}
\end{coro}

%%%%%%%%%%%%%%%%%%
\medskip
%%%%%%%%%%%%%%%%%%

\begin{nota} 
\label{Sec3.Rmk.locality}
\begin{itemize}
	\item[(a)] Corollary \ref{Sec3.Cor.LclProp} improves the result in \cite{Kha14}, where the radius of the ball 
		for the locality property is $r(x_0)=6\,\dist(x_0;\,K)/\sqrt{\lambda}$ 
		for $\lambda>2$.
	\item[(b)] In \cite{Kha14}, it was also established that $x_0\in \mathbb{R}^n$ is a stationary point 
		of $C^l_\lambda(\dist^2(\cdot;\,K))$ if and only if $x_0\in \co[K(x_0)]$. We will see that
		the `only if' part of this result is a consequence of arguments from the proof of Theorem \ref{Sec3.Thm.TghtApprx}: see Remark \ref{Sec3.Thm.TghtApprx.Rk}.
\end{itemize}
\end{nota}

%%%%%%%%%%%%%%%%%%
\medskip
%%%%%%%%%%%%%%%%%%

Theorem \ref{Sec3.Thm.TghtApprx} can also be combined with Proposition \ref{Sec2.Th.RegLwTr} to yield a regularity property of the distance transform, which
can be viewed as a weak version of the Lusin theorem for the squared-distance function. 
% We will see in Example \ref{Sec3.Ex.ShrpEst} below that estimate \eqref{Sec3.Cor.Reg.Eq} is in fact sharp.
% It is interesting to compare, for instance, this result with the Luzin theorem obtained by Liu \cite{Liu77} for Sobolev functions.

\begin{coro}\label{Sec3.Cor.Reg} 
	Assume $\lambda>0$.
	Let $K\subset \mathbb{R}^n$ be a non-empty closed set and $V_{\lambda,K}$ the neighbourhood of $M_K$
	defined by \eqref{Sec3.Def.Neigh}. Then 
	\[
		\dist^2(\cdot;\,K)\in C^{1,1}(\mathbb{R}^n\setminus V_{\lambda,K})\,.
	\]	
	Furthermore, for all
	$x,\, y\in \mathbb{R}^n\setminus V_{\lambda,K}$ 
	\begin{equation}\label{Sec3.Cor.Reg.Eq}
		\frac{|D\dist^2(y,\, K) - D\dist^2(z,\, K)|}{
		|y-z|}\leq 2\max\{1,\, \lambda\}\,.
	\end{equation}
\end{coro}

%%%%%%%%%%%%%%%%%%
\medskip
%%%%%%%%%%%%%%%%%%

\begin{nota} It follows from Corollary \ref{Sec3.Cor.Reg} and Stepanov's Theorem \cite{Hei05,Mal99} that the Hessian of
$\dist^2(\cdot;\, K)$ exists almost everywhere in 
$\mathbb{R}^n\setminus\overline{M}_K$. \end{nota}

%%%%%%%%%%%%%%%%%%
\medskip
%%%%%%%%%%%%%%%%%%

Both estimate \eqref{Sec3.Thm.TghtApprx.Eq1} in Theorem \ref{Sec3.Thm.TghtApprx} and estimate 
\eqref{Sec3.Cor.Reg.Eq} for the Lipschitz constant 
in Corollary \ref{Sec3.Cor.Reg} (when $\lambda\geq 1$) are, in fact, sharp, as the following example shows.

\begin{ex}\label{Sec3.Ex.ShrpEst} 
Consider $K=(-1,\, 1)^c=\R\setminus (-1,\,1)\subset \R$. Then $M_K=\{0\}$ {\rm and for } $\lambda>0$,
\[
	C^l_\lambda(\dist^2(\cdot;\, K))(x)=\left\{
					\begin{array}{ll}
				\displaystyle	\frac{\lambda}{\lambda + 1} - \lambda x^2,& \displaystyle |x| \leq  \frac{1}{\lambda + 1},\\[1.5ex]
				\displaystyle	\dist^2(x;\, K),			  & \displaystyle \frac{1}{\lambda + 1}\leq |x| \leq 1,\\[1.5ex]
				\displaystyle	0,					  & \displaystyle |x|\geq 1.
					   \end{array}
				   \right.
\]  
As a result,  $\dist^2(x;\, K)=C^l_\lambda(\dist^2(\cdot;\, K))(x)$ if and only if $|x| \geq  1/(\lambda + 1)$. 
Moreover, for this example we have that $\dist(x;\, M_K)=|x|$, whereas $\dist(x;\,K)=1-|x|$ if $|x|\leq 1$ and $\dist(x;\,K)=0$ if $|x|\geq 1$.
Therefore  \eqref{Sec3.Thm.TghtApprx.Eq1} holds if $|x|\geq 1$. 
If $0<|x|<1$,  \eqref{Sec3.Thm.TghtApprx.Eq1} is given by
\[
	\lambda|x|\geq 1-|x|\Leftrightarrow |x|\geq 1/(1+\lambda)\,,
\]
which shows that  estimate \eqref{Sec3.Thm.TghtApprx.Eq1} is sharp for $K=(-1,1)^c$.

\medskip
Furthermore, for $x_\pm :=\pm 1/(\lambda+1)$,
\[
	D\dist^2(x_\pm,\, K)=\mp4\lambda/(1+\lambda),\quad |x_+-x_-|=2/(1+\lambda)
\]
so that
\[
	\frac{|D\dist^2(x_+,\, K)-D\dist^2(x_-,\, K)|}{|x_+-x_-|}=2\lambda\,.
\]
Thus we can also conclude that estimate \eqref{Sec3.Cor.Reg.Eq} of Corollary \ref{Sec3.Cor.Reg}  is sharp when $\lambda\geq 1$.
\end{ex}

%%%%%%%%%%%%%%%%%%
\medskip
%%%%%%%%%%%%%%%%%%

We cover next the case of a bounded open set $\Omega\subset \mathbb{R}^n$, giving a
 precise statement about  equality of the medial axis maps $M_{\lambda}(x;\,\Omega^c)$ and $M_{\lambda}(x;\,\partial\Omega)$, followed by a modification of Corollary \ref{Sec3.Cor.Reg}.

\begin{prop}\label{Sec3.Pro.MACmp} 
	Suppose $\Omega\subset \mathbb{R}^n$ is a non-empty bounded open set and let $\lambda>0$. Then
	\begin{equation}\label{Eq.Cmp.LwTr.Open}
		C^l_\lambda(\dist^2(\cdot;\,\partial\Omega))(x)=C^l_\lambda(\dist^2(\cdot;\,\Omega^c))(x),\qquad x\in\Omega,
	\end{equation}
	and consequently,
	\begin{equation}\label{eq.SecMR.In} 
		M_{\lambda}(x;\,\partial\Omega)=M_{\lambda}(x;\,\Omega^c),\qquad x\in\Omega.
	\end{equation}
\end{prop}

%%%%%%%%%%%%%%%%%%
\medskip
%%%%%%%%%%%%%%%%%%

\begin{nota}
\begin{itemize}
	\item[$(a)$] Property \eqref{eq.SecMR.In} of the medial axis map is important in many practical 
		situations. For example, in image processing, the objects $\Omega$  of which we wish to find 
		the medial axis might be defined by taking a threshold from a greyscale image, that is,
		as a suplevel set of the image function.
		The object is then represented by a binary image rather than by its boundary. Therefore, in this case,
		it might be more convenient for us to compute numerically the medial axis map $M_{\lambda}(x;\,\Omega^c)$ 
		rather than $M_{\lambda}(x;\,\partial\Omega)$.
	\item[$(b)$] It is worth noting the different qualitative behaviour of the convex envelope and  the 
		compensated convex transform  that appears in 
		\eqref{Eq.Cmp.LwTr.Open}. For $\lambda>0$, the left hand side of \eqref{Eq.Cmp.LwTr.Open} 
		is always positive in $\mathbb{R}^n\setminus\Omega$ whereas the right hand side equals zero in $\Omega^c$.  
		By setting $\lambda=0$, the left hand side of \eqref{Eq.Cmp.LwTr.Open} reduces to the 
		the convex envelope of $\dist^2(x;\,\partial\Omega)$, which vanishes in the convex hull
		$\co[\overline\Omega]$ of the closure of $\Omega$, whereas the right hand side of \eqref{Eq.Cmp.LwTr.Open} 
		gives the convex envelope of $\dist^2(x;\,\Omega^c)$, which is identically zero in $\mathbb{R}^n$.
\end{itemize}
\end{nota}

%%%%%%%%%%%%%%%%%%
\medskip
%%%%%%%%%%%%%%%%%%

For a bounded open set $\Omega$, Corollary \ref{Sec3.Cor.Reg} modifies as follows.

\begin{coro}\label{Sec3.Cor.Reg.OpenDmn} 
	Let $\Omega\subset \mathbb{R}^n$ be a bounded non-empty open set. 
	Then 
	\begin{equation}\label{Sec3.Cor.Reg.OpenDmn.Eq0} 
		\dist^2(\cdot;\,\Omega^c)\in C^{1,1}(\Omega\setminus W_{\lambda,\Omega})\,, 
	\end{equation}		
	where
	\begin{equation}\label{Sec3.Cor.Reg.OpenDmn.Eq1}  
		W_{\lambda,\Omega}:=\left\{x\in \Omega,\; \dist(x; M_{\Omega^c}) \leq \frac{\diam(\Omega)}{\lambda}\right\}\,,
	\end{equation}
and $\diam(\Omega)$ is the diameter of $\Omega$.
Furthermore, for all $x,\, y\in \Omega\setminus W_{\lambda,\Omega}$, 
\begin{equation}\label{Sec3.Cor.Reg.OpenDmn.Eq2}
	\frac{|D\dist^2(y,\, \Omega^c) - D\dist^2(z,\, \Omega^c)|}{|y-z|}\leq 2\max\{1,\, \lambda\}\,.
\end{equation}
\end{coro}

A consequence of Corollary \ref{Sec3.Cor.Reg.OpenDmn} is that outside 
any neighbourhood of $W_{\lambda,\Omega}$, $\dist^2(\cdot;\, \Omega^c)$ is a $C^{1,1}$ function. 
However, we also notice that $\overline{M_K}$ can have positive $n$-dimensional Lebesgue measure. 
Therefore the measure of $W_{\lambda,\Omega}$ might not be small even when $\lambda>0$ is large. 
Corollary \ref{Sec3.Cor.Reg} and Corollary \ref{Sec3.Cor.Reg.OpenDmn} also demonstrate that  
the lower transform can be viewed as a $C^{1,1}$ extension of the squared-distance function 
from the set $V^c_{\lambda,K}$, on which $\dist^2(\cdot; K) = C^l_{\lambda}(\dist^2(\cdot; K))$, to $\R^n$ and from 
$\Omega\setminus W_{\lambda,\Omega}$ to $\Omega$, respectively.

%%%%%%%%%%%%%%%%%%
\medskip
%%%%%%%%%%%%%%%%%%

 Theorem \ref{Sec3.Thm.TghtApprx} showed that if $x \not\in \overline{M_K}$,
then $M_{\lambda}(x; K)=0$ when $\lambda $ is sufficiently large. We now further explore the relationship between the medial axis map $M_{\lambda}(\cdot; K)$
and the medial axis $M_K$, 
both establishing  $\lambda$-independent positive upper and lower bounds on $M_{\lambda}(x;\,K)$ whenever $x\in M_K$, and fully characterizing the limit of $M_{\lambda}(\cdot; K)$ as $\lambda \to \infty$. 
The following geometric structure will play a key r\^{o}le in both Theorem \ref{Sec3.Teo.BndMMA} and Theorem \ref{Sec3.Teo.LimMMA}.

\begin{defi}\label{Sec3.Defi.LimMMA} 
Let $K\subset \R^n$ be a non-empty closed set and for $x\in \R^n$, denote by $K(x)$
the set defined by \eqref{Eq.Def.SetK}, that is, $K(x)=\overline B(x;\,r(x))\cap K$, and 
denote by $\co[K(x)]$ the convex hull of $K(x)$. The {\bf quadratic multiscale medial axis landscape map of $K$} is defined for $x \in \R^n$ by
	\begin{equation}\label{Eq.Def.MMALM} 
		M_{\infty}(x;\,K):=\dist^2(x;\,K)-\dist^2(x;\,\co[K(x)])\,.
	\end{equation}
\end{defi}

%%%%%%%%%%%%%%%%%%
\medskip
%%%%%%%%%%%%%%%%%%

It is straightforward to see that $M_{\infty}(x; K)=0$ if $x \not\in M_K$ but $M_{\infty}(x; K)>0$ for all $x \in M_K$.
Indeed, if $x \not\in M_K$, then there exists $y \in K$ such that $K(x) = \{y\} = \co[K(x)]$ and $\dist^2(x; K)= |x-y|^2$, thus $M_{\infty}(x; K)=0$. 
On the other hand, if $x \in M_K$, then there exist distinct $y_1, y_2 \in K(x)$, so since $\frac{y_1+ y_2}{2} \in 
\co[K(x)]$, we have
$M_{\infty}(x; K) \geq \frac{1}{2}|x-y_1|^2 + \frac{1}{2} |x-y_2|^2 - \left| x - \frac{y_1+ y_2}{2} \right|^2 = \frac{1}{4} |y_1-y_2|^2 \; > \; 0.$

%%%%%%%%%%%%%%%%%%
\medskip
%%%%%%%%%%%%%%%%%%

The next result establishes key bounds on $M_{\lambda}(\cdot, K)$. 

\begin{teo}\label{Sec3.Teo.BndMMA} 
Let $K\subset \mathbb{R}^n$ be a non-empty closed set, and denote by $M_K$ the medial axis of $K$ and by 
$M_{\infty}(x;\,K)$ the quadratic multiscale medial axis landscape map defined by \eqref{Eq.Def.MMALM}.
\begin{itemize}
\item[$(i)$] For every $\lambda>0$ and every $x\in M_K$,  
	\begin{equation}\label{Eq.Bnd.MMAM} 
		M_{\infty}(x;\,K)\leq M_{\lambda}(x;\,K)\leq \dist^2(x;\,K)\,.
	\end{equation}
\item[$(ii)$] For every $\lambda>0$ and for every $x\in \R^n$, 
	\begin{equation} 
		0\leq M_{\lambda}(x;\,K)\leq \dist^2(x;\,K)\,.
	\end{equation}
\end{itemize}
\end{teo}

%%%%%%%%%%%%%%%%%%
\medskip
%%%%%%%%%%%%%%%%%%

The lower bound in \eqref{Eq.Bnd.MMAM} can be expressed in terms of the separation
angle which has been used, for instance, in \cite{SFM07},
for a local geometrical characterization of the medial axis.

\begin{defi} \label{Sec3.Def.OpnAng}
Let $K\subset \mathbb{R}^n$ be a non-empty closed set and denote by $M_K$ the medial axis of $K$.
For $x\in M_K$, let $y_1,\,y_2\in K(x)$ and denote by  $\angle[y_1-x,y_2-x]$
the angle between the two non-zero vectors $y_1-x$ and $y_2-x$, taken between $0$ and $\pi$,
i.e. $\angle[y_1-x,y_2-x]=\cos^{-1}\frac{(y_1-x)\cdot(y_2-x)}{|y_1-x||y_2-x|}$. 
We then define the separation angle  $\theta_x$ for $x\in M_K$ as follows:
\begin{equation}\label{Def.SepAngl}
	\theta_x= \max \Big\{\angle[y_1-x,y_2-x], \quad y_1,\, y_2\in K(x)\Big\}\,.
\end{equation}
\end{defi}

%%%%%%%%%%%%%%%%%%
\medskip
%%%%%%%%%%%%%%%%%%

\begin{nota}\label{Rmrk.SepAngl}
	Recall from Remark \ref{Rk.cpt} that $K(x)$ is compact and hence the supremum of the set 
	$ \{\angle[y_1-x,y_2-x], \quad y_1,\, y_2\in K(x)\}$ is realized   by a pair of distinct points of $K(x)$.
\end{nota}

%%%%%%%%%%%%%%%%%%
\medskip
%%%%%%%%%%%%%%%%%%

\begin{prop}\label{Sec3.Pro.OpAng} 
Let $K\subset \mathbb{R}^n$ be a non-empty closed set,  and denote by $M_K$ the medial axis of $K$.
Then for every $\lambda>0$ and $x\in M_K$, 
\begin{equation}\label{Sec3.Eq.OpAng}  
	\sin^2(\theta_x/2)\dist^2(x;\,K)\leq M_{\lambda}(x;\,K)\leq \dist^2(x;\,K)\,.
\end{equation}
\end{prop}

%%%%%%%%%%%%%%%%%%
\medskip
%%%%%%%%%%%%%%%%%%

\begin{nota}
\begin{itemize}
	\item[$(a)$] For a limit point $x\in \overline{M_K}\setminus M_K$, we have either $x\in K$ or $\theta_x=0$, and the 
		estimate \eqref{Sec3.Eq.OpAng} is valid also in this case.
	\item[$(b)$] Since the bounds in \eqref{Sec3.Eq.OpAng} are independent of $\lambda$, we have that for $x\in M_K$, 
		\begin{equation} 
			\sin(\theta_x/2)\,\,\dist^2(x;\,K)\leq \liminf_{\lambda\to+\infty}M_{\lambda}(x;\,K)
				\leq \limsup_{\lambda\to+\infty}\, M_{\lambda}(x;\,K)\leq \dist^2(x;\,K)\,.
		\end{equation}
	\item[$(c)$] Since for $x\in M_K$, there are at least 
		two different points $y_1,\, y_2\in \partial B(x;\,r(x))\cap  K$, it follows that
		$\theta_x>0$, and hence, comparing with \eqref{Sec3.Eq.OpAng}, we have
		$M_{\lambda}(x;\,K)>0$ for $x\in M_K$. 
\end{itemize}
\end{nota}

%%%%%%%%%%%%%%%%%%
\medskip
%%%%%%%%%%%%%%%%%%

The next result gives the limit behaviour of $M_{\lambda}(x; \, K)$ as $\lambda \to \infty$ for $x \not\in M_K$. 

\begin{prop}\label{Sec3.LimInfSup} 
Let $K\subset \mathbb{R}^n$ be a non-empty closed set and denote by $M_K$ the medial axis of $K$. 
Assume $\lambda>0$ and denote by $M_{\lambda}(x;\,K)$ the medial axis map of $K$ of scale $\lambda$.
Then for $x\in \R^n\setminus M_K$, 
\begin{equation}\label{Sec3.Eq.LimInfSup}  
		\lim_{\lambda\to+\infty}M_{\lambda}(x;\,K)=0\,.
\end{equation}
\end{prop}   

%%%%%%%%%%%%%%%%%%
\medskip
%%%%%%%%%%%%%%%%%%

 \begin{nota} Theorem \ref{Sec3.Thm.TghtApprx}
	clearly implies (\ref{Sec3.Eq.LimInfSup}) in the special case that $x \not\in \overline{M_K}$. 
	But the set $\overline{M_K}\setminus(M_K\cup K)$ may not be empty and we do not know whether 
	$M_\lambda(x,\,K)>0$  for all $\lambda>0$ for $x\in \overline{M_K}\setminus(M_K\cup K)$. 
	So Proposition \ref{Sec3.LimInfSup} is needed  for the proof of the general characterization 
	of $\lim_{\lambda \to \infty} M_{\lambda}(\cdot; \, K)$ that will be given in Theorem \ref{Sec3.Teo.LimMMA}. 	
\end{nota}

%%%%%%%%%%%%%%%%%%
\medskip
%%%%%%%%%%%%%%%%%%

\begin{nota}
By Proposition \ref{Sec3.Pro.OpAng} and Proposition \ref{Sec3.LimInfSup}, we have
\[
	\displaystyle\underset{\lambda\to+\infty}{\liminf}\, M_{\lambda}(x;\,K)\,
	\left\{\begin{array}{ll}
		\displaystyle >0\,,	& \displaystyle  \text{if }x\in M_K	\,,\\[1.5ex]
		\displaystyle = 0\,,	& \displaystyle  \text{if }x\not\in M_K	\,.
		\end{array}\right.
\]
This implies that the limit of $M_{\lambda}(x;\,K)$ can extract exactly the medial axis of $K$. 
If we apply a slightly weaker scaling to the ridge transform, say $(1+x)^{\alpha}$
for $0<\alpha<1$, and define $M_{\lambda}^{\alpha}(x;\,K)=(1+x)^{\alpha}R_\lambda(\dist^2(\cdot;\,K))(x)$,
we have
\[
	\displaystyle\underset{\lambda\to+\infty}{\liminf}\, M_{\lambda}^{\alpha}(x;\,K)\,
	\left\{\begin{array}{ll}
		\displaystyle +\infty\,,	& \displaystyle  \text{if }x\in M_K	\,,\\[1.5ex]
		\displaystyle = 0\,,	& \displaystyle  \text{if }x\not\in M_K	\,,
		\end{array}\right.
\]
that is, $M_{\lambda}^{\alpha}(x;\,K)$ approaches the indicator function of $M_K$ \cite{HUL01} as $\lambda$ 
becomes large. 
\end{nota}

%%%%%%%%%%%%%%%%%%
\medskip
%%%%%%%%%%%%%%%%%%

We can now characterize the limit of $M_\lambda(x;\,K)$ as $\lambda\to+\infty$ for all $x \in \R^n$. 

\begin{teo}\label{Sec3.Teo.LimMMA} 
	Suppose $K\subset \mathbb{R}^n$ be a non-empty closed set.
	Then for every $x\in \mathbb{R}^n$,
\begin{equation} 
	\lim_{\lambda\to+\infty}M_\lambda(x;\,K)=M_{\infty}(x;\,K)\,.
\end{equation}
\end{teo}

%%%%%%%%%%%%%%%%%%
\medskip
%%%%%%%%%%%%%%%%%%

From Theorem \ref{Sec3.Thm.TghtApprx}, Corollary \ref{Sec3.Cor.Reg} 
and Theorem \ref{Sec3.Teo.LimMMA}, it follows that when $\lambda>0$ is increasing, the support of 
$M_\lambda(\cdot;\,K)$ is contained in a shrinking neighbourhood of $\overline{M_K}$ 
and approaches the multiscale medial axis landscape map $M_{\infty}(x;\, K)$. The numerical advantage of studying 
$M_\lambda(\cdot;\,K)$ as an approximation of the multiscale medial axis landscape map $M_{\infty}(\cdot;\, K)$ is that it 
relies only on the computation of the lower compensated convex transform of the squared distance transform, 
whose construction is local by virtue of Corollary \ref{Sec3.Cor.LclProp}, whereas the computation of $M_{\infty}(\cdot;\, K)$
is difficult because we need to evaluate the convex hull $\co[K(x)])$.

%%%%%%%%%%%%%%%%%%
\medskip
%%%%%%%%%%%%%%%%%%

\begin{nota}
A further consequence of Theorem \ref{Sec3.Teo.LimMMA} is that for every fixed $x\in \R^n$ and for every 
	non-empty closed set $K\subset\R^n$, the family of lower transforms 
	$\lambda\mapsto C^l_\lambda(\dist^2(\cdot;\, K))(x)$ is `differentiable' at infinity. 
	If we let $\epsilon=1/\lambda$, $g(\epsilon;\,x)=C^l_{1/\epsilon}(\dist^2(\cdot;\, K))(x)$, and
	$g(0;\,x)=\lim_{\epsilon\to 0}g(\epsilon;\,x)=\dist^2(x;\, K)$, then
	\[
		\lim_{\epsilon\to0+}\frac{g(\epsilon;\,x)-g(0;\,x)}{\epsilon}=-M_{\infty}(x;\,K)\,,
	\]
	so we have the asymptotic expansion
	\begin{equation}
		C^l_\lambda(\dist^2(\cdot;\, K))(x)= \dist^2(x;\, K)-
		\frac{M_{\infty}(x;\,K)}{1+\lambda}+o\left(\frac{1}{1+\lambda}\right)
	\end{equation}
	when $\lambda\to\infty$. 
	\end{nota}
	
%%%%%%%%%%%%%%%%%%
\medskip
%%%%%%%%%%%%%%%%%%

\begin{nota}
In general, $x\mapsto \dist^2(x;\, \co[K(x)])$ is not continuous in $\R^n$ as $x$ 
approaches the medial axis $M_K$. But from Theorem \ref{Sec3.Teo.LimMMA}, we can show that for every $x\in\R^n$,
\begin{equation}\label{Sec3.AsyExp.Eq01}
	\lim_{\lambda\to+\infty}\Big\{(1+\lambda)C^l_\lambda(\dist^2(\cdot;\, K))(x)-
	\lambda\dist^2(x;\, K)\Big\}=\dist^2(x;\, \co[K(x)])
\end{equation}
using the equality
\begin{equation}\label{Sec3.AsyExp.Eq02}
	\dist^2(x;\, K)-M_\lambda(x;\, K)=(1+\lambda)C^l_\lambda(\dist^2(\cdot;\, K))(x)-\lambda\dist^2(x;\, K)\,.
\end{equation}
For  large $\lambda>0$, \eqref{Sec3.AsyExp.Eq01} can be viewed as an approximation of 
$\dist^2(x;\, \co[K(x)])$ by continuous functions. The function $\dist(x;\, \co[K(x)])$ 
has been used, for instance, for surface reconstruction when $K\subset\R^3$ is finite \cite{Dey06}.
While it is difficult in general to calculate  $\dist^2(x;\, \co[K(x)])$ directly, we see from \eqref{Sec3.AsyExp.Eq01} that
the numerical computation of $C^l_\lambda(\dist^2(\cdot;\, K))(x)$, whose evaluation involves only
local convex envelope calculations because of Corollary \ref{Sec3.Cor.LclProp}, offers an easy
 approximation of $\dist^2(x;\, \co[K(x)])$.
 \end{nota}

%%%%%%%%%%%%%%%%%%
\medskip
%%%%%%%%%%%%%%%%%%

We conclude this section by observing briefly that, based on the estimates of Theorem \ref{Sec3.Teo.BndMMA} and 
Proposition \ref{Sec3.Pro.OpAng}, it is reasonable to define an alternative medial
axis map by taking the square root of $M_{\lambda}(x;\,K)$.

\begin{defi} 
	We define the {\bf multiscale medial axis map of linear growth} (linear
	medial axis map for short) by 
	\begin{equation} 
		M^1_\lambda(x;\,K):=\sqrt{M_{\lambda}(x;\,K)}
	\end{equation}
	for $x\in \mathbb{R}^n$.
\end{defi}

From Proposition \ref{Sec3.Pro.OpAng}, we obtain that the height of this linear medial axis map 
is `proportional' to the distance function itself; that is, for $\lambda>0$, we have 
\begin{equation}
	\sin(\theta_x/2)\,\dist(x;\,K)\leq M^1_\lambda(x;\,K)\leq \dist(x;\,K)\quad\text{for }x\in M_K
\end{equation}
and
\begin{equation}
	0\leq M^1_\lambda(x;\,K)\leq \dist(x;\,K)\quad\text{for }x\in \R^n\,.
\end{equation}

Note that the linear medial axis map $M^1_\lambda(x;\,K)$ is different from a definition based on the lower 
compensated convex transform for the distance function $\dist(x;\,K)$ itself, 
i.e. based on $R_{\lambda}(\dist(\cdot;\,K))$. 
Of course, we can define such maps  using the $p$-distance function $\dist^p(x;\,K)$ for any 
$1\leq p<\infty$. But in this paper, we focus mainly on the medial axis map $M_{\lambda}(\cdot; K)$ defined using the 
squared-distance function, i.e., for $p=2$, 
in which case the geometry of $M_{\lambda}(x;\,K)$ is  easy to control. 
For instance, as we will see in the next section, $M_{\lambda}(x;\,K)$ has the same height along the parts 
of the medial axis generated by two points.
This is a key property when one looks for approximate medial axes by applying the Voronoi diagram method
of finite $\epsilon$-samples. 

%%%%%%%%%%%%%%%%%%%%%%%%%%%%%%%%%%%%%%%%%%%%%%%%%%%%%%%%%%%%%%%%%%%%%%%%%%%%%%%%%%%%%%%%%%%%%%%%%%%%%%%%%%%%%%%
%%%%%%%%%%%%%%%%%%%%%%%%%%%%%%%%%%%%%%%%%%%%%%%%%%%%%%%%%%%%%%%%%%%%%%%%%%%%%%%%%%%%%%%%%%%%%%%%%%%%%%%%%%%%%%%
%%%%%%%%%%%%%%%%%%%%%%%%%%%%%%%%%%%%%%%%%%%%%%%%%%%%%%%%%%%%%%%%%%%%%%%%%%%%%%%%%%%%%%%%%%%%%%%%%%%%%%%%%%%%%%%

\section{Hausdorff Stability}\label{SecHAU}
Quantifying the instability of the medial axis is of fundamental importance for both theory and computation.
This aspect becomes more and more relevant in practice nowadays, given that point clouds are increasingly being used 
for geometric modeling over a wide range of applications. Moreover, there are computational 
approaches, such as the Voronoi diagram method, which search for a continuous approximation of the medial axis 
of a shape starting from subsets of the Voronoi diagram of a sample of the shape boundary. 
The presence of noise on the boundary, and/or the discrete character of samples of the boundary shape thus call 
for methods that permit the control of the parts of the medial axis which are not stable.
In this section, we will discuss how this aspect is tackled by the multiscale medial axis map.
In the first part of the section, we examine the values of $M_{\lambda}(\cdot;\,K)$ when the distance of the 
point to the boundary of the set is achieved by two points, whereas in the second part we discuss
the Hausdorff stability of $M_{\lambda}(\cdot;\,K)$.

\begin{prop}\label{Sec4.Pro.MA2Pnt} 
Assume $n\geq 2$ and let $\{e_1,\,\dots,e_n\}$ be an orthonormal basis of the  Euclidean space $\mathbb{ R}^n$. Let
	$K=\{-\alpha e_1,\,\alpha e_1\}$, where $\alpha>0$. We write $y=e_2y_2+\cdots+e_ny_n\in\mathbb{ R}^{n-1}$ and
	represent the point $xe_1+y\in \R^n$ as the pair $(x,y)\in \R\times \R^{n-1}$, which therefore denotes the point
	$(x,y_2,\dots, y_n)\in\mathbb{ R}^n$. Then for every $\lambda>0$, we have
\begin{equation}\label{Eq.MA2Pnt}
	M_{\lambda}((x,y);\,K)=\left\{\begin{array}{ll}
		\displaystyle (1+\lambda)^2\left(|x|-\frac{\alpha}{1+\lambda}\right)^2\,,& 
						\displaystyle	|x|\leq  \frac{\alpha}{1+\lambda}\,,\\[1.5ex]
		\displaystyle 0\,, & \displaystyle |x|\geq \frac{\alpha}{1+\lambda}.
	\end{array}\right.
 \end{equation}
 \end{prop}
 
\begin{nota} 
The medial axis map $M_{\lambda}((x,y);\,K)$  reaches its maximum on the medial axis of $K$, at the point 
$(0,y)=y_2e_2+\cdots +y_ne_n\in \mathbb{R}^{n-1}$, attaining the value $M_{\lambda}((0,y);\,K)= \alpha^2$.
Note that $\alpha>0$ is half the distance between the two points $-\alpha e_1$ and $\alpha e_1$ of $K$. 
Another important observation is that 
$M_{\lambda}((x,y);\,K)$ is a function of the $x$-variable only, and does not change its height along 
the $y$ direction. Therefore, on branches of the medial axis generated by two points, the height remains 
the same.
If $\alpha$ is small (equivalently, the two points in $K$
are close to each other), the values of $M_{\lambda}((x,y);\,K)$ will be uniformly small. 
\end{nota}
 
We next give the Hausdorff stability property of the multiscale medial axis map, followed by some comments on 
 implications of this property for the localization of the medial axis of a domain.

\begin{teo}\label{Sec4.Teo.HAS} 
	Assume $\lambda>0$. Let $K,\, L\subset \R^n$ be non-empty compact sets. 
	Then as $L\to K$ under the Hausdorff distance,
	$M_{\lambda}(\cdot;\,L)\to M_{\lambda}(\cdot;\,K)$ uniformly in every fixed bounded set in $\mathbb{R}^n$.
	More precisely, if we let $\mu:=\dist_{\mathcal{H}}(K,\,L)$ be the Hausdorff distance between $K$ and $L$, 
	then for $x\in \mathbb{R}^n$, 
\begin{equation}\label{Sec4.Eq.Teo.HASLwtr} 
	\Big|C^l_\lambda(\dist^2(\cdot;\,K))(x)-C^l_\lambda(\dist^2(\cdot;\,L))(x)\Big|\leq 
	\mu\Big((\dist(x;\,K)+\mu)^2+1+\mu\Big)\,,
\end{equation} 
and
\begin{equation}\label{Sec4.Eq.Teo.HAS} 
	\Big|M_{\lambda}(x;\,K)-M_{\lambda}(x;\,L)\Big|\leq \mu(1+\lambda)
		\Big((\dist(x;\,K)+\mu)^2+2\dist(x;\,K)+2\mu+1\Big)\,.
\end{equation}
\end{teo}

\begin{nota}
\begin{itemize}
	\item[$(a)$]
		From Theorem \ref{Sec4.Teo.HAS} we also conclude that for any compact sets $K_1$, $K_2$
		\[
			\Big|M_{\lambda}(x;\,K_1)-M_{\lambda}(x;\,K_2)\Big|\leq \min_{i=1,2}\Big\{ \mu(1+\lambda)
			\Big((\dist(x;\,K_i)+\mu)^2+2\dist(x;\,K_i)+2\mu+1\Big)\Big\}\,,
		\]
		which shows that the medial axis map is uniformly continuous on compact sets with respect to the Hausdorff metric.
	\item[$(b)$]
		While the medial axis of $K$ is not a stable structure with respect to the Hausdorff distance, 
		its medial axis map $M_{\lambda}(x;\,K)$ is by contrast  a stable structure. This result
		complies with \eqref{Sec4.Eq.Teo.HAS} which shows that as $\lambda$ becomes large, 
		the bound in \eqref{Sec4.Eq.Teo.HAS} becomes large.
\end{itemize}
\end{nota}

As an immediate consequence of Theorem \ref{Sec4.Teo.HAS}
we have the following result, which relates the medial axis map of  
the boundary of a domain $\Omega$ with that of its $\epsilon$-samples $K_\epsilon$. 

\begin{coro}\label{Sec4.Cor.HAS}  
Assume $\lambda>0$. Let $\Omega\subset\mathbb{R}^n$ be a bounded open set with diameter $\diam(\Omega)$. 
Suppose $K_\epsilon\subset \mathbb{R}^n$ is a compact set such that 
$\dist_{\mathcal{H}}(\partial\Omega,\,K_\epsilon)\leq \epsilon $. Then for $x\in \R^n$
\begin{equation}
	\Big|C^l_\lambda(\dist^2(\cdot;\,\partial\Omega))(x)-C^l_\lambda(\dist^2(\cdot;\,K_\epsilon))(x)\Big|\leq 
	\epsilon\Big((\diam(\Omega)+\epsilon)^2+1+\epsilon\Big)\,,
\end{equation} 
and
\begin{equation} 
	\Big|M_{\lambda}(x;\,\partial\Omega)-M_{\lambda}(x;\,K_\epsilon)\Big|\leq 
	\epsilon(1+\lambda)\Big(
	(\diam(\Omega)+\epsilon)^2+2\diam(\Omega)+2\epsilon+1\Big)\,,
\end{equation} 
for all $x\in \overline\Omega$.
\end{coro}

%%%%%%%%%%%%%%%%%%
\medskip
%%%%%%%%%%%%%%%%%%

\begin{nota}\label{Sec4.Rem.VorPol}
\begin{itemize}
	\item[$(a)$] If we consider an $\epsilon$-sample $K_\epsilon$ of $\partial\Omega$, that is, a discrete set of points such that 
		$\dist_{\mathcal{H}}(\partial\Omega,K_\epsilon)\leq \epsilon $, Corollary \ref{Sec4.Cor.HAS} yields 
		a simple criteria that permits the 
		suppression of those parts of the Voronoi diagram of $K_{\epsilon}$ that are not related in the limit, 
		 as $\epsilon \to 0$, to the stable parts of the medial axis of $\Omega$. 
	\item[$(b)$] Since the medial axis of $K_{\epsilon}$ is the Voronoi diagram of $K_{\epsilon}$, 
		if $V_{\epsilon}$ denotes the set of all the vertices
		of the Voronoi diagram $\mathcal{V}or(K_{\epsilon})$ of $K_{\epsilon}$,  and  
		$P_{\epsilon}$ is the subset of $V_{\epsilon}$ formed by the poles of $\mathcal{V}or(K_{\epsilon})$ 
		introduced in \cite{AB99}, (i.e. those
		vertices of $\mathcal{V}or(K_{\epsilon})$ that converge to the medial axis of $\Omega$ as the sample 
		density approaches infinity), 
		then as a result of Proposition \ref{Sec3.LimInfSup}, for $\lambda>0$, we conclude that 
		\begin{equation}
			\lim_{\epsilon\to 0+}\,M_{\lambda}(x_{\epsilon};	\,K_{\epsilon})=0
			\quad\text{for }x_{\epsilon}\in V_{\epsilon}\setminus P_{\epsilon}\,.
		\end{equation}
		Since as $\epsilon \to 0+$, $K_{\epsilon}\to \partial \Omega$,
		and knowing that 
		$P_{\epsilon}\to M_{\Omega}$ \cite{ACK01b,BC00}, 
		then on the vertices of 
		$\mathcal{V}or(K_{\epsilon})$ that do not tend to $M_{\Omega}$, $M_{\lambda}(x_{\epsilon};\,K_{\epsilon})$ 
		must approach zero in the limit because of Proposition \ref{Sec3.LimInfSup}.
		As a result, in the context of the methods of approximating the medial axis starting from the Voronoi diagram of a sample
		(such as those described in \cite{ACK01b,Dey06,DZ04,SP08}), the use of the multiscale medial axis map offers 
		an alternative and much easier 
		tool to construct continuous approximations to the medial axis with guaranteed convergence as $\epsilon \to 0+$.
\end{itemize}
\end{nota}

%%%%%%%%%%%%%%%%%%
\medskip
%%%%%%%%%%%%%%%%%%

\begin{figure}[htbp]
\centerline{
	$\begin{array}{cc}
		\includegraphics[height=0.35\textwidth]{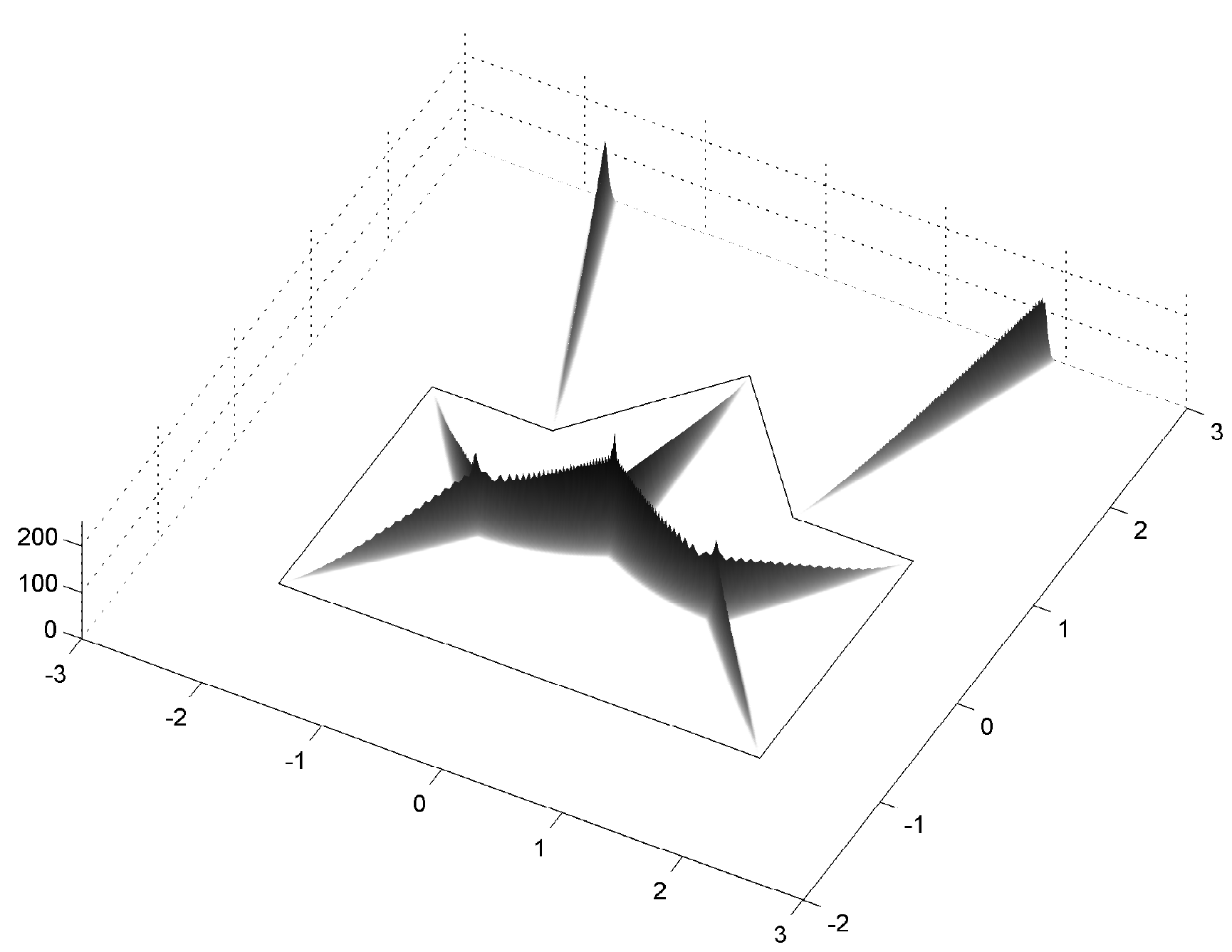}&
		\includegraphics[height=0.35\textwidth]{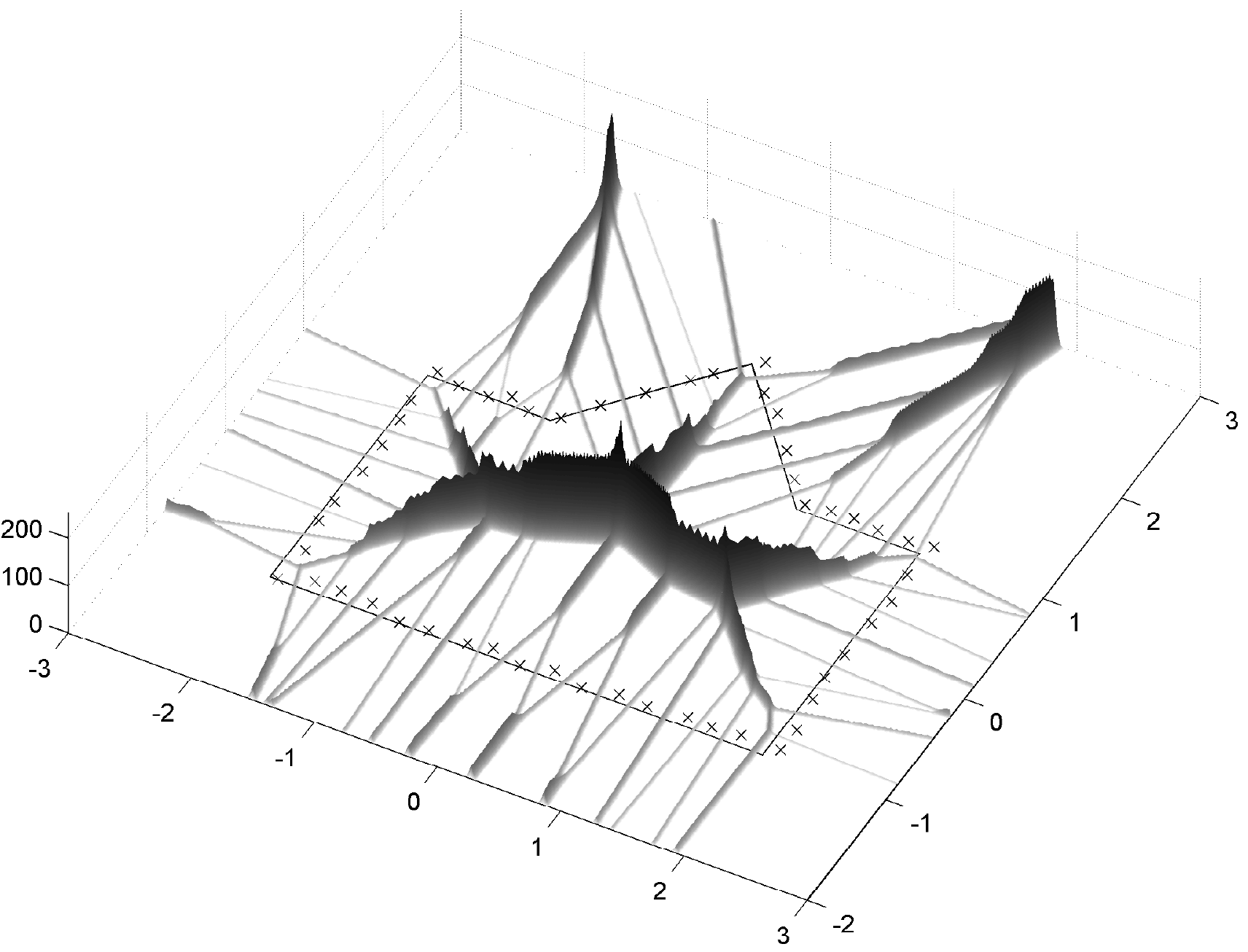}\\
		(a) & (b)\\
		\includegraphics[height=0.35\textwidth]{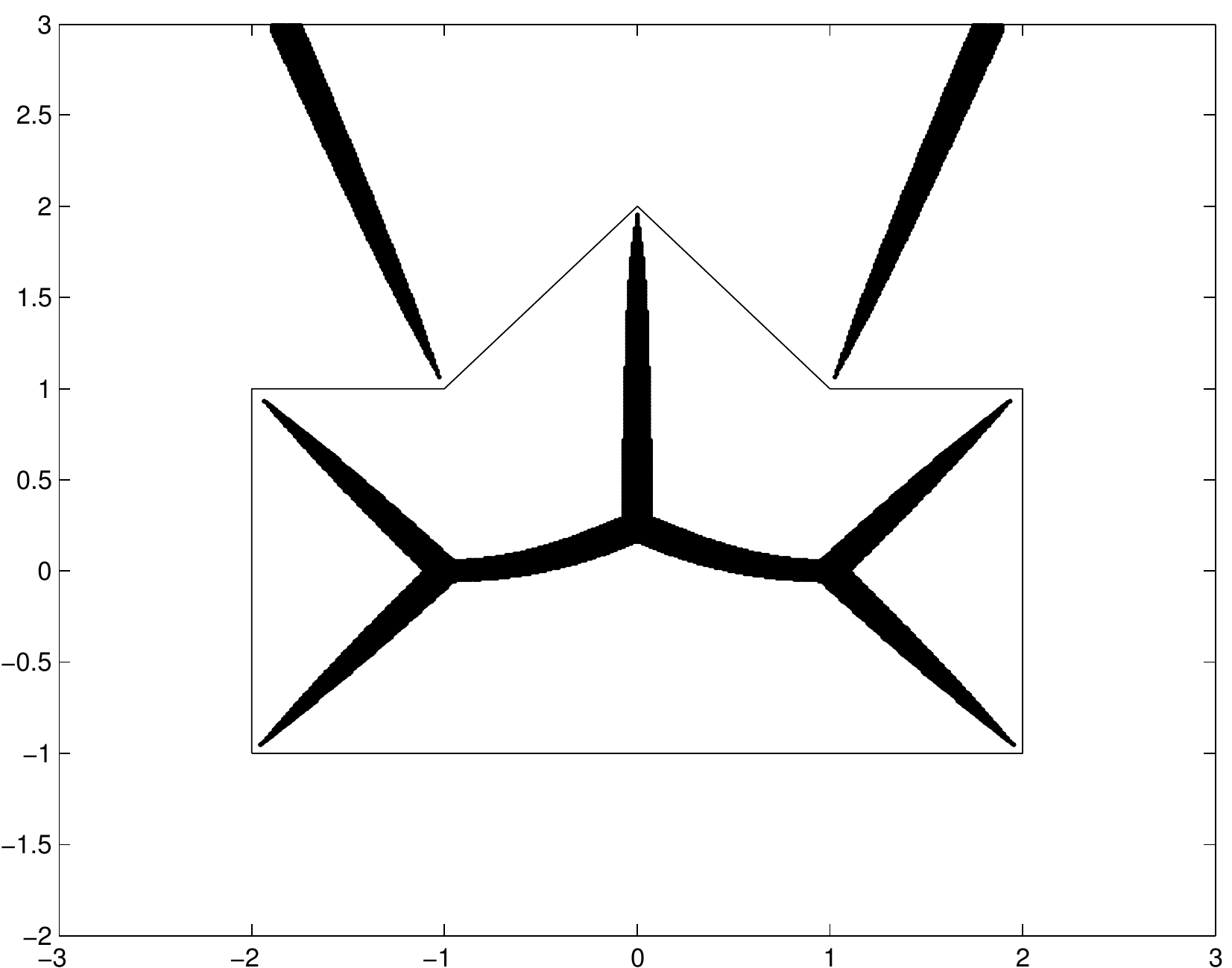}&
		\includegraphics[height=0.35\textwidth]{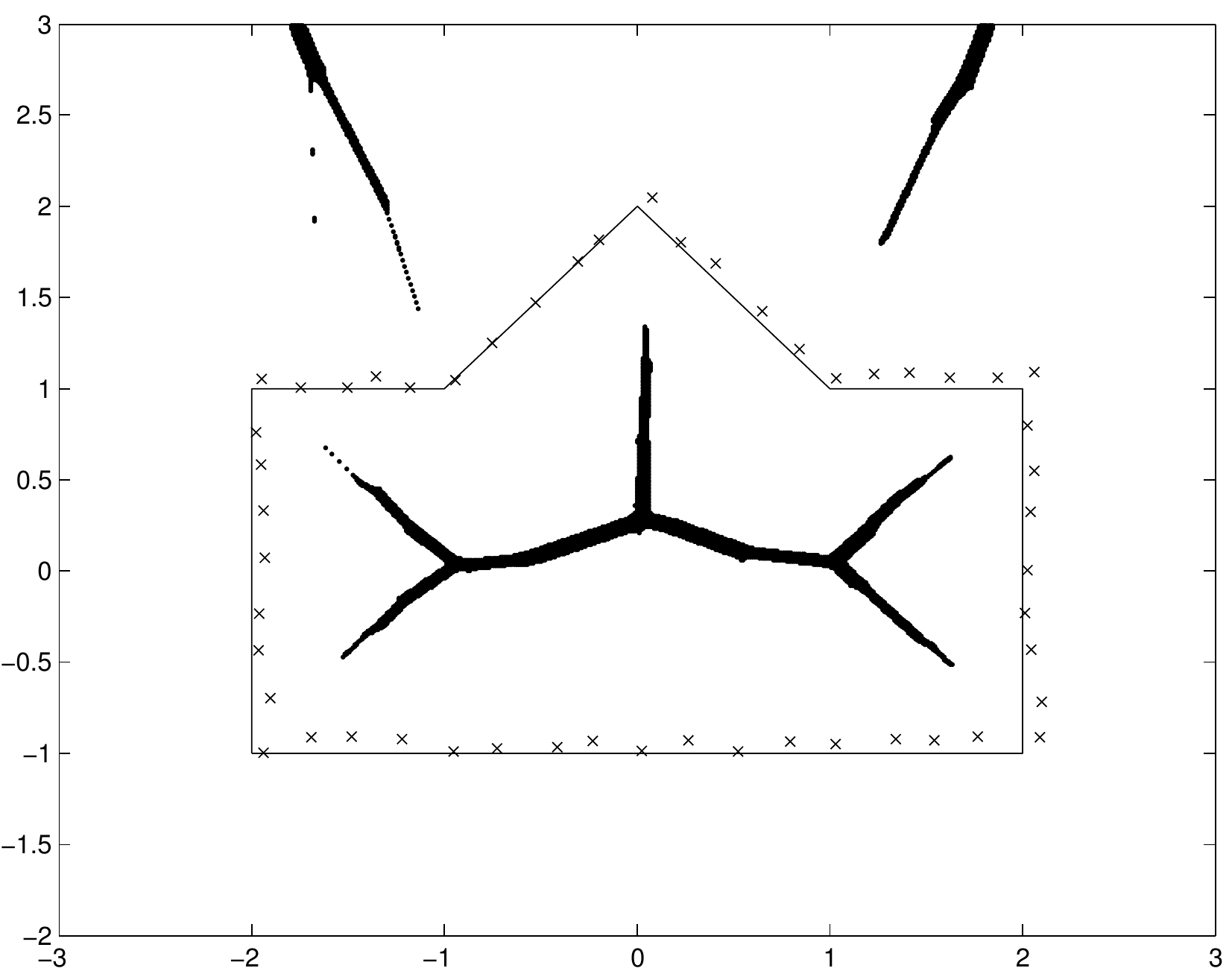}\\
		(c) & (d)
	\end{array}$
}
\caption{\label{Sec4.FigHaus} 
	Multiscale Medial Axis Map of a nonconvex domain $\Omega$ 
	and of an $\epsilon$-sample $K_{\epsilon}$ of its boundary.
	$(a)$ Graph of $M_{\lambda}(\cdot;\,\Omega)$ for $\lambda=5$;
	$(b)$ Graph of $M_{\lambda}(\cdot;\,K_{\epsilon})$; 
	$(c)$ Support of $M_{\lambda}(\cdot;\,\Omega)$;
	$(d)$ Suplevel set of $M_{\lambda}(x;\,K_{\epsilon})$ for a threshold equal 
	to $\displaystyle 0.15\max_{x\in \R^2} \{M_{\lambda}(x;\,K_{\epsilon})\}$.
	}
\end{figure}

With the aim of giving insight into the implications of the Hausdorff stability of $M_{\lambda}(x;\,\partial\Omega)$
and Corollary \ref{Sec4.Cor.HAS}, we display in Figure \ref{Sec4.FigHaus} the graph of the
multiscale medial axis map of a non-convex domain $\Omega$ and of an $\epsilon$-sample $K_{\epsilon}$
of its boundary.
Inspection of the graph of $M_{\lambda}(x;\,\partial\Omega)$ and $M_{\lambda}(x;\,K_{\epsilon})$, 
displayed in Figure \ref{Sec4.FigHaus}$(a)$ and Figure \ref{Sec4.FigHaus}$(b)$, reveals
that both functions take comparable values along the main branches of $M_{\Omega}$. Also,
$M_{\lambda}(x;\,K_{\epsilon})$ takes small values along the secondary branches, generated by the 
sampling of the boundary of $\Omega$. These values can therefore be filtered out
by simple thresholding so that a stable approximation of the medial axis of $\Omega$ can be
computed. This can be appreciated by looking at Figure \ref{Sec4.FigHaus}$(d)$, which displays a suplevel
set of $M_{\lambda}(x;\,K_{\epsilon})$ that appears to be a reasonable approximation of the support of
$M_{\lambda}(x;\,\partial\Omega)$ shown in Figure \ref{Sec4.FigHaus}$(c)$.

%%%%%%%%%%%%%%%%%%%%%%%%%%%%%%%%%%%%%%%%%%%%%%%%%%%%%%%%%%%%%%%%%%%%%%%%%%%%%%%%%%%%%%%
%%%%%%%%%%%%%%%%%%%%%%%%%%%%%%%%%%%%%%%%%%%%%%%%%%%%%%%%%%%%%%%%%%%%%%%%%%%%%%%%%%%%%%%
%%%%%%%%%%%%%%%%%%%%%%%%%%%%%%%%%%%%%%%%%%%%%%%%%%%%%%%%%%%%%%%%%%%%%%%%%%%%%%%%%%%%%%%

\section{Examples of Exact Medial Axis Maps and Their Supports}\label{SecEx}

In this section we illustrate the behaviours of our multiscale medial axis map for some $2d$ geometric objects
$K$, for which it is possible to obtain an explicit analytical expression for 
$M_{\lambda}(x;\,K)$.  
Thanks to the translation and the partial rotation invariance property of the convex envelope 
\cite[Proposition 2.3, 2.10]{ZOC14}, it is then possible to derive an explicit analytical expression for 
$M_{\lambda}(x;\,K)$ in the case that $K$ is a $3d$-solid obtained by, for instance, rotations or translations  
of the models considered in 
this section. For the sake of conciseness, we leave the derivations to interested readers. 

Though the derivation here is limited only to $2d$ geometric models, 
these models retain, nevertheless, their basic geometric features, because they are able
to show that $M_{\lambda}(x;\,K)$ can, in fact, provide an accurate 
and stable way to find $M_K$, the medial axis of $K$, and represents likewise an effective tool
to analyze the geometry and structure of $M_K$. We will also see how it is possible to select either 
the main stable parts of $M_K$ or to locate its fine parts by using suplevel sets of 
$M_{\lambda}(x;\,K)$.

%%%%%%%%%%%%%%%%%%%%%%%%%%%%%%%%%%%
\begin{ex}\label{Sec5.Ex01}
We consider the case of a four-point set $K\subset \R^2$ defined as follows. Let $b, \epsilon >0$
with $\epsilon \in (0,\,1)$, set $c=\epsilon b$. Define then $K=\{(b,c),\,(b,-c),\,(-b,c),\,(-b,-c)\}$.
For this set, we have  
\begin{equation}\label{Ex01.dstK4}
	\dist^2((x,y);\, K)=(|x|-b)^2+(|y|-c)^2\,,
\end{equation}
and, after some lengthy calculations based on the construction of affine functions, we can show that the lower transform can
be expressed as follows
\begin{equation}\label{Ex01.LwTrK4}
	C^l_{\lambda}(\dist^2(\cdot;\,K))(x,y)=(1+\lambda)g(x,y)+\frac{\lambda}{1+\lambda}(b^2+c^2)-\lambda(x^2+y^2)
\end{equation}
where the auxiliary function $g=g(x,y)$ is a continuous piecewise quadratic function defined as follows
\begin{equation}\label{Ex01.AxFnctK4}
	g(x,y)=\left\{
			\begin{array}{ll} 
		\displaystyle	0				&
					\displaystyle \text{if }\;|x|\leq \frac{b}{1+\lambda},\;\;|y|\leq \frac{c}{1+\lambda};\\[2ex]
		\displaystyle	\left(|x|-\frac{b}{1+\lambda}\right)^2	&
					\displaystyle  \text{if }\;|x|\geq \frac{b}{1+\lambda},\;\;|y|\leq \frac{c}{1+\lambda};\\[2ex]
		\displaystyle	\left(|y|-\frac{c}{1+\lambda}\right)^2	&
					\displaystyle  \text{if }\;|x|\leq \frac{b}{1+\lambda,}\;\;|y|\geq \frac{c}{1+\lambda};\\[2ex]
		\displaystyle	\left(|x|-\frac{b}{1+\lambda}\right)^2 + \left(|y|-\frac{c}{1+\lambda}\right)^2&
					\displaystyle  \text{if }\;|x|\geq \frac{b}{1+\lambda},\;\;|y|\geq \frac{c}{1+\lambda}\,.			
			\end{array}
			\right.
\end{equation}
The multiscale medial axis map $M_{\lambda}(x;\,K)$ is then computed using the definition \eqref{Eq.Def.MMA}.
In particular, for this example, after some algebraic rearrangements, it is possible to show that 
since all four points in $K$ lie on a circle centerd at the origin, the medial axis map of $K$ 
can be expressed as
\begin{equation}\label{Ex01.MMAK4}
	M_{\lambda}((x,y);\,K)=(1+\lambda)^2\Big(\dist^2((x,y);\,K/(1+\lambda))-\dist^2((x,y);\,\co(K)/(1+\lambda))\Big)
	\quad (x,y)\in\mathbb{R}^2\,,
\end{equation}
where for $Z\subset \R^n$ and $\alpha\in\R\setminus\{0\}$, we use the notation $Z/\alpha$ to denote the set 
$\{w\in \R^n: w=z/\alpha\text{ for }z\in Z\}$. 
By a closer inspection of \eqref{Ex01.MMAK4}, we can make then the following observations:
\begin{itemize}
\item[$(i)$] The support of $M_{\lambda}(\cdot;\,K)$ is 
	\[
		\sprt(M_{\lambda}(\cdot;\,K))=\left\{(x,y)\in \mathbb{R}^2:\, |x|\leq \frac{2b}{1+\lambda}\;\text{ or }\; 
		|y|\leq \frac{2\epsilon b}{1+\lambda}\right\}\,.
	\]
	The `thickness' of the support for the main branch $y$-axis is, therefore, $2b/(1+\lambda)$ while that for 
	the minor branch $x$-axis is $2\epsilon b/(1+\lambda)$. 
\item[$(ii)$] The height of the medial axis map along the main branch $y$-axis when $|y|\geq \epsilon b/(1+\lambda)$ is 
	$b^2$ while the value along $x$-axis when $|x|\geq b/(1+\lambda)$ is $\epsilon^2 b^2$.
\item[$(iii)$] At the only Voronoi vertex $0$ of the Voronoi diagram of $K$, the value is 
	$M_{\lambda}(0;\,K)=b^2(1+\epsilon^2)$.
\end{itemize}
Figure \ref{Sec5.Fig4pt} displays the graph of $M_{\lambda}(x;\,K)$ as given by \eqref{Ex01.MMAK4} for different values 
of $\lambda$ and for the set $K$ defined by $b=2$ and $\epsilon=0.5$. 
For each value of $\lambda$, we can easily verify the presence of two scales in $M_{\lambda}(x;\,K)$:
a strong one which is reflected by the values of $M_{\lambda}(x;\,K)$ along the $y-$axis generated by the two-point
set $\{(-b,\,c),\,(b,\,c)\}$ and a weak one captured by the value of $M_{\lambda}(x;\,K)$ along the
$x-$axis generated by the two-point set $\{(b,\,-c),\,(-b,\,-c)\}$.
In agreement  with our theoretical results, we also verify that the support of the continuous function
$M_{\lambda}(x;\,K)$ contains the medial axis of $K$, given by the Voronoi diagram of $K$
in this case, and such support shrinks to $M_K$ as $\lambda$ increases. 
\begin{figure}[htbp]
\centerline{
	$\begin{array}{cc}
		\includegraphics[height=0.35\textwidth]{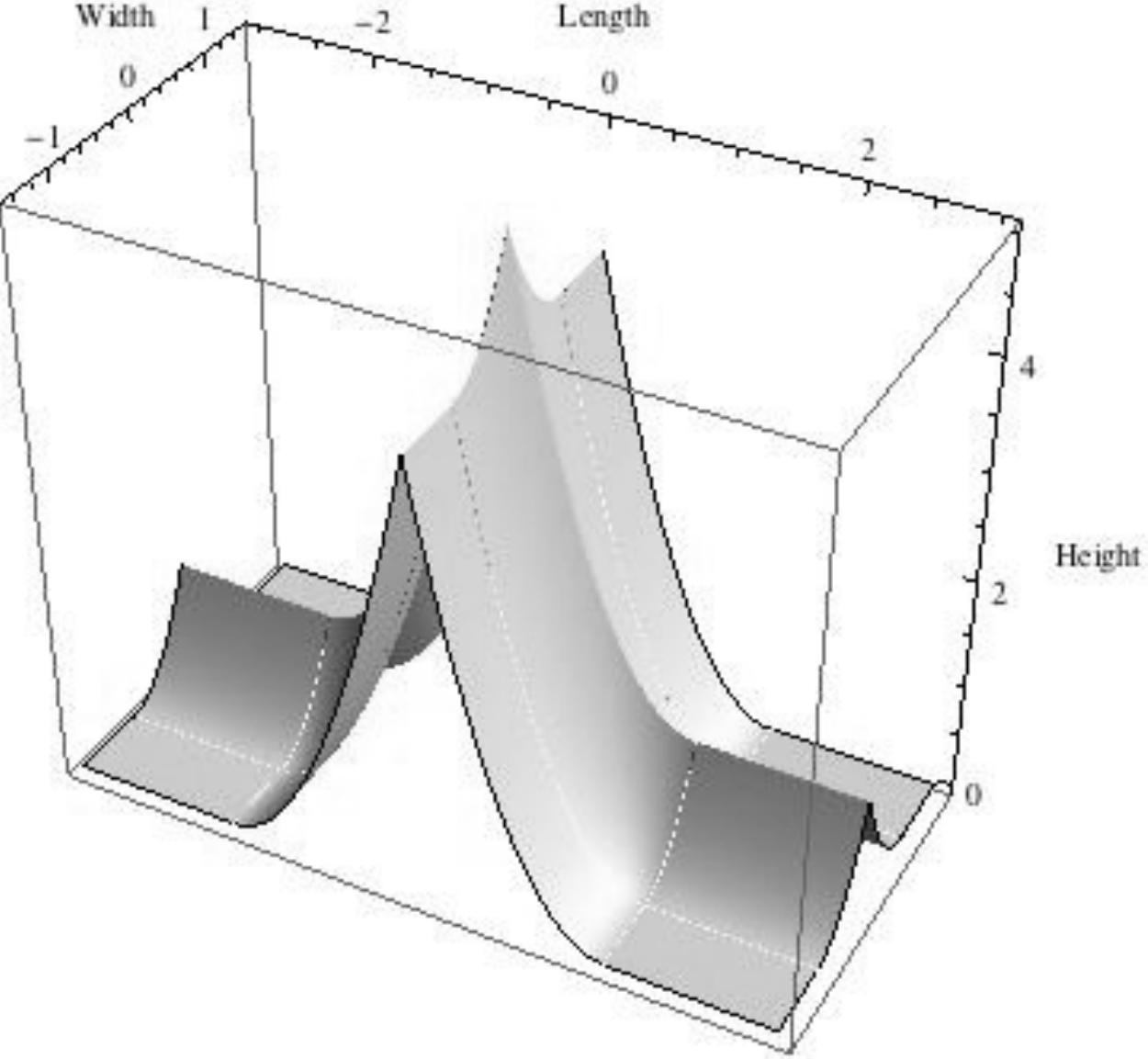}&
		\includegraphics[height=0.35\textwidth]{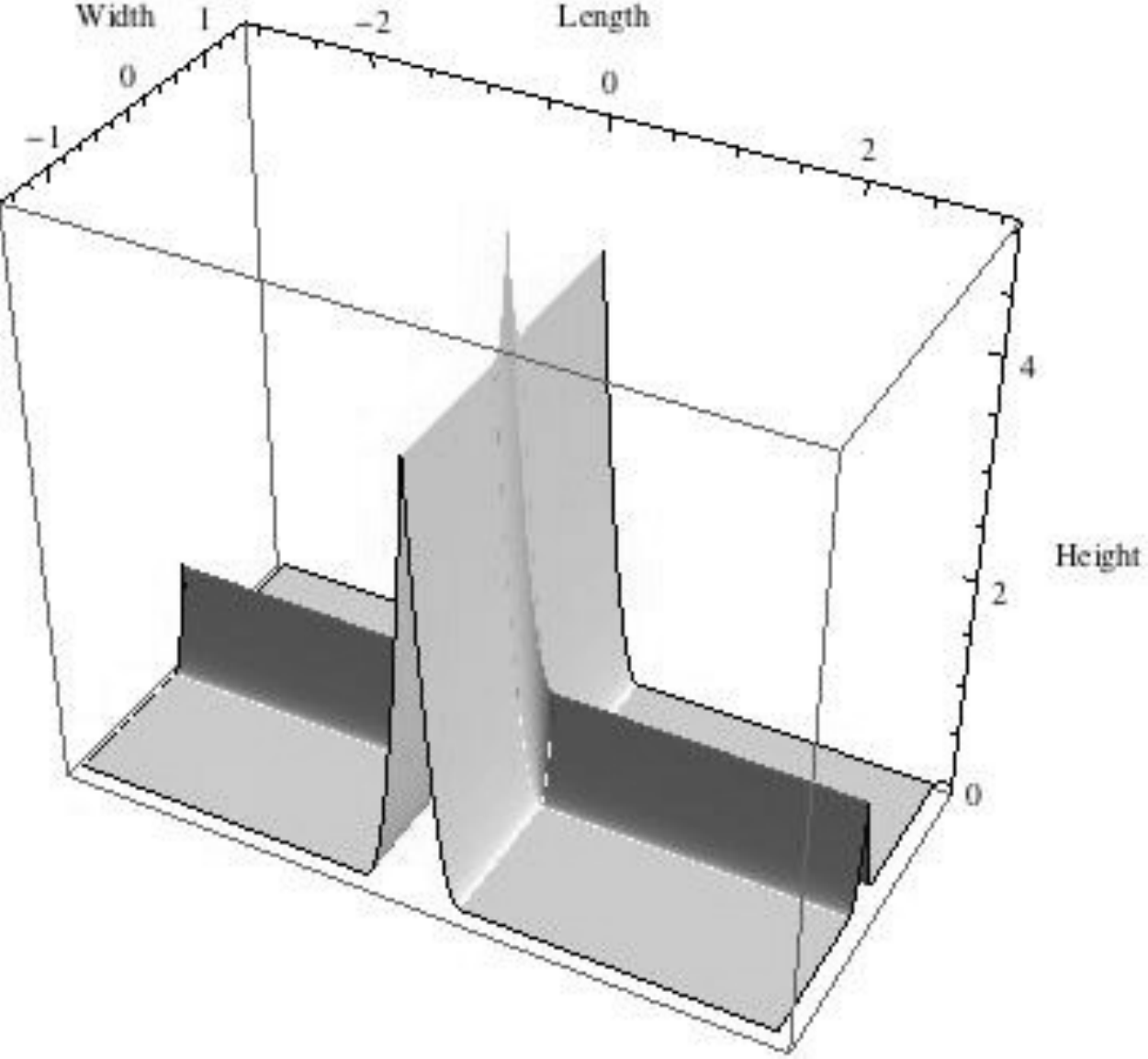}\\
		(a) & (b)\\
		\includegraphics[height=0.35\textwidth]{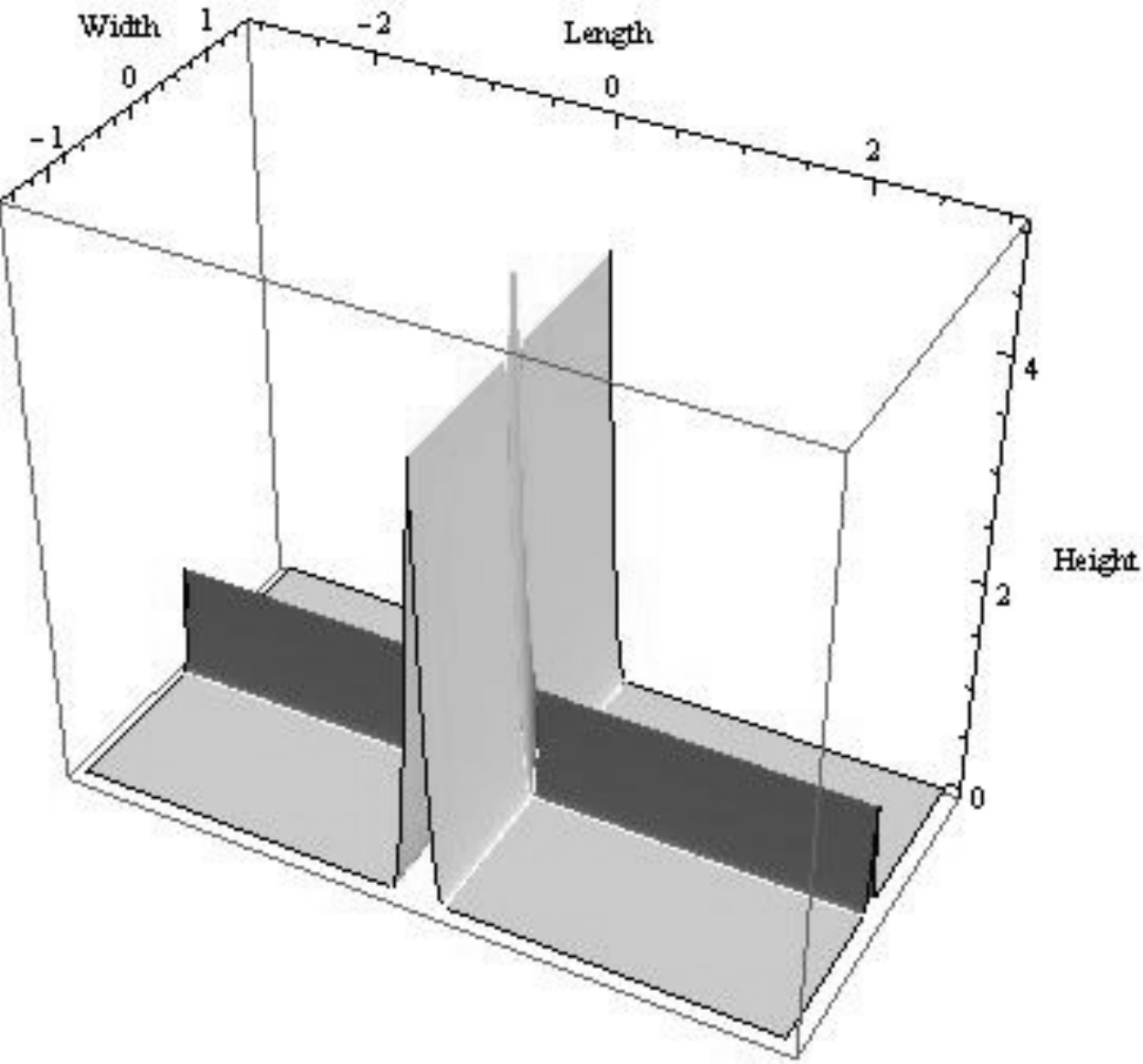}&
		\includegraphics[height=0.35\textwidth]{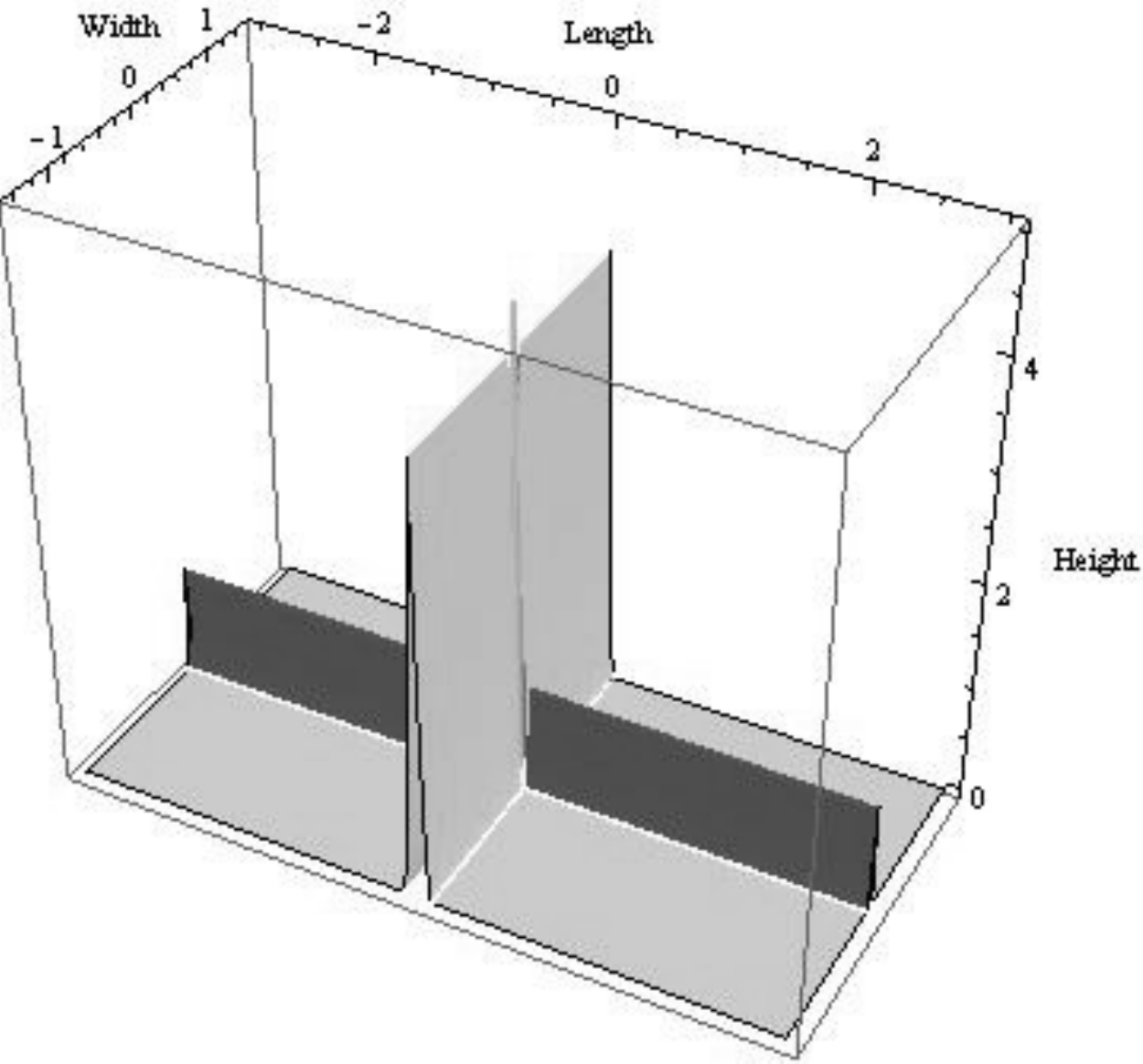}\\
		(c)& (d)
	\end{array}$
}
\caption{\label{Sec5.Fig4pt} 
	Graph of $M_{\lambda}(\cdot;\,K)$ for the four-point set $K=\{(b,c),\,(b,-c),\,(-b,c),\,(-b,-c)\}$ 
	with $b=2$ and $c=\epsilon b=1$, for different values
	of $\lambda$: $(a)$ $\lambda=0.25$; $(b)$ $\lambda=4$; $(c)$ $\lambda=8$; $(d)$ $\lambda=15$.
} 
\end{figure}
 
\end{ex}

%%%%%%%%%%%%%%%%%%%%%%%%%%%%%%%%%%%

\begin{ex}\label{Sec5.Ex2}
In this example, we consider first the case of the open set $\Omega_s=(-r, \infty)\times (-r,\,r)$ with $r>0$, whose results will be
used to construct the multiscale medial axis map of a rectangular domain. 
By inspection, we can easily infer that
\begin{equation}\label{Ex02.DstFnctStrip}
	\dist^2((x,y);\,\partial \Omega_s)=\left\{
			\begin{array}{ll} 
		\displaystyle	\left(|y|-r\right)^2	&
					\displaystyle \text{if }\;x\geq 0;\\[1.5ex]
		\displaystyle	\left(|y|-r\right)^2	&
					\displaystyle \text{if }\;-r\leq x\leq 0, \;\;|y|\geq |x|;\\[1.5ex]
		\displaystyle	\left(|x|-r\right)^2	&
					\displaystyle \text{if }\;-r\leq x\leq 0, \;\;|y|\leq |x|;\\[1.5ex]
		\displaystyle	\left(|x|-r\right)^2	&
					\displaystyle \text{if }\;x\leq -r, \;\;|y|\leq r;\\[1.5ex]
		\displaystyle	\left(|x|-r\right)^2+\left(|y|-r\right)^2	&
					\displaystyle \text{if }\;x\leq -r, \;\;|y|\geq r\,,
			\end{array}
			\right.
\end{equation}
whereas the lower transform, obtained  after lengthy calculations based on the construction of affine functions,
is given, for $(x,y)\in \Omega_s$ by
\begin{equation}\label{Ex02.LwTrStrip}
			\begin{array}{l}
	\displaystyle	C^l_{\lambda}(\dist^2(\cdot;\,\partial \Omega_s))(x,y)=\\[1.5ex]
		\phantom{xxxx}\left\{
			\begin{array}{ll} 
		\displaystyle	g(x+r,y+r)	&
					\displaystyle \text{if }\;-r\leq x\leq 0, \;\; -r\leq y, \;\; x+y\leq -\frac{r}{1+\lambda};\\[1.5ex]
		\displaystyle	g(x+r,r-y)	&
					\displaystyle \text{if }\;-r\leq x\leq 0, \;\; -r\leq -y, \;\; x-y\leq -\frac{r}{1+\lambda};\\[1.5ex]
		\displaystyle	r^2\frac{\lambda}{1+\lambda}-\lambda (x^2+y^2)	&
					\displaystyle \text{if }\; x\leq 0, \;\; -\frac{r}{1+\lambda}\leq x+y, \;\; -\frac{r}{1+\lambda}\leq x-y;\\[1.5ex]
		\displaystyle	r^2\frac{\lambda}{1+\lambda}-\lambda y^2	&
					\displaystyle \text{if }\; x\geq 0, \;\; |y|\leq \frac{r}{1+\lambda};\\[1.5ex]
		\displaystyle	\left(|y|-r\right)^2	&
				\displaystyle \text{if }\; x\geq 0,  \;\; \frac{r}{1+\lambda}\leq |y|;\\[1.5ex]
		\displaystyle	\left(|y|-r\right)^2	&
				\displaystyle \text{if }\; x\geq -r, \;\; |y|\geq r;\\[1.5ex]
		\displaystyle	\left(|x|-r\right)^2	&
				\displaystyle \text{if }\; x\leq -r, \;\; |y|\leq r;\\[1.5ex]
		\displaystyle	\left(|x|-r\right)^2+\left(|y|-r\right)^2	&
				\displaystyle \text{if }\; x\leq -r, \;\; |y|\geq r\,,
			\end{array}\right.
			\end{array}
\end{equation}
where the auxiliary function $g=g(x,y)$ is a continuous piecewise quadratic function defined as follows
\begin{equation}\label{Ex02.AxFnctStrip}
\begin{array}{l}
	\displaystyle	g(x,y)=\\[2ex]
		\phantom{xx}\left\{
		\begin{array}{l} 
			\displaystyle	x^2	\qquad
				\displaystyle \text{if }\;x\leq \frac{\lambda}{1+\lambda},	\;\; x\geq 0, \;\; y\geq 0;\\[1.5ex]
			\displaystyle	y^2	\qquad
				\displaystyle \text{if }\;y\leq \frac{x\lambda}{1+\lambda}, \;\; y\geq 0,\;\; x\geq 0;\\[1.5ex]
			\displaystyle	\frac{\lambda}{1+\lambda}\left((x+y)\frac{1+\lambda}{2\lambda+1}\right)^2-
			\displaystyle	\lambda \left(\left(x-(x+y)\frac{1+\lambda}{2\lambda+1}\right)^2+
				\displaystyle	\left(y-(x+y)\frac{1+\lambda}{2\lambda+1}\right)^2\right)\\[1.5ex]
			\phantom{xxxxxx}	\displaystyle \text{if }\;\frac{x\lambda}{1+\lambda}\leq y \leq \frac{x(1+\lambda)}{\lambda}, 
					\;\; y\geq 0, \;\; x\geq 0;\\[1.5ex]
			\displaystyle	x^2	\qquad
				\displaystyle \text{if }\;x\leq 0, \;\; y\geq 0 ;\\[1.5ex]
			\displaystyle	y^2	\qquad
				\displaystyle \text{if }\;x\geq 0, \;\; y\leq 0 ;\\[1.5ex]
			\displaystyle	x^2+y^2	\qquad
				\displaystyle \text{if }\;x\leq 0, \;\; y\leq 0 \,.
			\end{array}\right.
			\end{array}
\end{equation}

The multiscale medial axis map of $\Omega_s$, $M_{\lambda}((x,y);\,\Omega_s)$, is obtained by applying 
definition \eqref{Eq.Def.MMA} and by taking into account \eqref{Ex02.DstFnctStrip} and \eqref{Ex02.LwTrStrip}.
By exploiting properties of the lower transform with respect to symmetry and 
translation of axis, we can then easily obtain the analytical expression for the 
multiscale medial axis map of a rectangular domain. If, for instance, we consider the open bounded set 
$\displaystyle \Omega=\left(-\left(r+\frac{r}{2}\right),r+\frac{r}{2}\right)\times \left(-r,r\right)$,
then it is not difficult to show that 
\begin{equation}
	M_{\lambda}((x,y);\,\Omega)=\left\{\begin{array}{ll}
				\displaystyle M_{\lambda}((x+\frac{r}{2},y);\,\Omega_s)& \text{if }x\leq 0\,,\\[1.5ex]
				\displaystyle M_{\lambda}((-x+\frac{r}{2},y);\,\Omega_s)& \text{if }x\geq 0\,.
				\end{array}\right.
\end{equation}
Figure \ref{Sec5.FigRct}$(a)$ displays the support of $M_{\lambda}((x,y);\,\Omega)$ which is a neighbourhood of the medial axis,
whereas Figure \ref{Sec5.FigRct}$(b)$ depicts the graph of $M_{\lambda}((x,y);\,\Omega)$. 
For the points $(x,y)\in M_{\Omega}$ with $\theta_x=\pi$,
it follows that $M_{\lambda}((x,y);\,\Omega)=\dist^2((x,y);\,\partial \Omega)$, so implying that
 in this sense, the upper bound in \eqref{Sec3.Eq.OpAng} is 
sharp.

\begin{figure}[htbp]
\centerline{
	$\begin{array}{cc}
		\includegraphics[height=0.35\textwidth]{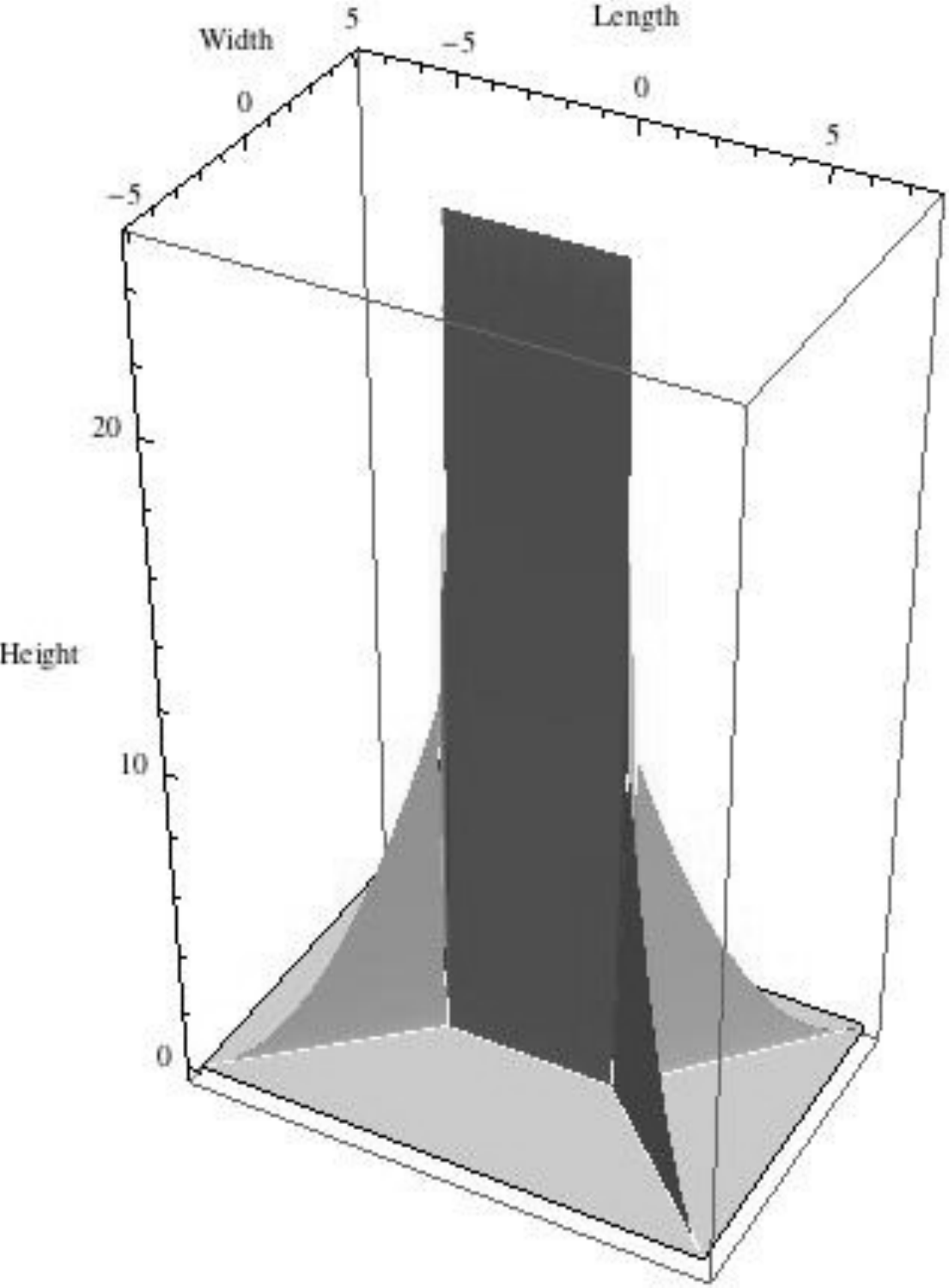}&
		\includegraphics[width=0.30\textwidth]{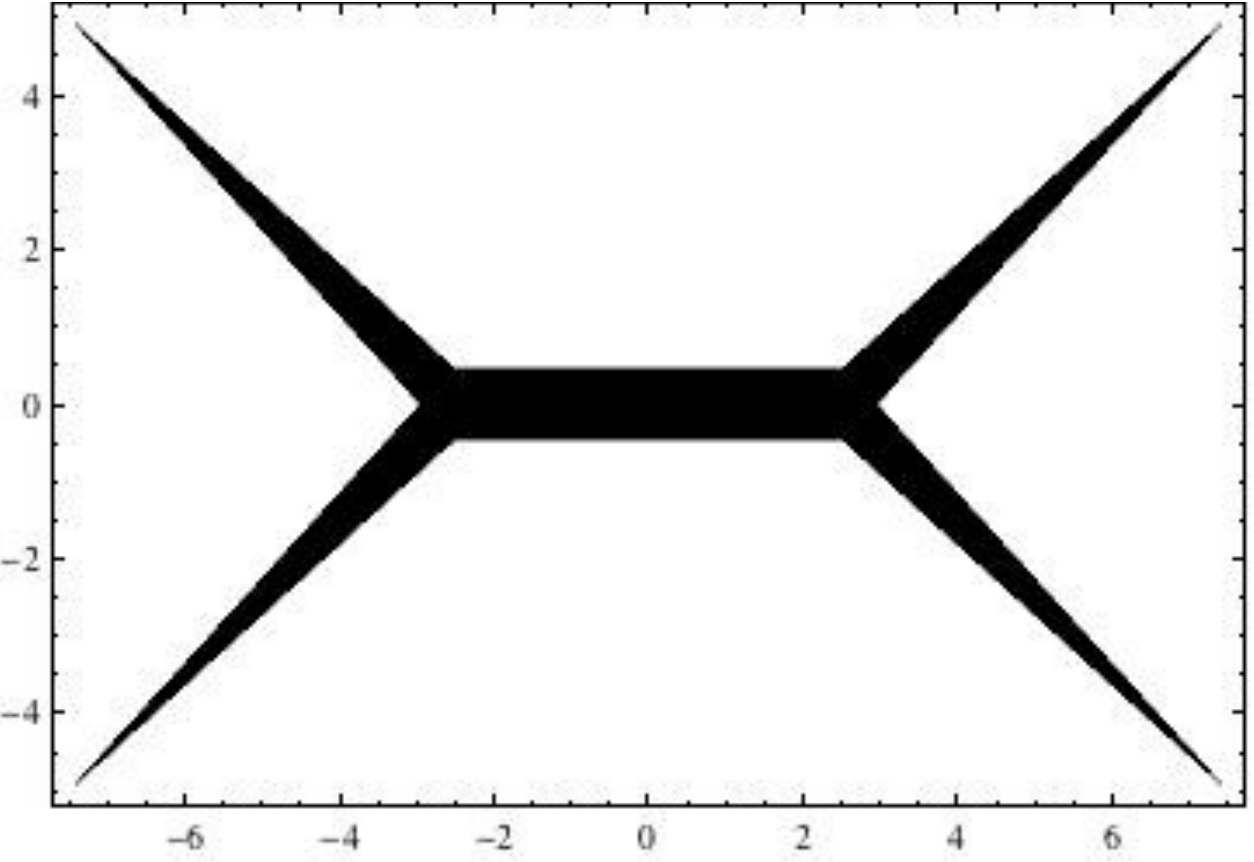}\\
		(a)&(b)
	\end{array}$
}
	\caption{\label{Sec5.FigRct}
	Multiscale medial axis map of $\Omega=\left(-\left(r+\frac{r}{2}\right),r+\frac{r}{2}\right)\times \left(-r,r\right)$ 
	for $\lambda=10$:
	$(a)$ Graph of $M_{\lambda}((x,y);\,\Omega)$ for $\lambda=10$;
	$(b)$ Domain $\Omega$ displayed along with the support of $M_{\lambda}((x,y);\,\Omega)$.}  
\end{figure}

\end{ex}

%%%%%%%%%%%%%%%%%%%%%%%%%%%%%%%%%%%

\begin{ex}\label{Sec5.Ex3}
We consider now the oval shaped domain $\Omega\subset \R^2$ made by the union of two semi-circles with 
center at the points $(-r/2,\,0)$ and $(r/2,\,0)$ and radius $r/2$, respectively, and the rectangle $(-r/2,r/2)\times (-r/2,r/2)$. 
For this domain, it is  not 
difficult to verify that 
\begin{equation}\label{Ex03.DstFnctOval}
	\dist^2((x,y);\,\partial \Omega)=\left\{
			\begin{array}{ll} 
		\displaystyle	\left(|y|-r\right)^2	&
					\displaystyle \text{if }\;|x|\leq \frac{r}{2};\\[1.5ex]
		\displaystyle	\left(\sqrt{\left(x+\frac{r}{2}\right)^2+y^2}-r\right)^2	&
					\displaystyle \text{if }\;x+\frac{r}{2}\leq 0;\\[1.5ex]
		\displaystyle	\left(\sqrt{\left(x-\frac{r}{2}\right)^2+y^2}-r\right)^2	&
					\displaystyle \text{if }\;x-\frac{r}{2}\geq 0
			\end{array}\right.
\end{equation}
whereas the lower transform, obtained  after some lengthy calculations based on the construction of affine functions,
is given, for $(x,y)\in \Omega$, by
\begin{equation}\label{Ex03.LwTrOval}
			\begin{array}{l}
	\displaystyle	C^l_{\lambda}(\dist^2(\cdot;\,\partial \Omega))(x,y)=\\[1.5ex]
		\phantom{xxxx}\left\{
			\begin{array}{ll} 
	\displaystyle	\frac{\lambda r^2}{1+\lambda}-\lambda\left(\left(x+\frac{r}{2}\right)^2+y^2	\right)
		&	\displaystyle \text{if }\; \sqrt{\left(x+\frac{r}{2}\right)^2+y^2}\leq \frac{r}{1+\lambda},\;\; x+\frac{r}{2}\leq 0;\\[1.5ex]
	\displaystyle	\left(\sqrt{\left(x+\frac{r}{2}\right)^2+y^2}-r	\right)^2
		&	\displaystyle \text{if }\; \sqrt{\left(x+\frac{r}{2}\right)^2+y^2}\geq \frac{r}{1+\lambda},\;\; x+\frac{r}{2}\leq 0;\\[1.5ex]
	\displaystyle	\frac{\lambda r^2}{1+\lambda}-\lambda\left(\left(x-\frac{r}{2}\right)^2+y^2	\right)
		&	\displaystyle \text{if }\; \sqrt{\left(x-\frac{r}{2}\right)^2+y^2}\leq \frac{r}{1+\lambda},\;\; x-\frac{r}{2}\geq 0;\\[1.5ex]
	\displaystyle	\left(\sqrt{\left(x-\frac{r}{2}\right)^2+y^2}-r	\right)^2
		&	\displaystyle \text{if }\; \sqrt{\left(x-\frac{r}{2}\right)^2+y^2}\geq \frac{r}{1+\lambda},\;\; x-\frac{r}{2}\geq 0;\\[1.5ex]
	\displaystyle	\frac{\lambda r^2}{1+\lambda}-\lambda y^2
		&	\displaystyle \text{if }\; |y|\leq \frac{r}{1+\lambda},\;\; |x|\leq \frac{r}{2};\\[1.5ex]
	\displaystyle	\left(|y|-r\right)^2
		&	\displaystyle \text{if }\; |y|\geq \frac{r}{1+\lambda},\;\; |x|\leq \frac{r}{2}\,.
\end{array}\right.
			\end{array}
\end{equation}
Figure \ref{Sec5.Fig.MMAMOval}$(a)$ displays the graph of $M_{\lambda}((x,y);\,\Omega)$ obtained by applying definition \eqref{Eq.Def.MMA} 
where we account for \eqref{Ex03.DstFnctOval} and \eqref{Ex03.LwTrOval}, whereas Figure \ref{Sec5.Fig.MMAMOval}$(b)$ shows the support of 
$M_{\lambda}((x,y);\,\Omega)$ along with the domain $\Omega$. 

\begin{figure}[htbp]
\centerline{$\begin{array}{cc}
	\includegraphics[width=0.30\textwidth]{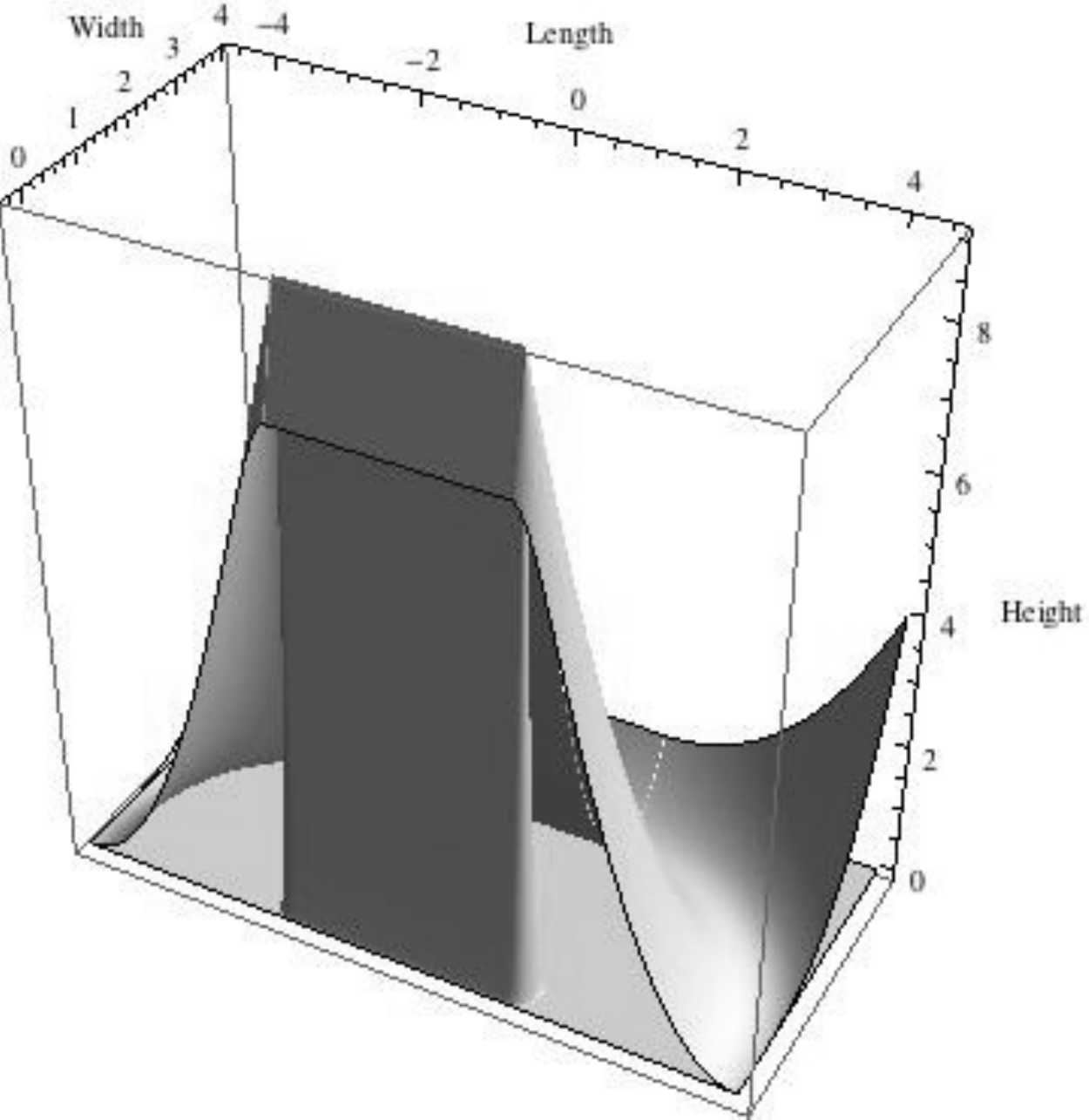}&
	\includegraphics[width=0.30\textwidth]{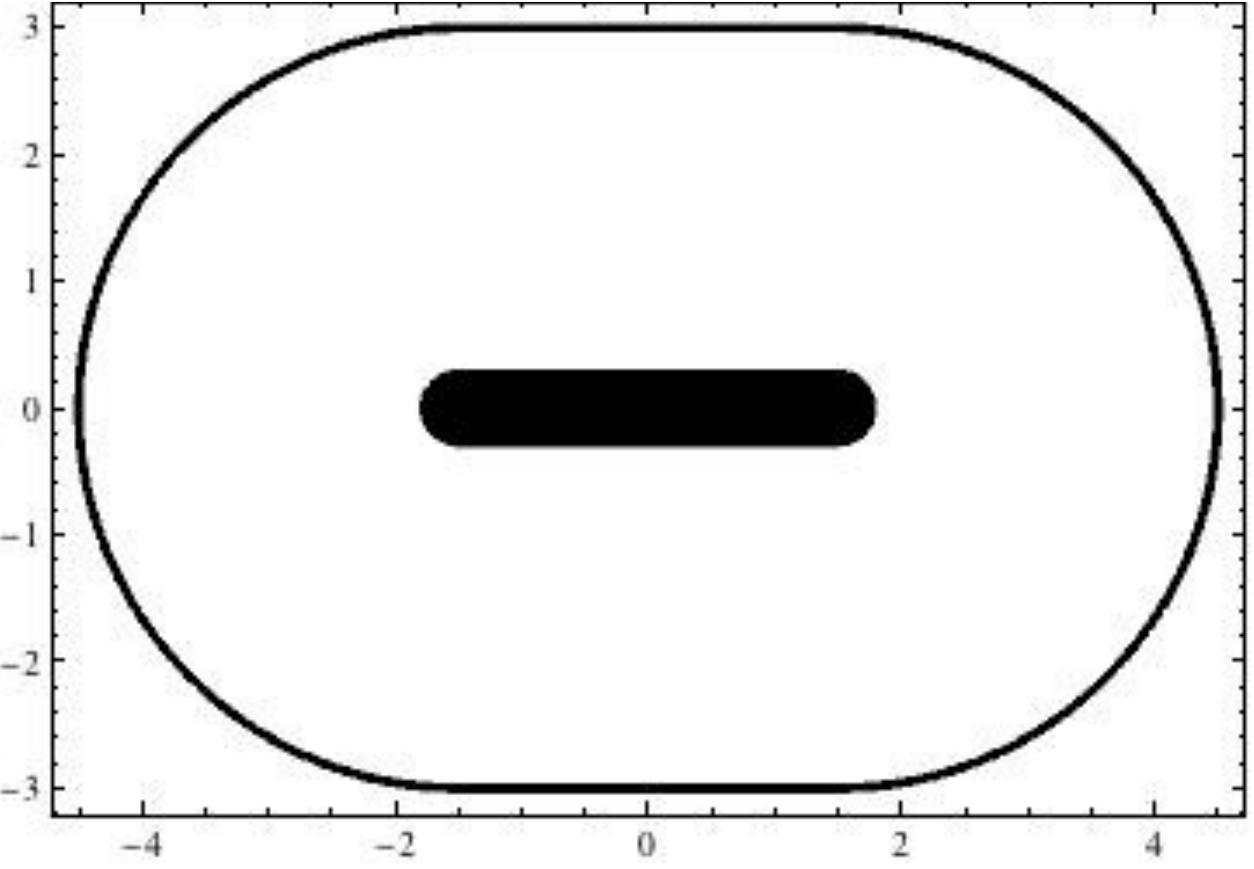}\\
	(a)&(b)
	\end{array}$
	}  
\caption{\label{Sec5.Fig.MMAMOval} 
	Multiscale medial axis map for the oval shaped domain:		
	$(a)$ Graph of the medial axis map, $M_{\lambda}((x,y);\,\Omega)$, for $\lambda=10$ and of the squared distance function,
	$\dist^2((x,y);\,\partial \Omega)$ whose zero level set gives the boundary of the domain. 
	Only the restriction to $y\geq 0$ is displayed, given that  $M_{\lambda}((x,y);\,\Omega)\leq \dist^2((x,y);\,\partial \Omega)$;
	$(b)$ The oval shaped domain $\Omega$ displayed with the support of its medial axis map $M_{\lambda}((x,y);\,\Omega)$. 
}  
\end{figure} 
 
\end{ex}

%%%%%%%%%%%%%%%%%%%%%%%%%%%%%%%%%%%

In the following example we describe the behaviour of $M_{\lambda}(x;\,K)$ for the case of a discrete set $K$, 
sampled from a connected set, and evaluate the structure of $M_{\lambda}(x;\,K)$ as the sample density approaches infinity.

\begin{ex}\label{Sec5.Ex3.a}
We consider the geometric model of the uniform sampling of two parallel lines at distance $b$ to each other.
The points are taken equally spaced over each line at distance $c=\epsilon b$ with $\epsilon\in(0,1)$ measuring the 
sampling density, in the sense that as $\epsilon\to 0$, the sampling density on the two lines tends to infinity. 
The sampling of the two parallel lines is defined so that the discrete points are aligned along the $x-$
and $y-$axis as displayed in Figure \ref{Sec5.Fig.GeoModSmpl}$(a)$. For such a sample $K$, we can then use the results obtained 
in Example \ref{Sec5.Ex01} for the 
four-point set which we refer to as $K_{4p}$. It is  not difficult to show that, for $(x,y)\in \R^2$, and 
$i\in \mathbb{Z}$ such that $|y+i\epsilon b|\leq \epsilon b$, 
\begin{equation}
	M_{\lambda}((x,y);\,K)=M_{\lambda}((x,y+i\epsilon b);\,K_{4p})\,,
\end{equation}
where $M_{\lambda}((x,y+i\epsilon b);\,K_{4p})$ is the multiscale medial axis map for the four-point set discussed in 
the Example \ref{Sec5.Ex01}.

\begin{figure}[htbp]
$\begin{array}{lll}
	\includegraphics[height=0.35\textwidth]{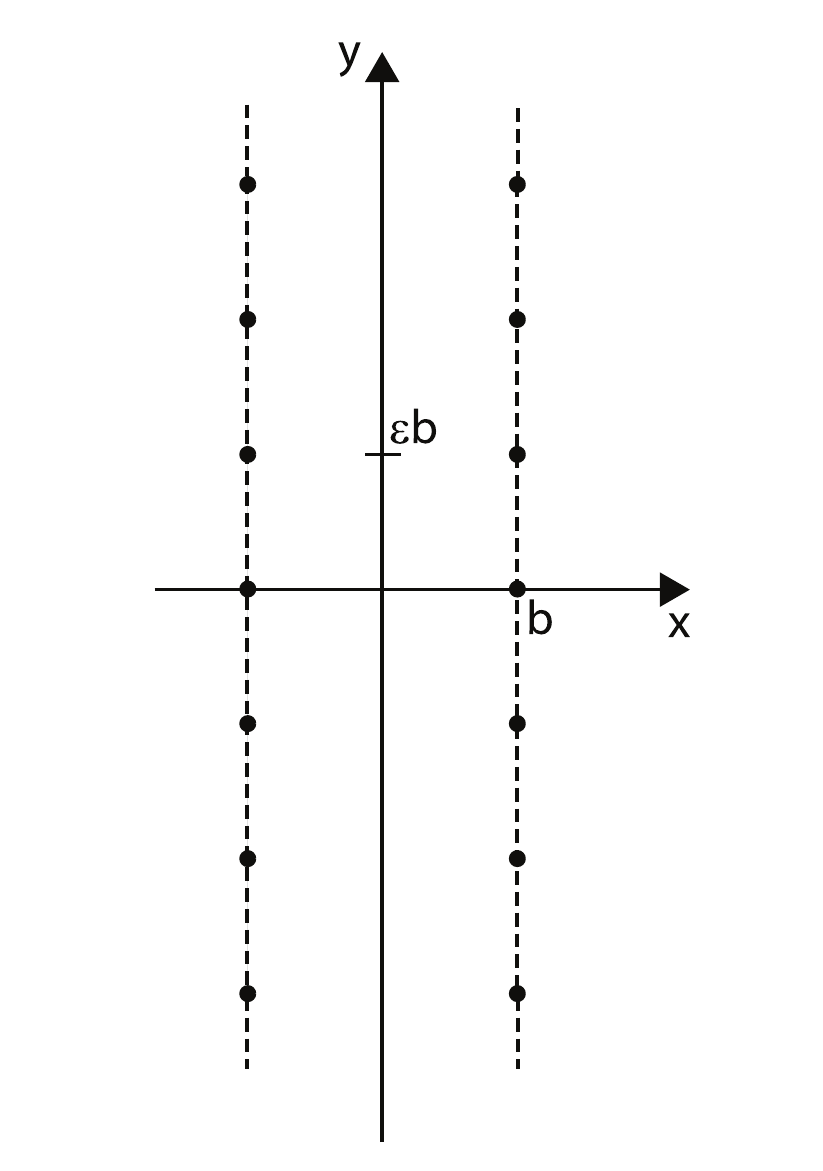}
	&
	\includegraphics[width=0.30\textwidth]{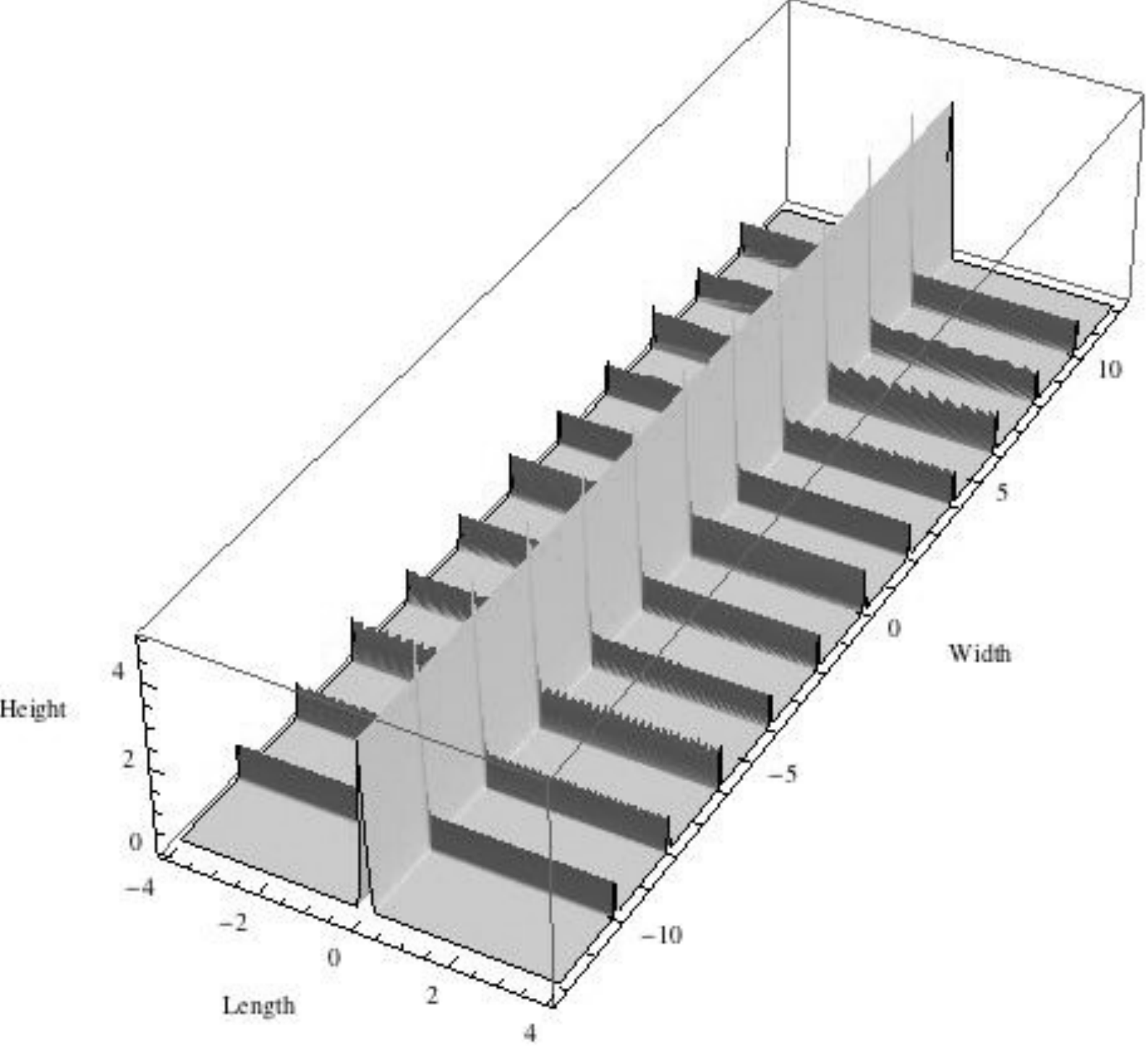}
	&
	\includegraphics[width=0.40\textwidth]{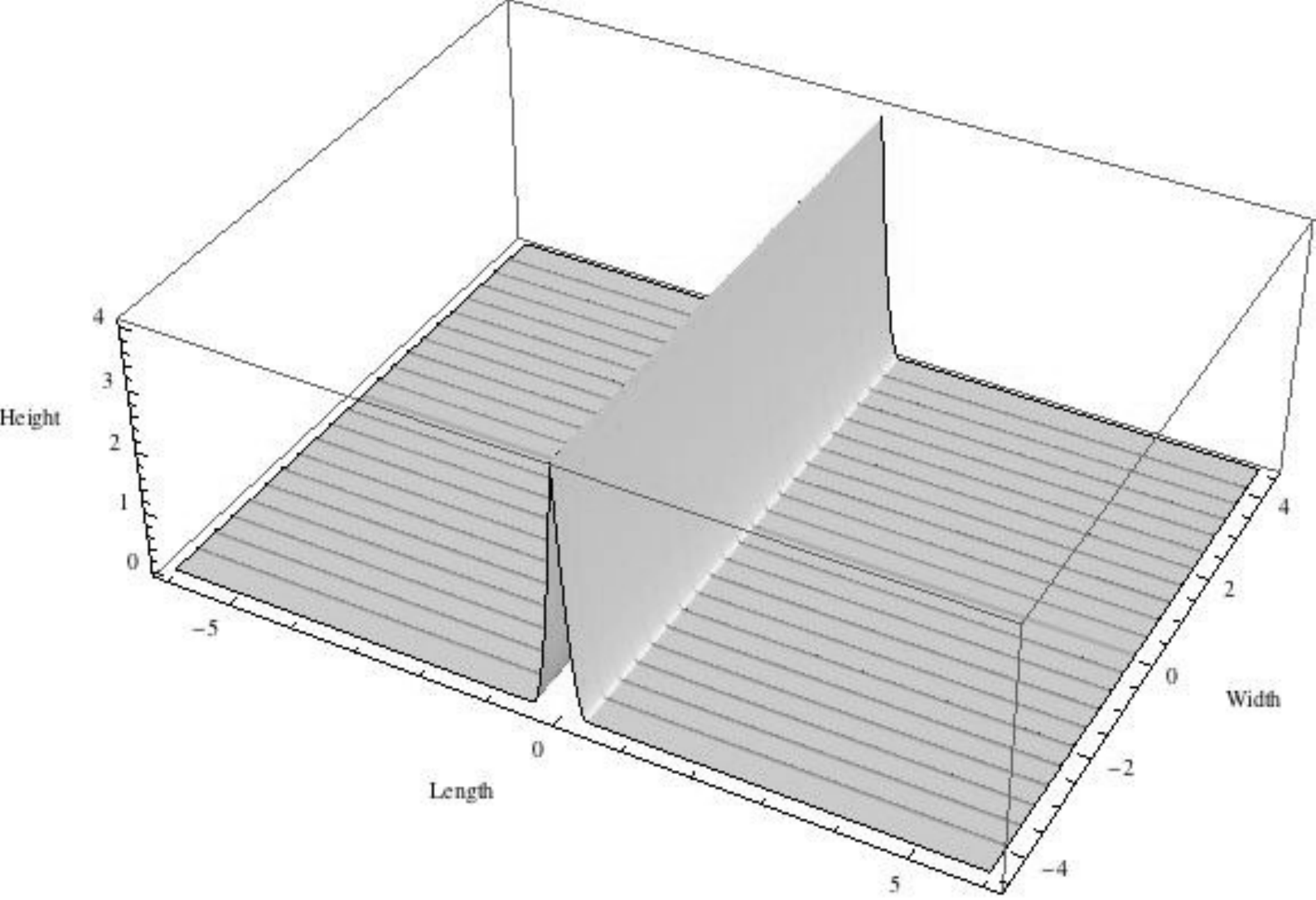}\\
	(a)&(b)&(c)
\end{array}
$
\caption{\label{Sec5.Fig.GeoModSmpl} 
	$(a)$ Geometric model of a set $K$ representing a uniform $\epsilon-$sample of two parallel lines.
	Graph of the multiscale medial axis map for different sample density:
	$(b)$ $b=2$; $\epsilon=0.5$; $\lambda=10$;
	$(c)$ $b=2$; $\epsilon=0.1$; $\lambda=10$.}
\end{figure} 

Figure \ref{Sec5.Fig.GeoModSmpl}$(b)$ and \ref{Sec5.Fig.GeoModSmpl}$(c)$ display the graph of $M_{\lambda}((x,y);\,K)$ 
for $b=2$, $\lambda=10$ and for different different sampling density $\epsilon$, 
$\epsilon=0.5$ and $\epsilon=0.1$, respectively.

The comparison of the graphs of $M_{\lambda}((x,y);\,K)$ for the two different sample densities of the same object, shows how 
the value of $M_{\lambda}((x,y);\,K)$ along the minor branches of the medial axis of the discrete set $K$ attenuates as
$\epsilon \to 0$. Consistently  with the finding obtained for the four-point set,
we have that the value of $M_{\lambda}((x,y);\,K)$ along the minor branches is proportional to $\epsilon^2$. It then 
follows that by setting a threshold not lower than such a value, we can single out the stable part of $M_K$, which
provides an approximation for (and in this case, is in fact coincident with) the medial axis of the two parallel lines.

\end{ex}

%%%%%%%%%%%%%%%%%%%%%%%%%%%%%%%%%%%

In this last example, we analyze a model of perturbations of the boundary domain represented by staircase-like
piecewise affine curves. This effect is very common, for instance, in digital images and is the source
of unrealistic medial axis branches, which are usually not desirable. Common practice in this case is to perform 
a boundary smoothing prior to any image processing operation. We will show that this is not needed with the multiscale medial axis map.
The fine structure of the medial axis corresponding to the irregularities of the boundary is indeed captured by the multiscale medial 
axis map and can be filtered out.
We will verify this statement on a prototype model of this boundary domain perturbation, by
showing that the height of the medial axis map on such branches can be very small if the 
stair like effect is small. 

%%%%%%%%%%%%%%%%%%%%%%%%%%%%%%%%%%%

\begin{ex}\label{Sec5.Ex4}
Assume $c>0$ and let us consider the set $K=\{(x,y)\in\mathbb{R}^2:\, x\geq 0,\, y\geq 0,\, x+y\leq c\}\cap
\{(x,y)\in\mathbb{R}^2,\,  x+y\geq c\}$ displayed in Figure \ref{Sec5.Fig.GeoModStep}, which is  used 
as a prototype of one single step perturbation of a boundary domain. 

\begin{figure}[htbp]
\centerline{\includegraphics[height=0.30\textwidth]{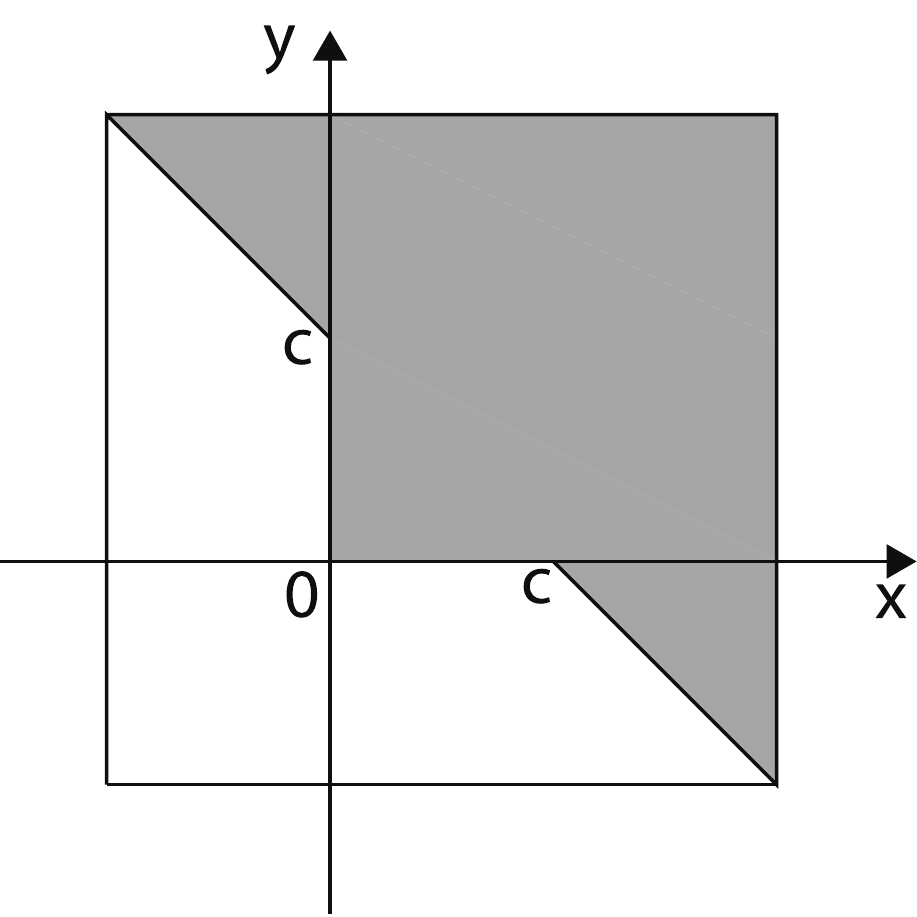}
	}  
\caption{\label{Sec5.Fig.GeoModStep} 
	Geometric model of one stair-like boundary perturbation, as in a
	digitized domain. 
}
\end{figure} 
 
The squared distance function $\dist^2((x,y);\, K^c))$ to the complement of $K$  and  the 
lower transform $C^l_\lambda(\dist^2(\cdot;\, K^c))(x,y)$ are then given, respectively, by
\begin{equation}\label{Ex.DstFncOneStep}
	\dist^2((x,y);\, K^c)=\left\{
		\begin{array}{ll} 
\displaystyle		x^2		&\text{if } x \leq y ,\quad 0 \leq x \leq c,\quad 0 \leq y \leq c,\\[1.5ex]
\displaystyle		y^2		&\text{if }  x \geq y, \quad 0 \leq x \leq c,\quad 0 \leq  y \leq c,\\[1.5ex]
\displaystyle		x^2 + (y-c)^2	&\text{if } y \geq c, \quad 0 \leq y - x \leq c,\\[1.5ex]
\displaystyle		(x-c)^2 + y^2	&\text{if } x \geq c,\quad 0 \leq x - y \leq c,\\[1.5ex]
\displaystyle		\frac{1}{2}(x + y-c)^2	
						&\text{if } x + y \geq c,     \quad y - x \geq c,\\[1.5ex]
\displaystyle		\mbox{}			&\text{or } x + y \geq c,\quad  x - y \geq c,\\[1.5ex]
\displaystyle		0			&\text{otherwise}\,;
		\end{array}\right.
\end{equation}
and
\begin{equation}\label{Ex.LwTrOneStep}
\begin{array}{l}	
     \displaystyle C^l_\lambda(\dist^2(\cdot;\, K^c))(x,y)=\\
	\phantom{xx}\left\{\begin{array}{l}
\displaystyle 		x^2	\qquad\text{if }  x \leq  \frac{y\lambda}{1 + \lambda},\quad  0 \leq x \leq c, \quad 0 \leq y \leq c,\\[1.5ex]
\displaystyle 		y^2	\qquad\text{if }  y \leq  \frac{x\lambda}{1 + \lambda},\quad 0 \leq x \leq c,\quad 0 \leq y \leq c,\\[1.5ex]
\displaystyle 		\lambda (1 + \lambda) \left(\frac{x + y}{1 + 2 \lambda}\right)^2 - 
\displaystyle 			\lambda \left(x - \frac{(1 + \lambda) (x + y)}{1 + 2 \lambda}\right)^2 - 
\displaystyle 			 \lambda \left(y - \frac{(1 + \lambda)(x + y)}{1 + 2 \lambda}\right)^2\\[1.5ex]
\displaystyle 		\qquad\mbox{}\quad\text{if }  \frac{\lambda x}{1 + \lambda} \leq 
\displaystyle 			y \leq \frac{(1 + \lambda) x}{\lambda},\quad 0 \leq x \leq c,\quad 0 \leq y \leq c,\quad
\displaystyle 			x + y \leq \frac{(1 + 2 \lambda)c}{1 + \lambda},\\[1.5ex]
\displaystyle 		\frac{1}{2}(x + y-c)^2; \qquad\text{if } x + y \geq c,\quad y - x \geq c, \quad\text{or}\quad  x + y \geq c,
\displaystyle 			\quad x - y \geq c,\\[1.5ex]
\displaystyle 		x^2 + (y - c)^2 \qquad\text{if } y \geq c,\quad \frac{c}{1 + \lambda} \leq y - x \leq c,\\[1.5ex]
\displaystyle 		(x - c)^2 + y^2 \qquad\text{if } x \geq c,\quad \frac{c}{1 + \lambda} \leq x - y \leq c,\\[1.5ex]
\displaystyle 		\frac{1}{2}\left(\frac{\lambda c^2}{1 + \lambda} - \lambda (y - x)^2 + (x + y - c)^2\right)\\[1.5ex]
\displaystyle 		\phantom{xxxxxxxx}\text{if }  x + y \geq (1 + 2 \lambda) \frac{c}{1 + \lambda},\quad  
\displaystyle 			\frac{c}{1 + \lambda} \geq y - x \geq -\frac{c}{1 + \lambda},\\[1.5ex]
\displaystyle 		x^2 + y^2 \qquad\text{if } x \leq 0,\quad y \leq 0,\\[1.5ex]
\displaystyle 		0 \qquad \text{otherwise}\,.
	\end{array}\right.
\end{array}
\end{equation}
The medial axis map $M_{\lambda}((x,y);\,K^c)$ is then obtained from \eqref{Eq.Def.MMA} using \eqref{Ex.DstFncOneStep} and 
\eqref{Ex.LwTrOneStep}. The graph of  $M_{\lambda}((x,y);\,K^c)$ for $\lambda=9$ and step size $c=1$ is shown in  
Figure \ref{Sec5.Fig.SingStp}$(a)$, whereas Figure \ref{Sec5.Fig.SingStp}$(b)$ displays its support.
By inspecting the graph of $M_{\lambda}((x,y);\,K^c)$, we observe that after an initial increase near the corner tip, 
$M_{\lambda}((x,y);\,K^c)$ keeps a constant value along $M_K$, with this value proportional to the square of the step size.
It is not difficult to verify the following
\begin{equation}
	\lim_{\lambda\to+\infty}M_{\lambda}((x,y);\,K)=
		M_{\infty}((x,y);\,K)
		=
		\left\{\begin{array}{ll}
		\displaystyle \frac{x^2}{2}\,,	& \displaystyle  \text{if }0\leq x\leq c,\quad x=y	\,,\\[1.5ex]
		\displaystyle \frac{c^2}{2}\,,	& \displaystyle  \text{if }x\geq c ,	 \quad x=y	\,,\\[1.5ex]
		\displaystyle 0\,,		& \displaystyle  \text{if }x\not = y\,,
	\end{array}\right.
\end{equation}
with uniform convergence if $y=x$ and $x\geq 0$.\\

\begin{figure}[htbp]
\centerline{$\begin{array}{cc}
	\includegraphics[width=0.35\textwidth]{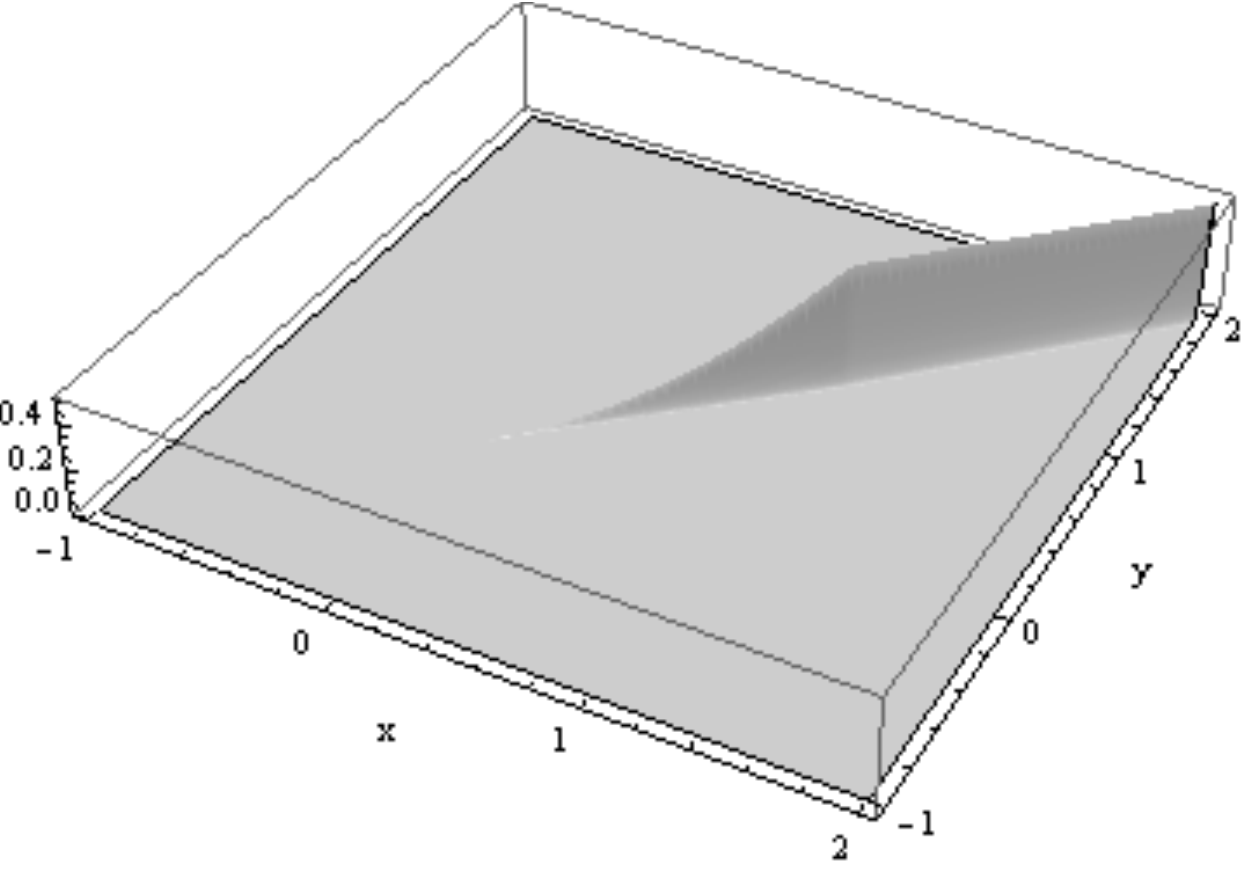}&
	\includegraphics[width=0.30\textwidth]{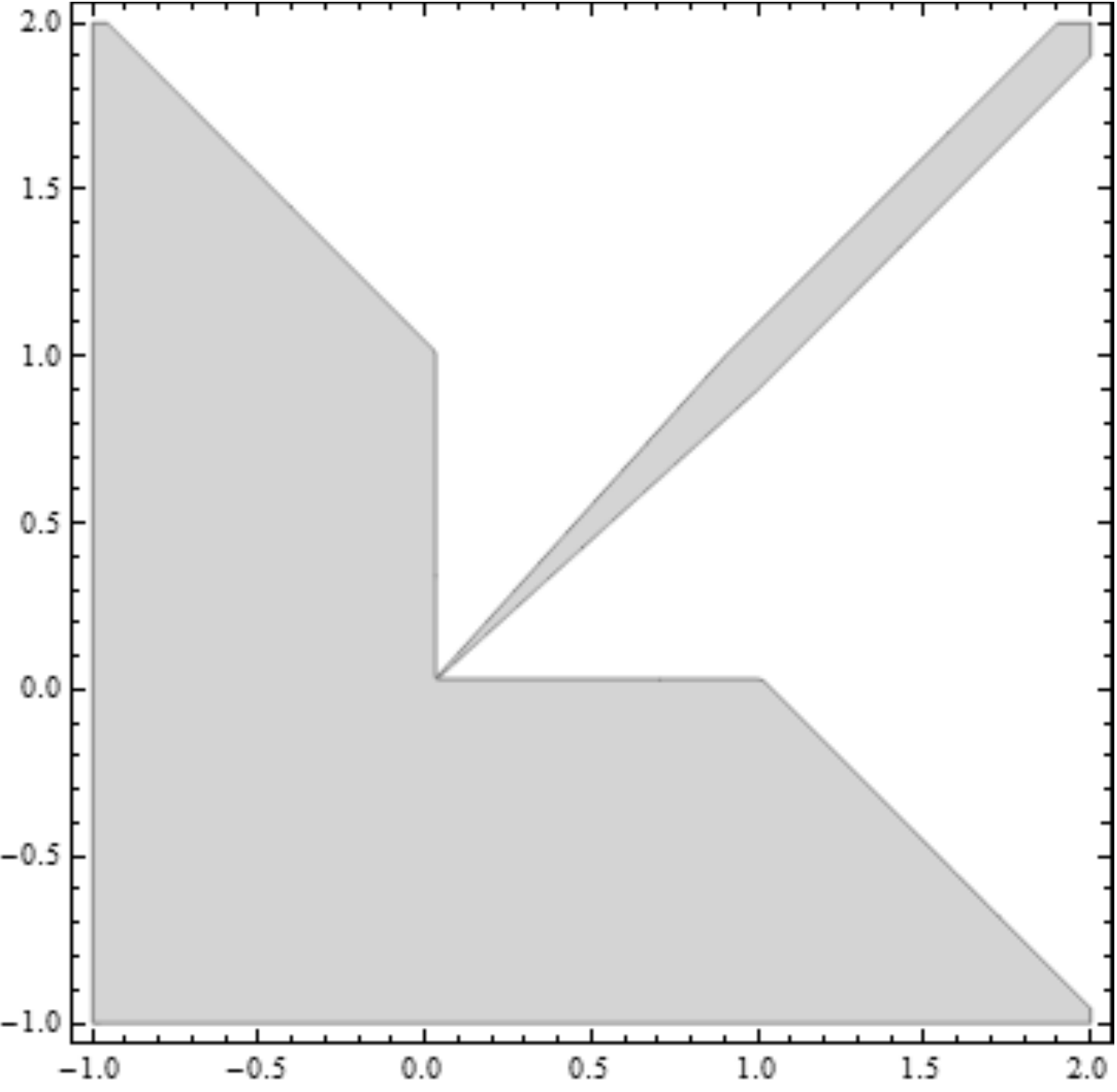}\\
	(a)&(b)
	\end{array}
	$}
\caption{\label{Sec5.Fig.SingStp}
	Multiscale medial axis map of the set $K$ displayed in Figure \ref{Sec5.Fig.GeoModStep}, prototype of boundary perturbation as
	occurring in digitized images.
	$(a)$ Graph of  $M_{\lambda}((x,y);\,K^c)$ for $\lambda=10$; 
	$(b)$ Support of $M_{\lambda}((x,y);\,K^c)$ displayed together with the set $K$.}  
\end{figure}

Despite its simplicity, this basic model elucidates the behaviour of $M_{\lambda}((x,y);\,K_s^c)$ for a set $K_s$ 
with a stair like boundary profile as, for instance, the one displayed in Figure \ref{Sec5.Fig.GeoModStrc}$(a)$.

\begin{figure}[htbp]
\centerline{$\begin{array}{ccc}
	\includegraphics[height=0.30\textwidth]{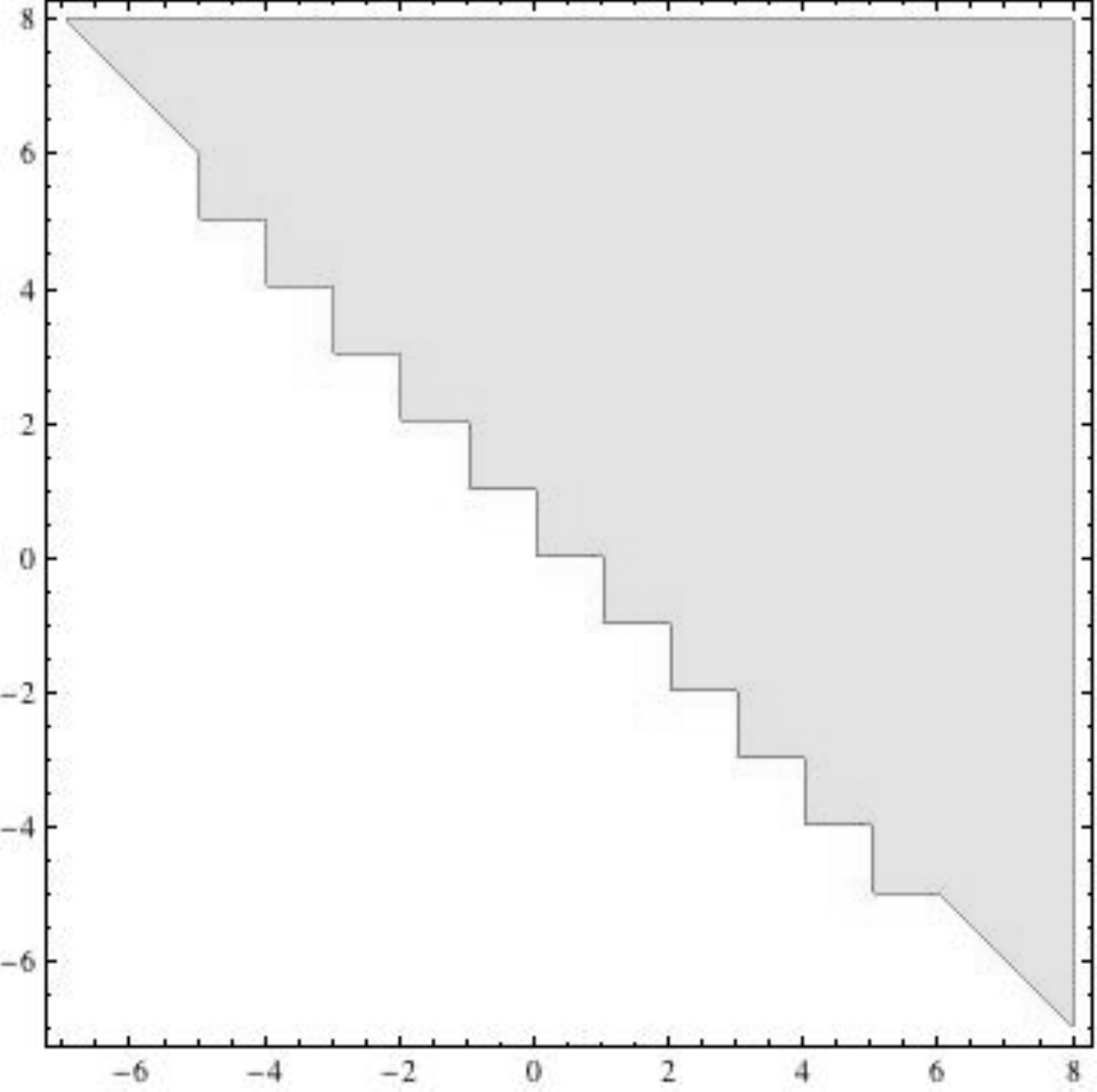}
	&
	\includegraphics[height=0.30\textwidth]{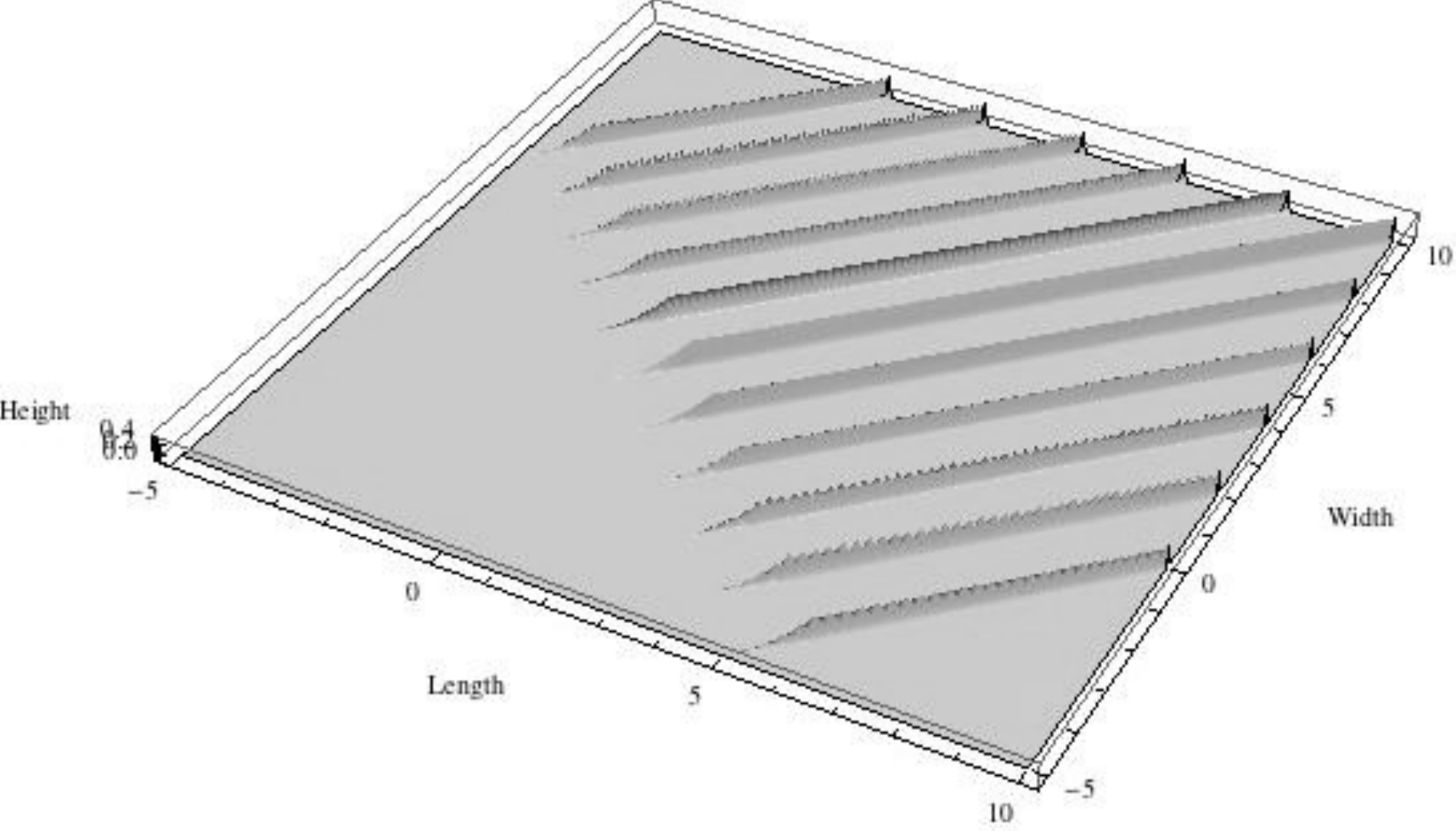}
	&
	\includegraphics[height=0.30\textwidth]{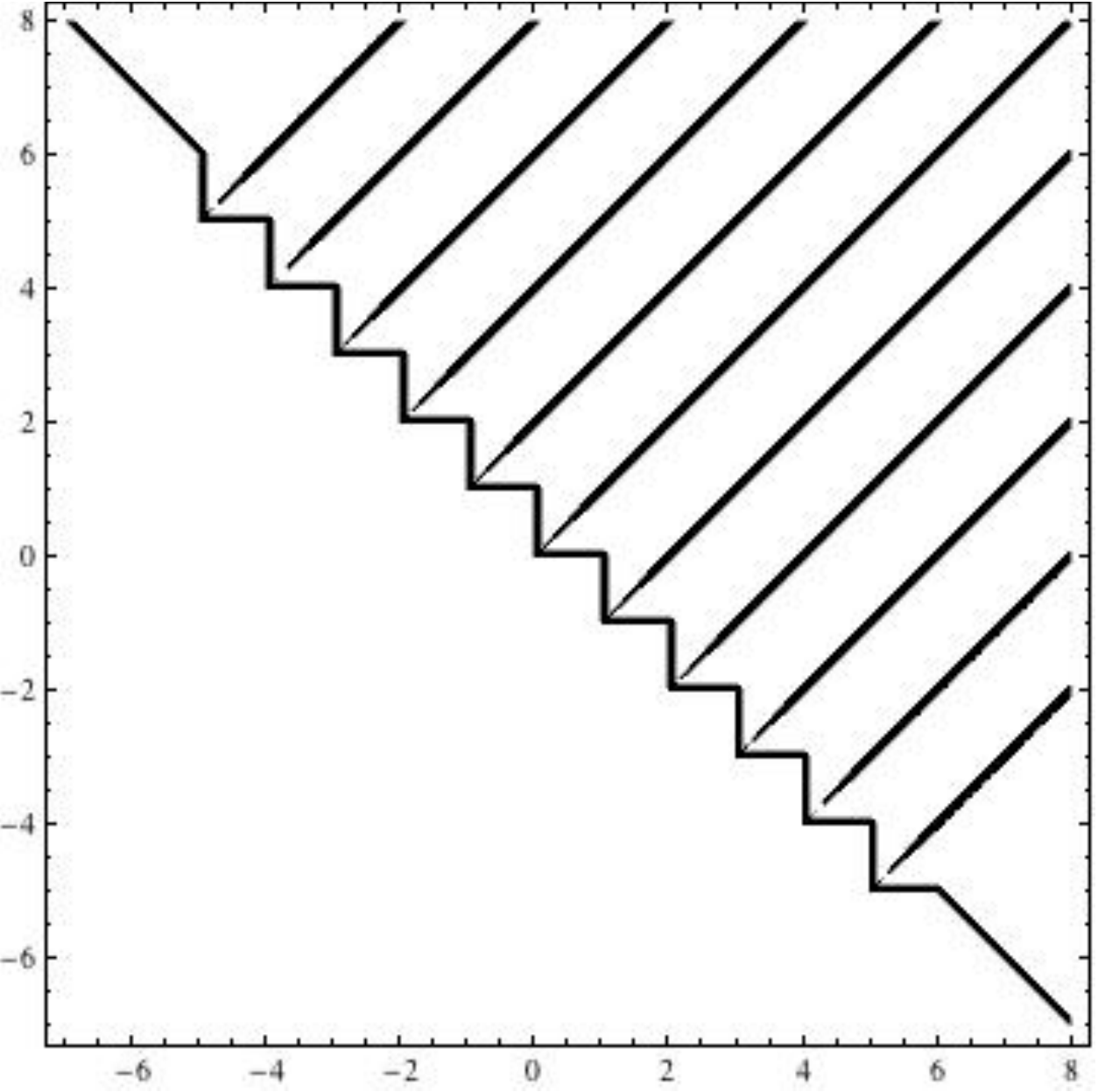}\\
	(a)&(b)&(c)
\end{array}$}
\caption{\label{Sec5.Fig.GeoModStrc} 
	$(a)$ Geometric model of a prototype of staircase type perturbation of the boundary domain $K_s$.
	$(b)$ Graph of  $M_{\lambda}((x,y);\,K_s^c)$ for $\lambda=10$;  
	$(c)$ Support of $M_{\lambda}((x,y);\,K_s^c)$ displayed together with the set $K$.
	}
\end{figure} 

\begin{figure}[htbp]
	\centerline{$\begin{array}{ccc}
		\includegraphics[height=0.35\textwidth]{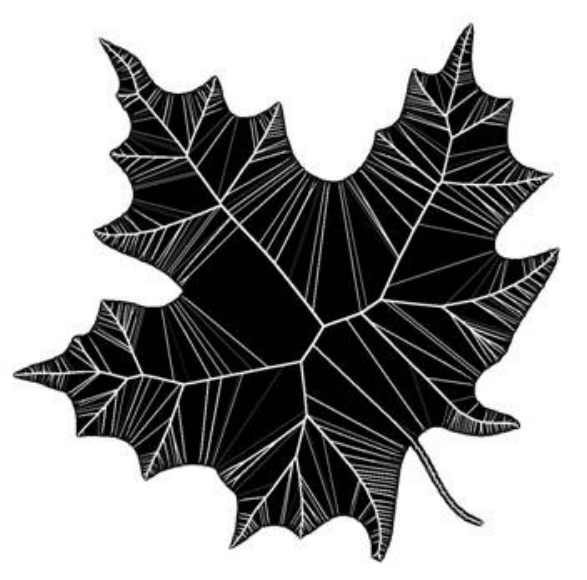}&\phantom{ccccccc} &
		\includegraphics[height=0.35\textwidth]{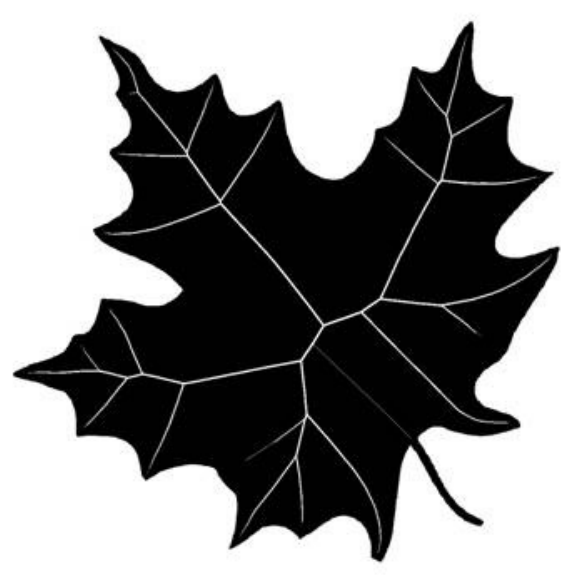}\\
		(a)&&(b)
		\end{array}$
		}
\caption{\label{Sec5.Fig.MapleLeaf}
	Multiscale medial axis map as obtained by the numerical implementation of $M_{\lambda}((x,y);\,K)$: 
	$(a)$ Support of the multiscale medial axis map $M_{\lambda}((x,y);\,K)$ for the digital image of 
	a maple leaf, for $\lambda=10$. All the fine branches generated by the steps on the boundary are displayed.; 
	$(b)$ Suplevel set of $M_{\lambda}((x,y);\,K)$ corresponding to a threshold equal to one
	displaying only stable parts of the medial axis.}  
\end{figure}  

For such a set $K_s$, it is not difficult to verify that for $(x,y)\in \R^2$, let $i\in \mathbb{Z}$ such that 
$|(x-ic)-(y+ic)|\leq c$, then 
\begin{equation}
	M_{\lambda}((x,y);\,K_s^c)=M_{\lambda}((x-ic,y+ic);\,K^c)
\end{equation}
with $M_{\lambda}((x,y);\,K^c)$ corresponding to the one step boundary domain perturbation discussed at the beginning of this example.
Figure \ref{Sec5.Fig.GeoModStrc}$(b)$ contains the graph of $M_{\lambda}((x,y);\,K_s)$, whereas 
Figure \ref{Sec5.Fig.GeoModStrc}$(c)$ shows its support, displayed together 
with the set $K_s$. The height of the ridges along $M_{K_s}$ depends only on the gap size $c$, in particular, it is 
proportional to $c^2$. It follows, therefore, that by setting the threshold larger than $c^2$, 
the corresponding suplevel set of $M_{\lambda}((x,y);\,K_s)$ will filter out all  minor branches of $M_{K_s}$, generated by the 
step-stair like  boundary. 

As an application of these concepts, we show in Figure \ref{Sec5.Fig.MapleLeaf} 
the results of the numerical realization of $M_{\lambda}((x,y);\,K)$
for the digital image of a maple leaf, where we can note the effects just discussed. In particular, Figure \ref{Sec5.Fig.MapleLeaf}$(a)$ depicts
the support of $M_{\lambda}((x,y);\,K)$ with the display of all  fine branches created by the step-like irregularities of the boundary domain, 
whereas Figure \ref{Sec5.Fig.MapleLeaf}$(b)$ shows the suplevel set of $M_{\lambda}((x,y);\,K)$ corresponding to a threshold  
equal to one which singles out only the neighbourhood of  stable parts of $M_K$.

\end{ex}

%%%%%%%%%%%%%%%%%%%%%%%%%%%%%%%%%%%%%%%%%%%%%%%%%%%%%%%%%%%%%%%%%%%%%%%%%%%%%%%%%%%%%%%
%%%%%%%%%%%%%%%%%%%%%%%%%%%%%%%%%%%%%%%%%%%%%%%%%%%%%%%%%%%%%%%%%%%%%%%%%%%%%%%%%%%%%%%
%%%%%%%%%%%%%%%%%%%%%%%%%%%%%%%%%%%%%%%%%%%%%%%%%%%%%%%%%%%%%%%%%%%%%%%%%%%%%%%%%%%%%%%

\section{Proofs of Main Results}\label{SecPrf}

%%%%%%%%%%%%%%%%%%%%%%%%%%%%%%%%%%%%%%%%%%%%%%%%%%%%%%%%%%%%%%%%%%%%%%%%%%%
%%%%%%%%%%%%%%%%%%%%%%%%%%%%%%%%%%%%%%%%%%%%%%%%%%%%%%%%%%%%%%%%%%%%%%%%%%%
%%%%%%%%%%%%%%%%%%%%%%%%%%%%%%%%%%%%%%%%%%%%%%%%%%%%%%%%%%%%%%%%%%%%%%%%%%%

\noindent {\bf Proof of Lemma \ref{Sec2.Lem.UpDif}:}
The existence of an affine function $\ell$ such that \eqref{Sec2.Lem.UpDif.Eq}$(i)$ and 
\eqref{Sec2.Lem.UpDif.Eq}$(ii)$ hold is  well known  (see e.g. \cite[Remark 2.1]{Zha08a} or \cite{HUL01}). 
The claim that $2\leq k\leq n+1$ comes from the 
Carath\'eodory Theorem \cite[Cor. 17.1.5 ]{Roc66} and the fact that $\co[f](0)<f(0)$. 
Also it is easy to see that $x_i$'s can be made distinct and $x_i\neq 0$. 
For the proof of \eqref{Sec2.Lem.UpDif.Eq}$(iii)$, observe that $f$ is upper 
semidifferentiable, $\ell\leq f$ and $\ell(x_i)=f(x_i)$. 
By \cite[Cor. 2.5]{BKK00}, it thus follows  that $f$ is differentiable at $x_i$ and 
$Df(x_i)=D\ell(x_i)=2a$.
The proof of \eqref{Sec2.Lem.UpDif.Eq}$(iv)$ is obtained from the definition of the convex envelope.
For \eqref{Sec2.Lem.UpDif.Eq}$(v)$, we have, by definition of the convex envelope, 
that $\ell\leq \co[f]$ and $\ell(0)=\co[f](0)$. Since by hypothesis, $\co[f]\in C^{1,1}(\R^n)$, 
 we can conclude that $2a=D\ell(0)=D\co[f](0)$, again using  \cite[Cor. 2.5]{BKK00}.
\hfill\qed\\

%%%%%%%%%%%%%%%%%%%%%%%%%%%%%%%%%%%%%%%%%%%%%%%%%%%%%%%%%%%%%%%%%%%%%%%%%%%
%%%%%%%%%%%%%%%%%%%%%%%%%%%%%%%%%%%%%%%%%%%%%%%%%%%%%%%%%%%%%%%%%%%%%%%%%%%
%%%%%%%%%%%%%%%%%%%%%%%%%%%%%%%%%%%%%%%%%%%%%%%%%%%%%%%%%%%%%%%%%%%%%%%%%%%

\noindent {\bf Proof of Proposition \ref{Sec2.Th.RegLwTr}:}
It is known \cite[Theorem 3.1]{Zha08a} that $C^l_\lambda(\dist^2(\cdot;\,K))\in C^{1,1}(\mathbb{R}^n)$, so
we only need to improve the estimate of the Lipschitz constant obtained in \cite[pag. 755]{Zha08a}, namely
 $8+10\lambda$. 
From the definition of the lower transform, we have 
\begin{equation}\label{Sec2.Th.RegLwTr.Prf.Eq1}
	C^l_\lambda(\dist^2(\cdot;\,K))(x+y)-C^l_\lambda(\dist^2(\cdot;\,K))(x)-
	(DC^l_\lambda(\dist^2(\cdot;\,K))(x))\cdot y\geq -\lambda|y|^2,
\end{equation}
for all $x,\, y\in \mathbb{R}^n$. Now $\dist^2(\cdot;\,K)$ is a $2$-semiconcave function \cite[Prop. 2.2.2]{CSi04}, that is, 
$x\mapsto |x|^2-\dist^2(x,\,K):=g(x)$ is a convex function. So if we let  $f_\lambda(x):=(\lambda+1)|x|^2-g(x)$, then
by Lemma \ref{Sec2.Lem.KZLem43}, 
we have
\begin{equation}\label{Sec2.Th.RegLwTr.Prf.Eq2}
	\co[f_\lambda](x+y)-\co[f_\lambda](x)-(D\co[f_\lambda](x))\cdot y\leq (\lambda+1)|y|^2,
\end{equation}
for all $x,\, y\in \mathbb{R}^n$. We also have, for $x\in \mathbb{R}^n$, that
\[
	C^l_\lambda(\dist^2(\cdot;\,K))(x)=\co[\dist^2(\cdot;\,K)+\lambda|\cdot|^2](x)-\lambda|x|^2=
	\co[f_\lambda(\cdot)](x)-\lambda|x|^2\,.
\]
By \eqref{Sec2.Th.RegLwTr.Prf.Eq2}, we obtain
 	\begin{eqnarray}
		\lefteqn{C^l_\lambda(\dist^2(\cdot;\,K))(x+y)-C^l_\lambda(\dist^2(\cdot;\,K))(x)-
		(DC^l_\lambda(\dist^2(\cdot;\,K))(x))\cdot y} \label{Sec2.Th.RegLwTr.Prf.Eq3}\\
		&&=\co[f_\lambda](x+y)-\co[f_\lambda](x)-(D\co[f_\lambda](x))\cdot y-\lambda|y|^2\leq 
		(\lambda+1)|y|^2-\lambda|y|^2 = |y|^2\,. \nonumber
\end{eqnarray}

Combining \eqref{Sec2.Th.RegLwTr.Prf.Eq1} and \eqref{Sec2.Th.RegLwTr.Prf.Eq3}, we have
\[
	-\lambda|y|^2\leq C^l_\lambda(\dist^2(\cdot;\,K))(x+y)-C^l_\lambda(\dist^2(\cdot;\,K))(x)-
	(DC^l_\lambda(\dist^2(\cdot;\,K))(x))\cdot y\leq |y|^2\,.
\]
Thus $x\mapsto C^l_\lambda(\dist^2(\cdot;\,K))(x)$ is both $2$-semiconcave and 
$2\lambda$-semiconvex. By \cite[Corollary 3.3.8]{CSi04}, we therefore conclude that 
$C^l_\lambda(\dist^2(\cdot;\,K))\in C^{1,1}(\R^n)$, and the Lipschitz constant of 
the gradient  $DC^l_\lambda(\dist^2(\cdot;\,K))$ is not greater than $2\max\{1,\,\lambda\}$.
\hfill\qed\\

%%%%%%%%%%%%%%%%%%%%%%%%%%%%%%%%%%%%%%%%%%%%%%%%%%%%%%%%%%%%%%%%%%%%%%%%%%%
%%%%%%%%%%%%%%%%%%%%%%%%%%%%%%%%%%%%%%%%%%%%%%%%%%%%%%%%%%%%%%%%%%%%%%%%%%%
%%%%%%%%%%%%%%%%%%%%%%%%%%%%%%%%%%%%%%%%%%%%%%%%%%%%%%%%%%%%%%%%%%%%%%%%%%%

\noindent {\bf Proof of Proposition \ref{Sec2.Th.LipGradLwTr}:}
Let $\lambda>0$. Without loss of generality, we may assume that $x_0=0$.
We consider two different cases, depending on the  values of $C^l_\lambda(\dist^2(\cdot;\,\,K))(0)$ and $\dist^2(0;\,K)$. 
\medskip

\textit{Case (i)}: $C^l_\lambda(\dist^2(\cdot;\,\,K))(0)<\dist^2(0;\,K)$.
\medskip

In this case, by definition of the lower transform and the convex envelope, we have
\[
	C^l_\lambda(\dist^2(\cdot;\,K))(0)=\co[\dist^2(\cdot;\,K)+\lambda|\cdot|^2](0).
\]
Since the function $\dist^2(\cdot;\,K)+\lambda|\cdot|^2$ is also continuous, 
upper semidifferentiable, coercive and
$\co[\dist^2(\cdot;\,K)+\lambda|\cdot|^2]\in C^{1,1}$, the assumptions of Lemma \ref{Sec2.Lem.UpDif} 
are satisfied. Let the affine function $\ell$  and $x_1,\dots,x_k\in\mathbb{R}^n$, 
$\lambda_1>0,\dots,\lambda_k>0$ with
$2\leq k\leq n+1$ be as given by Lemma \ref{Sec2.Lem.UpDif}, satisfying \eqref{Sec2.Lem.UpDif.Eq}$(i)$
to  \eqref{Sec2.Lem.UpDif.Eq}$(v)$. 
We have, by \eqref{Sec2.Lem.UpDif.Eq}$(ii)$ and \eqref{Sec2.Lem.UpDif.Eq}$(iii)$, that
\begin{equation}\label{Sec2.Th.LipGradLwTr.Prf.Eq1}
	2a\cdot x_i+b=\dist^2(x_i;\,K)+\lambda|x_i|^2,\quad 2a=D(\dist^2(\cdot;\, K)+\lambda|\cdot|^2)(x_i)\,.
	\end{equation}
Therefore $x\mapsto \dist^2(x;\, K)$ is differentiable at $x_i$ for $i=1,2,\dots,k$. 
By \cite[Lemma 8.5.12]{Hor83}, we see that $K(x_i):=\{y_i\}$ consists of a single element 
$y_i\in K$, so that $\dist^2(x_i;\,K)=|x_i-y_i|^2$ and
$D\dist^2(\cdot;\,K)(x_i)=2(x_i-y_i)$. 
Thus \eqref{Sec2.Th.LipGradLwTr.Prf.Eq1} reduces in this case to
\begin{equation}\label{Sec2.Th.LipGradLwTr.Prf.Eq2}
	2a\cdot x_i+b=|x_i-y_i|^2+\lambda|x_i|^2,\quad 2a=2(x_i-y_i)+2\lambda x_i\,.
\end{equation}
By \eqref{Sec2.Th.LipGradLwTr.Prf.Eq2}$_{2}$ and knowing that $\sum^k_{i=1}\lambda_i=1$ and $\sum^k_{i=1}\lambda_ix_i=0$, 
we find
\begin{equation}\label{Sec2.Th.LipGradLwTr.Prf.Eq3}
	a=\sum^k_{i=1}\lambda_i\Big((x_i-y_i)+\lambda x_i\Big)=-\sum^k_{i=1}\lambda_iy_i\,,
\end{equation}
whereas by \eqref{Sec2.Th.LipGradLwTr.Prf.Eq2}$_{1}$, \eqref{Sec2.Th.LipGradLwTr.Prf.Eq3} and 
the strict convexity of $|\cdot|^2$, we obtain
\begin{equation}\label{Sec2.Th.LipGradLwTr.Prf.Eq4}
\begin{split}
	b=\sum^k_{i=1}\lambda_i(2a\cdot x_i+b)&=\sum^k_{i=1}\lambda_i(|x_i-y_i|^2+\lambda|x_i|^2)\\
					      &>\left|\sum^k_{i=1}\lambda_i(x_i-y_i)\right|^2
					      =\left|-\sum^k_{i=1}\lambda_iy_i\right|^2=|a|^2\,.
\end{split}
\end{equation}
Hence $|a|^2<b$, that is,
\[
	\left|\frac{1}{2}DC^l_\lambda(\dist^2(\cdot;\,K))(0)\right|^2<C^l_\lambda(\dist^2(\cdot;\,K))(0)\,.
\]
\textit{Case (ii)}: $C^l_\lambda(\dist^2(\cdot;\,K))(0)=\dist^2(0;\,K)$.
\medskip

In this case,  we have 
\[ 
\co[\dist^2(\cdot;\,K)+\lambda|\cdot|^2](0) = C^l_\lambda(\dist^2(\cdot;\,K))(0)= \dist^2(0; K) + \lambda|0|^2,
\]
and
\[
	  \co[\dist^2(\cdot;\,K)+\lambda|\cdot|^2](x)
	\leq \dist^2(x;\,K)+\lambda|x|^2\,,
\]
for $x\in\R^n$. Since $\co[\dist^2(\cdot;\,K)+\lambda|\cdot|^2]\in C^{1,1}$ is convex and
$\dist^2(\cdot;\,K)+\lambda|\cdot|^2$ is upper-semidifferentiable, it follows from \cite[Corollary 2.5]{BKK00} that
$D(\dist^2(x;\,K)+\lambda|x|^2)$ exists at $0$, and 
\[
	D(\dist^2(\cdot;\,K)+\lambda|\cdot|^2)(0)=D\co[\dist^2(\cdot;\,K)+\lambda|\cdot|^2](0)\,.
\]
Again by  \cite[Lemma 8.5.12]{Hor83}, $K(0)=\{y_0\}$ with $y_0\in K$ the unique point
that realizes the distance of $0$ to $K$. So
\begin{equation}
\begin{split}
	|D\co[\dist^2(\cdot;\,K)+\lambda|\cdot|^2](0)|^2&=|D(\dist^2(\cdot;\,K)+\lambda|\cdot|^2)(0)|^2\\
							&=|-2y_0|^2=4\dist^2(0;\,K)=4\co[\dist^2(\cdot;\,K)+\lambda|\cdot|^2](0)\,.
\end{split}
\end{equation}
\hfill\qed\\

%%%%%%%%%%%%%%%%%%%%%%%%%%%%%%%%%%%%%%%%%%%%%%%%%%%%%%%%%%%%%%%%%%%%%%%%%%%
%%%%%%%%%%%%%%%%%%%%%%%%%%%%%%%%%%%%%%%%%%%%%%%%%%%%%%%%%%%%%%%%%%%%%%%%%%%
%%%%%%%%%%%%%%%%%%%%%%%%%%%%%%%%%%%%%%%%%%%%%%%%%%%%%%%%%%%%%%%%%%%%%%%%%%%

\noindent {\bf Proof of Theorem \ref{Sec3.Thm.TghtApprx}:}
Let $\lambda>0$. Without loss of generality, we may assume that $x_0=0$. We prove our result by establishing the contrapositive, and therefore suppose that
\begin{equation} \label{Sec3.Thm.TghtApprx.Prf.Eq1}
	C^l_\lambda(\dist^2(\cdot;\,K))(0)<\dist^2(0;\,K)
\end{equation}
and seek to prove that 
\[
	\lambda \dist(0;\, M_K)< \dist(0;\, K)\,.
\]
As in the proof of Proposition \ref{Sec2.Th.LipGradLwTr}, all the assumptions 
of Lemma \ref{Sec2.Lem.UpDif} are met for the function $f(\cdot) =\dist^2(\cdot; K) + \lambda|\cdot|^2$. Hence there exist an affine function
$\ell$,  points $x_1,\dots,x_k\in\mathbb{R}^n$, and $\lambda_1>0,\dots,\lambda_k>0$,
$2\leq k\leq n+1$,  that satisfy  \eqref{Sec2.Lem.UpDif.Eq}$(i)$ to 
\eqref{Sec2.Lem.UpDif.Eq}$(v)$, which ensures that \eqref{Sec2.Th.LipGradLwTr.Prf.Eq2}  holds
with $\{y_i\} = K(x_i)\}$, $i \in \{1, \ldots, k\}$. 
From \eqref{Sec2.Th.LipGradLwTr.Prf.Eq2}$_{2}$, we  have
\begin{equation} \label{Sec3.Thm.TghtApprx.Prf.Eq2}
	x_i=\frac{y_i+a}{1+\lambda}\,,
\end{equation}
which when substituted into \eqref{Sec2.Th.LipGradLwTr.Prf.Eq2}$_{1}$ yields
\[
	2a\cdot\left(\frac{y_i+a}{1+\lambda}\right)+b=\left|\frac{y_i+a}{1+\lambda}-y_i\right|^2+\lambda\left|
	\frac{y_i+a}{1+\lambda}\right|^2\,.
\]
A simple manipulation of this equation in $y_i$ then gives
\begin{equation} \label{Sec3.Thm.TghtApprx.Prf.Eq3}
	\left|y_i-\frac{a}{\lambda}\right|^2=
		\frac{(1+\lambda)|a|^2}{\lambda^2}+\frac{(1+\lambda)b}{\lambda}:=c^2,
		\quad i=1,2,\dots,k,
\end{equation}
and  \eqref{Sec3.Thm.TghtApprx.Prf.Eq2} and  \eqref{Sec3.Thm.TghtApprx.Prf.Eq3} together imply
\begin{equation} \label{Sec3.Thm.TghtApprx.Prf.Eq4}
	\left|x_i-\frac{a}{\lambda}\right|^2=\frac{c^2}{(1+\lambda)^2}.
	\quad i=1,2,\dots,k\,.
\end{equation}
Now from \eqref{Sec3.Thm.TghtApprx.Prf.Eq3}, it follows that the points $y_1,y_2,\dots, y_k\in K$ lie on the sphere 
$S(a/\lambda;\,c)=\{y\in \mathbb{R}^n,\, |y-a/\lambda|=c\}$, and since $x_i\not = x_j$ for $i \not = j$,
we also have that $y_i\neq y_j$ for $i\neq j$, by \eqref{Sec3.Thm.TghtApprx.Prf.Eq2}. 
We show next that the open ball $B(a/\lambda;\,c)$ does not intersect $K$, and hence 
$a/\lambda\in M_K$, the medial axis of $K$. 
We prove this claim by contradiction. Suppose $y^\ast\in K\cap B(a/\lambda;\, c)$, and define
\begin{equation}\label{Sec3.Thm.TghtApprx.Prf.Eq5}
	x^\ast=\frac{y^\ast+a}{1+\lambda}\,.
\end{equation}
Then we have, from \eqref{Sec2.Lem.UpDif.Eq}$(i)$, that
\begin{equation}\label{Sec3.Thm.TghtApprx.Prf.Eq6}
	2a\cdot x^\ast+b=\ell(x^\ast)\leq \dist^2(x^\ast,\, K)+
	\lambda|x^\ast|^2\leq |x^\ast-y^\ast|^2+\lambda|x^\ast|^2\,
\end{equation}
and by replacing \eqref{Sec3.Thm.TghtApprx.Prf.Eq5} into \eqref{Sec3.Thm.TghtApprx.Prf.Eq6}, it follows that
\[
	\left|y^\ast-\frac{a}{\lambda}\right|^2\geq c^2\,,
\]
which contradicts the assumption that $y^\ast\in B(a/\lambda;\,c)$. Hence 
\begin{equation}
\label{Sec3.Thm.TghtApprx.Prf.Eq6a}
B\left(\frac{a}{\lambda}; c\right) \cap K = \emptyset,
\end{equation}
and thus $a/\lambda\in M_K$.
By the strict convexity of $|\cdot|^2$ and the fact that $\sum^k_{i=1}\lambda_ix_i=0$, 
we then have, from \eqref{Sec3.Thm.TghtApprx.Prf.Eq4}, \eqref{Sec2.Th.LipGradLwTr.Prf.Eq4} and Proposition \ref{Sec2.Th.LipGradLwTr}, that
\begin{equation}
	\begin{split}
		\dist^2(0;\, M_K)\leq \left|\frac{a}{\lambda}\right|^2
					<\sum^k_{i=1}\lambda_i\left|x_i-\frac{a}{\lambda}\right|^2
		&=\frac{c^2}{(\lambda+1)^2}
					=\frac{1}{(\lambda+1)^2}\left(\frac{(1+\lambda)|a|^2}{\lambda^2}
					+\frac{(1+\lambda)b}{\lambda}\right)\\
		&=\frac{|a|^2}{\lambda^2(1+\lambda)}+\frac{b}{\lambda(1+\lambda)}
		<\frac{b}{\lambda^2}<
		\frac{\dist^2(0;\, K)}{\lambda^2}\,,
	\end{split}
\end{equation}
and hence
\[
	\lambda\dist(0;\, M_K)<\dist(0;\, K)\,.
\]	
This proves that if $0\notin \overline{M_K}$ and 
\[
	\lambda\geq \frac{\dist(0;\, K)}{\dist(0;\, M_K)}\,,
\]
then
$C^l_\lambda(\dist^2(\cdot;\,K))(0)\geq\dist^2(0;\,K)$.
Since we always have that $C^l_\lambda(\dist^2(\cdot;\,K))(0)\leq\dist^2(0;\,K)$, it can be concluded that
\[
	C^l_\lambda(\dist^2(\cdot;\,K))(0)=\dist^2(0;\,K)\,,
\]
which completes the proof. \hfill\qed\\

\begin{nota}
\label{Sec3.Thm.TghtApprx.Rk}
Recall that in Remark \ref{Sec3.Rmk.locality}(b), we noted that if $x_0$ is a critical point of 
	 $C^l_\lambda(\dist^2(\cdot;\,K))$, 
	then $x_0\in \co[K(x_0)]$. Translating $x_0$ to $0$, we can now see  that if $C^l_{\lambda}(\dist^2(\cdot;\,K))(0) < \dist^2(0; K)$, this follows from \eqref{Sec2.Th.LipGradLwTr.Prf.Eq3}, \eqref{Sec3.Thm.TghtApprx.Prf.Eq3} and \eqref{Sec3.Thm.TghtApprx.Prf.Eq6a}, since Lemma \ref{Sec2.Lem.UpDif}(v) implies that $a=0$ if $DC^l_{\lambda}(\dist^2(\cdot;\,K))(0)=0$, whereas if $C^l_{\lambda}(\dist^2(\cdot;\,K))(0) = \dist^2(0; K)$,  the arguments in the proof of Proposition \ref{Sec2.Th.LipGradLwTr}(ii) yield
	that $0 \in K$ if $DC^l_{\lambda}(\dist^2(\cdot;\,K))(0)=0$, thus clearly $0 \in  \co [K(0)]$ in this case also. 
\end{nota}

\medskip
%%%%%%%%%%%%%%%%%%%%%%%%%%%%%%%%%%%%%%%%%%%%%%%%%%%%%%%%%%%%%%%%%%%%%%%%%%%
%%%%%%%%%%%%%%%%%%%%%%%%%%%%%%%%%%%%%%%%%%%%%%%%%%%%%%%%%%%%%%%%%%%%%%%%%%%
%%%%%%%%%%%%%%%%%%%%%%%%%%%%%%%%%%%%%%%%%%%%%%%%%%%%%%%%%%%%%%%%%%%%%%%%%%%

\noindent {\bf Proof of Corollary \ref{Sec3.Cor.Sprt}:}
Note first that (\ref{Sec3.prop.vlk}) and (\ref{Sec3.subset}) together yield that
$$\bigcap_{0<\lambda<+ \infty} \sprt(M_{\lambda}(\cdot; K)) \subset \bigcap_{0<\lambda< + \infty} V_{\lambda, K} = \overline{M_K}.$$
On the other hand, suppose $x_0 \in \R^N$ is such that $M_{\lambda}(x_0; K)=0$. Then we have that $\dist^2(x_0; K) =
C^{l}_{\lambda} (\dist^2(\cdot; K))(x_0)$, so as argued in the proof of Proposition \ref{Sec2.Th.LipGradLwTr} (ii), it follows that $x_0 \not\in M_K$. Thus
or all $\lambda >0$, $M_K \subset \{x \in \R^N: M_{\lambda}(x; K)>0\}$,
which implies that $$\overline{M_K} \subset \bigcap_{0<\lambda <+\infty} \sprt M_{\lambda}( \cdot; K).$$
\hfill\qed\\

%%%%%%%%%%%%%%%%%%%%%%%%%%%%%%%%%%%%%%%%%%%%%%%%%%%%%%%%%%%%%%%%%%%%%%%%%%%
%%%%%%%%%%%%%%%%%%%%%%%%%%%%%%%%%%%%%%%%%%%%%%%%%%%%%%%%%%%%%%%%%%%%%%%%%%%
%%%%%%%%%%%%%%%%%%%%%%%%%%%%%%%%%%%%%%%%%%%%%%%%%%%%%%%%%%%%%%%%%%%%%%%%%%%

\noindent {\bf Proof of Corollary \ref{Sec3.Cor.LclProp}:}
We only need to consider the case where $C^l_\lambda(\dist^2(\cdot;\,K))(x_0)<\dist^2(x_0;\,K)$, since
otherwise the claim is clearly true. Without loss of generality, assume that $x_0=0$. 
As in the proof of Theorem \ref{Sec3.Thm.TghtApprx},  the assumptions of Lemma \ref{Sec2.Lem.UpDif} are satisfied for the function $f(\cdot) =\dist^2(\cdot; K) + \lambda|\cdot|^2$, so
there exist an affine function  $\ell(x) = 2a\cdot x+b$,  points $x_1,\dots,x_k\in\mathbb{R}^n$, and $\lambda_1>0,\dots,\lambda_k>0$, 
$2\leq k\leq n+1$, that satisfy  \eqref{Sec2.Lem.UpDif.Eq}$(i)$ to \eqref{Sec2.Lem.UpDif.Eq}$(v)$,
so that 
\eqref{Sec3.Thm.TghtApprx.Prf.Eq4} holds. Moreover, from Proposition \ref{Sec2.Th.LipGradLwTr} and the assumption that
$C^l_\lambda(\dist^2(\cdot;\,K))(0)<\dist^2(0;\,K),$  it follows that
$|a|^2<b< \dist^2(0; K).$ Hence, for each for $i=1,2,\dots,k$, we have
\[
	|x_i|-\frac{|a|}{\lambda}\leq\left |x_i-\frac{a}{\lambda}\right|=
	\frac{c}{1+\lambda}<\frac{\sqrt{b}}{\lambda}\,,
\]
and so
\[
	|x_i|<\frac{|a|}{\lambda}+\frac{\sqrt{b}}{\lambda}<\frac{2\sqrt{b}}{\lambda}<\frac{2\dist(0;\,K)}{\lambda}\,.
\]
Thus $x_i\in B(0;\,r(0))$ for $i=1,2,\dots,k$, where $r(0)=2\dist(0;\,K)/\lambda$.
The conclusion then follows from the definitions of the convex envelope and of $C^l_\lambda(\dist^2(\cdot;\,K))(0)$.
\hfill\qed\\

%%%%%%%%%%%%%%%%%%%%%%%%%%%%%%%%%%%%%%%%%%%%%%%%%%%%%%%%%%%%%%%%%%%%%%%%%%%
%%%%%%%%%%%%%%%%%%%%%%%%%%%%%%%%%%%%%%%%%%%%%%%%%%%%%%%%%%%%%%%%%%%%%%%%%%%
%%%%%%%%%%%%%%%%%%%%%%%%%%%%%%%%%%%%%%%%%%%%%%%%%%%%%%%%%%%%%%%%%%%%%%%%%%%

\noindent {\bf Proof of Corollary \ref{Sec3.Cor.Reg}:}
Since $\overline{M_K} \subset V_{\lambda, K}$, it follows from Theorem \ref{Sec3.Thm.TghtApprx} that
$\dist^2(x; K) = C^l_\lambda(\dist^2(\cdot;\,K))(x)$ for all $x$ in the set $\R^n \setminus V_{\lambda, K}$, which is
an open set because $V_{\lambda, K}$ is closed. Thus  the result is immediate from Proposition \ref{Sec2.Th.RegLwTr}. \hfill\qed\\

\noindent {\bf Proof of Proposition \ref{Sec3.Pro.MACmp}:}
Since $\partial \Omega\subset \Omega^c$, then 
\begin{equation}
	\dist^2(y,\Omega^c)\leq \dist^2(y,\partial \Omega)\quad\text{for }y\in \R^n\,,
\end{equation}
and by the ordering property of the lower transform,  
\begin{equation}
 	C^l_\lambda(\dist^2(\cdot;\,\Omega^c))(x)\leq C^l_\lambda(\dist^2(\cdot;\,\partial\Omega))(x) \quad 
 	\text{for } x\in \Omega.
\end{equation}
We want now to prove that the equality  actually holds for $x \in \Omega$. We will show this 
by a contradiction argument. Though the equality $\dist^2(y;\,\Omega^c)= \dist^2(y;\,\partial \Omega)$ 
holds for all $y\in \Omega$, we cannot straightforwardly deduce the equality of the lower transforms in $\Omega$.
Assume therefore, that at some point $x\in \Omega$ we have that 
\begin{equation}\label{Eq.MACmp.01}
 	C^l_\lambda(\dist^2(\cdot;\,\Omega^c))(x)< C^l_\lambda(\dist^2(\cdot;\,\partial\Omega))(x)\,.
\end{equation}
By the translation invariance of the distance and of the lower transform \cite[Proposition 2.10]{ZOC14a}, we can assume, 
without loss of generality, that $x=0\in\Omega$, so that \eqref{Eq.MACmp.01} becomes
\begin{equation}
 	C^l_\lambda(\dist^2(\cdot;\,\Omega^c))(0)< C^l_\lambda(\dist^2(\cdot;\,\partial\Omega))(0)\,,
\end{equation}
which is then equivalent to state that 
\begin{equation}\label{Eq.MACmp.02}
 	\co[\dist^2(\cdot;\,\Omega^c)+\lambda|\cdot|^2](0)<\co[\dist^2(\cdot;\,\partial\Omega)+\lambda|\cdot|^2](0)\,.
\end{equation}
Since the function $\dist^2(x,\partial\Omega)+\lambda|x|^2$ is coercive 
and is continuous,  by Proposition \ref{Prp.GlCnvx}$(ii)$ there exists an affine function $\ell(x)$ such that
\begin{equation}\label{Eq.MACmp.03}
 	\ell(x)\leq \dist^2(x;\,\partial\Omega)+\lambda|x|^2\quad\text{for all }x\in \mathbb{R}^n\,,
\end{equation}
and 
\begin{equation}\label{Eq.MACmp.04}
 	\ell(0)=\co[\dist^2(\cdot;\,\partial\Omega)+\lambda|\cdot|^2](0)\,.
\end{equation}
Note that by Proposition \ref{Prp.GlCnvx}$(i)$, $\ell(0)\geq 0$.
There must be a point $y_0\in (\overline\Omega)^c$ such that
\begin{equation}\label{Eq.MACmp.05}
 	\ell(y_0)> \dist^2(y_0;\,\Omega^c)+\lambda|y_0|^2=\lambda|y_0|^2.
\end{equation}
If we write our affine function as $\ell(x)=a\cdot x+b$ with $b=\ell(0)$ given by \eqref{Eq.MACmp.04} and $a\in\mathbb{R}^n$,
then \eqref{Eq.MACmp.05} reads as
\begin{equation}\label{Eq.MACmp.06}
	a\cdot y_0+b>\lambda|y_0|^2	
\end{equation}
Since for $x\in\partial\Omega$,  
\begin{equation}\label{Eq.MACmp.07}
 	\ell(x)\leq \dist^2(x;\,\partial\Omega)+\lambda|x|^2=
 	\dist^2(x;\,\Omega^c)+\lambda|x|^2\,,
\end{equation}
and for $t>0$ small enough, $ty_0\in \Omega$, it follows that there exists $0<t_0<1$ for which $t_0y_0\in \Omega^c$ and 
\begin{equation}\label{Eq.MACmp.08}
 	\ell(t_0y_0)=\dist^2(t_0y_0;\,\Omega^c)+\lambda|t_0y_0|^2=\lambda|t_0y_0|^2\,.
\end{equation}
This implies that 
\begin{equation}\label{Eq.MACmp.09}
	t_0a\cdot y_0+b=\lambda|t_0y_0|^2\,,
\end{equation}
that is, 
\begin{equation}\label{Eq.MACmp.10}
	a\cdot y_0= t_0\lambda|y_0|^2-b/t_0\,.
\end{equation}
If we substitute \eqref{Eq.MACmp.10} into \eqref{Eq.MACmp.06}, we obtain
\begin{equation}
 	t_0\lambda|y_0|^2-b/t_0+b>\lambda|y_0|^2\,,
\end{equation}
that is,
\begin{equation}
	b<-t_0\lambda |y_0|^2<0\,,
\end{equation}
which contradicts the fact that $b\geq 0$. Thus 
\begin{equation}
 	\ell(x)\leq \dist^2(x;\,\Omega^c)+\lambda|x|^2\quad\text{for all }x\in \R^n\,,
\end{equation}
hence,
\begin{equation}
	\co[\dist^2(\cdot;\,\Omega^c)+\lambda|\cdot|^2](0)\geq 
	\co[\dist^2(\cdot;\,\partial\Omega)+\lambda|\cdot|^2](0)\,,
\end{equation}
thus
\begin{equation}
	C^l_\lambda(\dist^2(\cdot;\,\Omega^c))(0)   \geq 
	C^l_\lambda(\dist^2(\cdot;\,\partial\Omega))(0)\,,
\end{equation}
which contradicts the initial assumption \eqref{Eq.MACmp.01}. \hfill \qed\\

%%%%%%%%%%%%%%%%%%%%%%%%%%%%%%%%%%%%%%%%%%%%%%%%%%%%%%%%%%%%%%%%%%%%%%%%%%%
%%%%%%%%%%%%%%%%%%%%%%%%%%%%%%%%%%%%%%%%%%%%%%%%%%%%%%%%%%%%%%%%%%%%%%%%%%%
%%%%%%%%%%%%%%%%%%%%%%%%%%%%%%%%%%%%%%%%%%%%%%%%%%%%%%%%%%%%%%%%%%%%%%%%%%%

\noindent {\bf Proof of Corollary \ref{Sec3.Cor.Reg.OpenDmn}:}
We only need to verify that $V_{\lambda,\Omega^c}\cap \Omega\subset W_{\lambda,\Omega}$. 
In fact,  if $x\in V_{\lambda,\Omega^c}\cap\Omega$, then 
$\lambda\dist(x;\, M_{\Omega^c}) \leq \dist(x;\Omega^c)\leq \diam(\Omega).$ Thus
$x\in W_{\lambda,\Omega}$ so that $\Omega\setminus W_{\lambda,\Omega}\subset \Omega\setminus
V_{\lambda,\Omega^c}$. The conclusion then follows from Corollary \ref{Sec3.Cor.Reg}.
\hfill\qed\\

%%%%%%%%%%%%%%%%%%%%%%%%%%%%%%%%%%%%%%%%%%%%%%%%%%%%%%%%%%%%%%%%%%%%%%%%%%%
%%%%%%%%%%%%%%%%%%%%%%%%%%%%%%%%%%%%%%%%%%%%%%%%%%%%%%%%%%%%%%%%%%%%%%%%%%%
%%%%%%%%%%%%%%%%%%%%%%%%%%%%%%%%%%%%%%%%%%%%%%%%%%%%%%%%%%%%%%%%%%%%%%%%%%%

\noindent {\bf Proof of Theorem \ref{Sec3.Teo.BndMMA}:}
Note first that clearly $M_{\lambda}(x; \,K) \geq 0$ for all $x \in \R^n$. 
We now prove the positive lower bound for $M_{\lambda}(x;\,K)$ when $x \in M_K$.
Let $x\in M_K$ and $r(x)=\dist(x;\,K)$. Since $K(x)\subset K$ then 
\begin{equation}\label{Eq.Prf.BndMMA.01}
	\dist^2(y;\,K)\leq \dist^2(y;\,K(x))\quad\text{for all }y\in \mathbb{R}^n\,,
\end{equation}
hence, by the ordering property of the lower transform, Proposition \ref{Sec2.Pro.OrdLwTr}, 
\begin{equation} \label{Eq.Prf.BndMMA.02}
	C^l_\lambda(\dist^2(\cdot;\,K))(y)\leq C^l_\lambda(\dist^2(\cdot;\,K(x)))(y)
	\quad\text{ for all }y\in \mathbb{R}^n\,.
\end{equation}
Now, by Proposition \ref{Sec2.Pro.LwTrSetK}, we have for $y=x$ that
\begin{equation}\label{Eq.Prf.BndMMA.03}
	C^l_\lambda(\dist^2(\cdot;\,K(x)))(x)=(1+\lambda)\dist^2(x;\,\co[K_{2,\lambda}(x)])+\frac{\lambda}{1+\lambda}r^2(x)\,.
\end{equation}
By the Carath\'eodory's theorem \cite{Roc66}, for every
$w=x+y\in \co[K_{2,\lambda}(x)]$, there are at most $n+1$ points $x+y_i/(1+\lambda)\in 
K_{2,\lambda}(x)$ with $x+y_i\in K(x)$, i.e. $x+y_i\in \partial K$  and $|y_i|=r(x)$,
and $\lambda_1\geq 0,\dots, \lambda_{n+1}\geq 0$ such that $\sum^{n+1}_{i=1}\lambda_i=1$ and 
$w=x+\sum^{n+1}_{i=1}\lambda_i y_i/(1+\lambda)$.
Thus we have
\begin{equation}\label{Eq.Prf.BndMMA.04}
	|x-w|^2=\left|\sum^{n+1}_{i=1}\lambda_i\frac{y_i}{1+\lambda}\right|^2=\frac{1}{(1+\lambda)^2}
	\left|\sum^{n+1}_{i=1}\lambda_i y_i\right|^2= 
	\frac{1}{(1+\lambda)^2} \left|x-\sum^{n+1}_{i=1}\lambda_i(x+y_i)\right|^2.
\end{equation}
which yields
\begin{equation}\label{Eq.Prf.BndMMA.05}
	\dist^2(x;\,\co[K_{2,\lambda}(x)])=\frac{1}{(1+\lambda)^2}\dist^2(x;\,\co[K(x)])\,.
\end{equation}
By substituting \eqref{Eq.Prf.BndMMA.05} into \eqref{Eq.Prf.BndMMA.03} we have
\begin{equation}\label{Eq.Prf.BndMMA.06}
	C^l_\lambda(\dist^2(\cdot;\,K(x)))(x)=\frac{1}{1+\lambda}\dist^2(x;\,\co[K(x)])+\frac{\lambda}{1+\lambda}r^2(x)\,.
\end{equation}
By comparing \eqref{Eq.Prf.BndMMA.02} and \eqref{Eq.Prf.BndMMA.06}, we finally obtain
\begin{equation*}\label{Eq.Prf.BndMMA.07}
	M_{\lambda}(x;\,K)\geq (1+\lambda)(\dist^2(x;\,K)-C^l_\lambda(\dist^2(x;\,K(x))))
			= \dist^2(x;\,K)-\dist^2(x;\,\co[K(x)])\,.
\end{equation*}
To find now an upper bound to $M_{\lambda}(x;\,K)$ that holds for all $x \in \R^n$, 
note first that if $x \in K$, then $M_{\lambda}(x;\,K)=\dist^2(x; \,K)=0$. Suppose now that
$x \not\in K$. Then $r(x)>0$ and $B(x;\, r(x)) \subset K^c$, so $K \subset B^c(x; r(x))$ and hence
 \begin{equation}\label{Eq.Prf.BndMMA.08}
	\dist^2(y;\,B^c(x,r(x)))\leq  \dist^2(y;\,K)\quad\text{for }y\in \R^n\,,
\end{equation}
and thus,  by the ordering property of the lower transform (Proposition \ref{Sec2.Pro.OrdLwTr}), 
\begin{equation}\label{Eq.Prf.BndMMA.09}
	C_{\lambda}^l(\dist^2(\cdot;\,B^c(x,r(x))))(y)\leq  C_{\lambda}^l(\dist^2(\cdot;\,K))(y)\quad\text{for }y\in \R^n\,.
\end{equation}
By Lemma \ref{Sec2.Lem.LwDst}, for $|y|\leq r(x)/(1+\lambda)$, after a simple translation of points and due to 
the invariance of the distance transform, we have
\begin{equation}\label{Eq.Prf.BndMMA.10}
	C_{\lambda}^l(\dist^2(\cdot;\,B^c(x,r(x))))(y)=\frac{\lambda}{1+\lambda}r^2(x)-\lambda|y-x|^2
\end{equation}
which for $y=x$ gives
\begin{equation}\label{Eq.Prf.BndMMA.11}
	C_{\lambda}^l(\dist^2(\cdot;\,B^c(x,r(x))))(x)=\frac{\lambda}{1+\lambda}r^2(x)\,,
\end{equation}
hence 
\begin{equation}\label{Eq.Prf.BndMMA.12}
	\dist^2(x;\,K)-C_{\lambda}^l(\dist^2(\cdot;\,B^c(x,r(x))))(x)=\frac{1}{1+\lambda}\dist^2(x;\,K)\,.
\end{equation}
By comparing \eqref{Eq.Prf.BndMMA.09} and \eqref{Eq.Prf.BndMMA.12}, we then conclude that
\begin{equation*}\label{Eq.Prf.BndMMA.13}
	\dist^2(x;\,K)-C_{\lambda}^l(\dist^2(\cdot;\,K))(x)\leq\dist^2(x;\,K)-C_{\lambda}^l(\dist^2(\cdot;\,B^c(x;\,r(x))))(x)
						=\frac{1}{1+\lambda}\dist^2(x;\,K)\,,
\end{equation*}
which completes the proof.\hfill\qed\\

%%%%%%%%%%%%%%%%%%%%%%%%%%%%%%%%%%%%%%%%%%%%%%%%%%%%%%%%%%%%%%%%%%%%%%%%%%%
%%%%%%%%%%%%%%%%%%%%%%%%%%%%%%%%%%%%%%%%%%%%%%%%%%%%%%%%%%%%%%%%%%%%%%%%%%%
%%%%%%%%%%%%%%%%%%%%%%%%%%%%%%%%%%%%%%%%%%%%%%%%%%%%%%%%%%%%%%%%%%%%%%%%%%%

\noindent {\bf Proof of Proposition \ref{Sec3.Pro.OpAng}:}
Let $x\in M_K$, $r(x)=\dist(x;\,K)>0$, and denote by 
$x_1,\, x_2\in\overline{B}(x;\,r(x))\cap K$ the points of
$K(x)$ that realize the separation angle $\theta_x$ at the point $x$. Thus
$|x-x_1|=|x-x_2|=r(x)$.  Since $\{x_1,x_2\}\subset K(x)$, then $x_1+t(x_2-x_1)\in \co[K(x)]$
for $t\in [0,1]$, which for $t=1/2$ gives
\begin{equation}\label{Eq.Prf.OpAng.01}
	\dist (x;\,\co[K(x)])\leq \dist (x;\,\frac{x_1+x_2}{2}) = \cos \left( \frac{\theta_x}{2} \right) \, r(x)\,.
\end{equation}
Thus
\begin{equation}
	\dist^2(x;\,\co[K(x)])\leq \dist^2(x;\,K)\,\cos^2\left(\frac{\theta_x}{2}\right)
\end{equation}
and hence
\begin{equation}
	\dist^2(x;\,K)-\dist^2(x;\,\co[K(x)])\geq 
	\dist^2(x;\,K)\,( 1- \cos^2 \left( \frac{\theta_x}{2} \right)) =  
	\dist^2(x;\,K)\, \sin^2\left(\frac{\theta_x}{2}\right)\,,
\end{equation}
as required. 
\hfill\qed\\

%%%%%%%%%%%%%%%%%%%%%%%%%%%%%%%%%%%%%%%%%%%%%%%%%%%%%%%%%%%%%%%%%%%%%%%%%%%%%%%%%%%%%
%%%%%%%%%%%%%%%%%%%%%%%%%%%%%%%%%%%%%%%%%%%%%%%%%%%%%%%%%%%%%%%%%%%%%%%%%%%
%%%%%%%%%%%%%%%%%%%%%%%%%%%%%%%%%%%%%%%%%%%%%%%%%%%%%%%%%%%%%%%%%%%%%%%%%%%

\noindent {\bf Proof of Proposition \ref{Sec3.LimInfSup}:}
If $x_0\in K$, clearly $M_\lambda(x_0;\,K)=0$ for all $\lambda>0$. So we may assume that $x_0\notin K$. 
Since $x_0\notin M_K$, $\dist^2(\cdot;\,K)$ is differentiable at $x_0$ \cite[Lemma 8.5.12]{Hor83}. 
Therefore for every $\epsilon>0$, there exists $\delta>0$ such that
\[
	|\dist^2(x_0+y;\,K)-\dist^2(x_0;\,K)-D\dist^2(x_0;\,K)\cdot y|\leq \epsilon|y|
\]
for $y\in \overline B(0;\,\delta)$. Now by the locality property  Corollary \ref{Sec3.Cor.LclProp}, we have
\[
	C^l_{\lambda}(\dist^2(\cdot;\, K))(x_0)=
		\co_{\overline B(x_0;\,r(x_0))}[\dist^2(\cdot;\, K)+\lambda|(\cdot)\, -\,x_0|^2](x_0)
\]
where $r(x_0)=2\dist(x_0;\,K)/\lambda>0$. Thus for $\lambda>0$ sufficiently large, $r(x_0)<\delta$. 
Since $x\mapsto\dist^2(x;\,K)+\lambda|x-x_0|^2$ is continuous and coercive, by Lemma \ref{Sec2.Lem.UpDif} 
and Corollary \ref{Sec3.Cor.LclProp}, there exist $x_1,\dots,x_k\in \overline B(x_0;\,r(x_0))$ and 
$\lambda_1>0,\dots,\lambda_k>0$ such that $\sum^k_{i=1}\lambda_i=1$,  $\sum^k_{i=1}\lambda_ix_i=x_0$
and
\[
	\begin{array}{l}
		\displaystyle	\co_{\bar B(x_0;\,r(x_0))}[\dist^2(\cdot;\, K)+\lambda|(\cdot)-x_0|^2](x_0)=
		\displaystyle	\sum^k_{i=1}\lambda_i\left(\dist^2(x_i;\,K)+\lambda|x_i-x_0|^2\right)\\
		\displaystyle	\geq \sum^k_{i=1}\lambda_i\left(\dist^2(x_0;\, K)+D\dist^2(x_0;\, K)\cdot (x_i-x_0)-
		\displaystyle	\epsilon|x_i-x_0|+\lambda|x_i-x_0|^2\right)\\
		\displaystyle	\geq \dist^2(x_0;\, K)-\frac{\epsilon^2}{4\lambda}\,,
	\end{array}
\]
as $|x_i-x_0|\leq r(x_0)<\delta$. Here we have also used the facts that 
$\sum^k_{i=1}\lambda_iD\dist^2(x_0;\, K)\cdot (x_i-x_0)=0$ and that  
$\lambda t^2-\epsilon t\geq -\epsilon^2/(4\lambda)$ for $t\in\R$. 
Since we also have $C^l_\lambda(\dist^2(\cdot;\,K))(x_0)\leq \dist^2(x_0;\,K)$, we have
\[
	0\leq M_\lambda(x_0;\,K)\leq \frac{(1+\lambda)\epsilon^2}{4\lambda}\,.
\]
Thus
\[
	0\leq \liminf_{\lambda\to+\infty}M_\lambda(x_0;\,K)\leq
	\limsup_{\lambda\to+\infty}M_\lambda(x_0;\,K)\leq \frac{\epsilon^2}{4}\,.
\]
Since $\epsilon>0$ is arbitrary, the conclusion follows.
\hfill\qed\\

%%%%%%%%%%%%%%%%%%%%%%%%%%%%%%%%%%%%%%%%%%%%%%%%%%%%%%%%%%%%%%%%%%%%%%%%%%%
%%%%%%%%%%%%%%%%%%%%%%%%%%%%%%%%%%%%%%%%%%%%%%%%%%%%%%%%%%%%%%%%%%%%%%%%%%%
%%%%%%%%%%%%%%%%%%%%%%%%%%%%%%%%%%%%%%%%%%%%%%%%%%%%%%%%%%%%%%%%%%%%%%%%%%%

\noindent {\bf Proof of Theorem \ref{Sec3.Teo.LimMMA}:}
We only consider the case $x_0\in M_K$. Again without loss of generality, we may assume 
that $x_0=0\in M_K$. Let $K_0=K(0)$ and $\dist(0;\, K)=\dist(0;\, K_0)=r_0>0$. 
Since $K_0\subset K$, we have $\dist^2(x;\,K)\leq \dist^2(x;\,K_0)$ for $x\in\R^n$, so that
\begin{equation}\label{Sec3.Teo.LimMMA.Prf.Eq1}
	C^l_\lambda(\dist^2(\cdot;\,K))(x)\leq
	C^l_\lambda(\dist^2(\cdot;\,K_0))(x)=\frac{\lambda r^2_0}{1+\lambda}+
	(1+\lambda)\dist^2(x;\, \co[K_0/(1+\lambda)])-
	\lambda|x|^2,
\end{equation}
for $x\in\mathbb{R}^n$. Therefore
\begin{equation}\label{Sec3.Teo.LimMMA.Prf.Eq2}
	\begin{split}
		M_\lambda(0;\,K)&\geq M_\lambda(0;\,K_0)=(1+\lambda)(\dist^2(0;\, K_0)-
					C^l_\lambda(\dist^2(\cdot;\,K_0))(0))\\
				&=r_0^2-\dist^2(0;\, \co[K_0])=M(0;\, K_0).
	\end{split}
\end{equation}

Next we establish lower bounds for $C^l_\lambda(\dist^2(\cdot;\,K))(x)$  using 
the locality property from Corollary \ref{Sec3.Cor.LclProp}. 
For  $0<\epsilon<r_0$ sufficiently small,  let $K_{0;\,\epsilon}\subset S(0;\,r_0)$ 
be the closed $\epsilon$-neighbourhood of $K_0$ on the sphere,  defined using the geodesic 
distance $\rho(x,\,y)$ on $S(0;\,r_0)$, that is
$K_{0;\,\epsilon}=\{y\in S(0;\,r_0),\; \rho(y;\,K_0)\leq \epsilon\},$
where $\rho(y;\,K_0)=\inf\{\rho(y,x),\, x\in K_0\}$. 

The aim of the following technical construction is to show that for $x$  in a small neighbourhood of $0$, 
$\dist^2(x;\, K_{0,2\epsilon})$ is a lower bound for $\dist^2(x;\, K)$.
For $\delta>0$,  define the closed neighbourhood
$K^\delta_{0,\epsilon}=\{ (r_0+t)y/|y|,\; y\in K_{0,\epsilon},\, 0\leq t\leq \delta\}$ and note that
 $K^\delta_{0,\epsilon}$ is clearly a compact set. 
Then it can easily be proved, using a contradiction argument, that for every $0<\epsilon<1$, there exists $0< \delta \leq  \epsilon^2$ such that 
$K\cap \overline B(0;\,r_0+\delta)\subset K^\delta_{0,\epsilon}$. 
Define also another compact set by $$V_{\epsilon,\delta}=K^\delta_{0,\epsilon}\cup K_{0,2\epsilon},$$
where $K_{0,2\epsilon}$ will be used to `shadow' $K^\delta_{0,\epsilon}$,
and the unbounded closed set $$W_{\epsilon,\delta}=V_{\epsilon,\delta}\cup B^c(0;\,r_0+\delta).$$
Clearly, $K\subset W_{\epsilon,\delta}$, so that
$\dist^2(x;\, K)\geq \dist^2(x;\, W_{\epsilon;\,\delta})$ for all $x\in \R^n$.\\
We claim that there exists $\eta>0$ sufficiently small such that 
\begin{equation}\label{Sec3.Teo.LimMMA.Prf.Eq3}
	\dist^2(x;\, W_{\epsilon;\,\delta})=\dist^2(x;\, K_{0,2\epsilon})
\end{equation}
for $x\in \bar B(0;\,\eta)$. We postpone the proof of \eqref{Sec3.Teo.LimMMA.Prf.Eq3} to the end 
and proceed first to assume that \eqref{Sec3.Teo.LimMMA.Prf.Eq3} holds.
Then for $\lambda>0$ sufficiently large, we have
\[
	\frac{2\dist(0; \, W_{\epsilon,\delta})}{\lambda}=
	\frac{2\dist(0; \, K_{0,2\epsilon})}{\lambda}=\frac{2r_0}{\lambda}<\eta\,.
\]
By the locality property (Corollary \ref{Sec3.Cor.LclProp}), we have
\begin{equation}
	\begin{array}{l}
	   \displaystyle C^l_\lambda(\dist^2(\cdot;\, W_{\epsilon,\delta}))(0)-
	   (\dist^2(\cdot;\, W_{\epsilon,\delta})+\lambda|\cdot|^2)(0)\\[1.5ex]
	\phantom{xxxxx}	 \displaystyle =   \co_{\overline B(0;\,\eta)}[\dist^2(\cdot;\, W_{\epsilon,\delta})+\lambda|\cdot|^2)](0)\\[1.5ex]
	\phantom{xxxxx}	 \displaystyle =  \co_{\overline B(0;\,\eta)}[\dist^2(\cdot;\, K_{0;2\epsilon})+\lambda|\cdot|^2](0)\\[1.5ex]
	\phantom{xxxxx}	 \displaystyle =  \co[ \dist^2(\cdot;\, K_{0,2\epsilon})+\lambda|\cdot|^2](0)\\[1.5ex]
	\phantom{xxxxx}		 \displaystyle  = C^l_\lambda(\dist^2(\cdot;\, K_{0,2\epsilon}))(0)\\[1.5ex]
	\phantom{xxxxx}	 \displaystyle =  \frac{\lambda r_0^2}{1+\lambda}+(1+\lambda)\dist^2\left(0;\,\co\left[\frac{K_{0,2\epsilon}}{1+\lambda}\right]\right)\,,
	\end{array}
\end{equation}
where we have used \eqref{Sec3.Teo.LimMMA.Prf.Eq3}  and Proposition \ref{Sec2.formula-set-on-sphere}. Thus we obtain
\[
	C^l_\lambda(\dist^2(\cdot;\, K))(0)\geq C^l_\lambda(\dist^2(\cdot;\, W_{\epsilon,\delta}))(0)=
	\frac{\lambda r_0^2}{1+\lambda}+(1+\lambda)\dist^2(0;\, \co[K_{0,2\epsilon}/(1+\lambda)]).
\]
As $\dist^2(0;\, K)=\dist^2(0;\, W_{\epsilon,\delta})=\dist^2(0;\, K_{0,2\epsilon})=r_0^2$, we then have 
\[
	M_\lambda(0;\, K)\leq r_0^2-\dist^2(0;\, \co[K_{0,2\epsilon}])\,.
\]
Therefore for sufficiently large $\lambda>0$,
\begin{equation} 
	M_{\infty}(0; \,K_0) \leq M_\lambda(0;\, K)\leq M_{\infty}(0;\, K_{0,2\epsilon}).
\end{equation}
Passing to the limit $\lambda\to+\infty$  then gives that for each fixed $\epsilon>0$ small,
\[
	M_{\infty}(0; \, K_0) \leq \liminf_{\lambda\to+\infty}M_\lambda(0;\, K)\leq 
	\limsup_{\lambda\to+\infty}M_\lambda(0;\, K)\leq M_{\infty}(0;\, K_{0,2\epsilon})\,.
\]
Since $K_0$ is compact and $K_{0,2\epsilon}\to K_0$ as $\epsilon\to 0$ under the 
Hausdorff distance in $\mathbb{R}^n$, we also have that 
$\co[K_{0,2\epsilon}] \to \co[K_0]$ as $\epsilon\to 0$ under the Hausdorff distance in $\R^n$.
Thus as $V\mapsto \dist^2(0;\, V)$ is continuous under the Hausdorff distance for compact sets $V\subset\R^n$ \cite{AT04}, 
it follows that $\lim_{\epsilon\to 0+}\dist^2(0;\, \co[K_{0,2\epsilon}])=\dist^2(0;\, \co[K_0])$, and hence $\lim_{\epsilon\to 0+} M_{\infty}(0;\, K_{0,2\epsilon})
= M_{\infty}(0; \, K_0)$. 
Hence $\lim_{\lambda\to+\infty}M_\lambda(0;\, K)$ exists, and
\[
	\lim_{\lambda\to+\infty}M_\lambda(0;\, K)=M_{\infty}(0;\,K_0)=M_{\infty}(0;\, K)\,.
\]
\medskip
It remains to prove \eqref{Sec3.Teo.LimMMA.Prf.Eq3}. First note that when $0<\eta<\delta/2$,
\begin{equation}
\label{ineq1}
	\dist(x;\, S(0;\,r_0+\delta))>\dist(x; K_{0,2\epsilon}),
\end{equation}
because $\dist(x;\, S(0;\,r_0+\delta))=r_0+\delta-|x|\geq r_0+\delta-\eta$
and $\dist(x;\, K_{0,2\epsilon})\leq \dist(0;\, K_{0,2\epsilon})+|x|\leq r_0+\eta$, so that \eqref{ineq1} holds if $ r_0+\eta<r_0+\delta-\eta$, which is equivalent to $2\eta<\delta$.

Now we show that $\dist(x;\, V_{\epsilon,\delta})=\dist(x;\, K_{0,2\epsilon})$. 
Given any point $z_0=(t_0+r_0)y_0/|y_0|\in V_{\epsilon,\delta}\setminus K_{0,2\epsilon}$, 
with $0<t_0\leq \delta$ and $y_0\in K_{0,\epsilon}$, we observe that a necessary condition 
for some $x\in \overline B(0;\,\eta)$ to reach the distance to $V_{\epsilon,\delta}$ at $z_0$, that is,
$\dist(x;\, V_{\epsilon,\delta})=|x-z_0|$, is that the line passing through $z_0$ and $x$ 
does not intersect $K_{0,2\epsilon}$.
Notice that for the point $y_0\in K_{0,\epsilon}$, the $\epsilon$-neighbourhood of $y_0$ 
in $S(0;\,r_0)$ under the geodesic distance $\rho$, given by 
$S_{y_0,\epsilon}:=\{w\in S(0;\,r_0),\,\rho(y_0,w)\leq \epsilon\}$ is contained in 
$K_{0,2\epsilon}$. Therefore if we draw a line passing through $z_0$ and the relative 
boundary of $S_{y_0,\epsilon}$ in $S(0;\,r_0)$ and we can show that the distance 
between the line and the origin $0$ is bounded below by a positive constant uniformly 
with respect to $y_0\in K_{0,\epsilon}$ and $0<t_0\leq \delta$, then we 
can find $0<\eta<\delta/2$, such that
$\dist(x;\, V_{\epsilon,\delta})= \dist(x;\, K_{0,2\epsilon})$
for $x\in \overline B(0;\,\eta)$.

Due to the symmetry of Euclidean balls and spheres, we only need to consider the 
case in $\R^2$ with $y_0=(r_0,0)$, $z=(r_0+t,0),$ where $0<t\leq \delta$ and 
$S_{y_0,\epsilon} =\{(r_0\cos\theta,\,r_0\sin\theta),\; -\epsilon/r_0\leq \theta\leq \epsilon/r_0\}$.
The distance between the line $L$ passing through $z$ and the 
boundary point $(r_0\cos(\epsilon/r_0),\,r_0\sin(\epsilon/r_0))$ 
and the origin $(0,0)$  is attained at a point of the form $(s,u(s))$, where
\[
	u(s)=\frac{r_0\sin(\epsilon/r_0)(r_0+t-s)}{r_0+t-r\cos(\epsilon/r_0) },
\]
so that the squared-distance between $(0,0)$ and a point $(s,u(s))$ in $L$ is
\[
	s^2+\frac{r^2_0\sin^2(\epsilon/r_0)(r_0+t-s)^2}{(r_0+t-r\cos(\epsilon/r_0))^2}\,,
\]
with the minimum point at
\[
	s_0=\frac{r_0^2 \sin^2(\epsilon/r_0)(r_0+t)}{r^2_0\sin^2(\epsilon/r_0)+(r_0+t-r\cos(\epsilon/r_0))^2}\,.
\]
The distance between $(0,0)$ and $L$ is
\[
	\sqrt{s_0^2+u^2(s_0)}\geq |s_0|\geq \frac{r_0^3\sin^2(\epsilon/r_0)}{2(r_0+\delta)^2}:=\eta_0>0\,.
\]
Therefore if we choose $0<\eta<\min\{\delta/2,\, \eta_0\}$, we have, for all $x\in \overline B(0;\,\eta)$, that
$\dist(x;\, V_{\epsilon,\delta})=\dist(x;\, K_{0,2\epsilon})$, and hence
\[
	\dist(x;\, W_{\epsilon,\delta})=\dist(x;\, K_{0,2\epsilon})
\]
for all $x\in \overline B(0;\,\eta)$.
\hfill\qed\\

%%%%%%%%%%%%%%%%%%%%%%%%%%%%%%%%%%%%%%%%%%%%%%%%%%%%%%%%%%%%%%%%%%%%%%%%%%%
%%%%%%%%%%%%%%%%%%%%%%%%%%%%%%%%%%%%%%%%%%%%%%%%%%%%%%%%%%%%%%%%%%%%%%%%%%%
%%%%%%%%%%%%%%%%%%%%%%%%%%%%%%%%%%%%%%%%%%%%%%%%%%%%%%%%%%%%%%%%%%%%%%%%%%%

\noindent {\bf Proof of Proposition \ref{Sec4.Pro.MA2Pnt}:}
This follows from the definition of $M_{\lambda}((x,y);\,K)$ and Lemma \ref{Sec3.Lem.2Pnt}.
\hfill \qed\\

%%%%%%%%%%%%%%%%%%%%%%%%%%%%%%%%%%%%%%%%%%%%%%%%%%%%%%%%%%%%%%%%%%%%%%%%%%%
%%%%%%%%%%%%%%%%%%%%%%%%%%%%%%%%%%%%%%%%%%%%%%%%%%%%%%%%%%%%%%%%%%%%%%%%%%%
%%%%%%%%%%%%%%%%%%%%%%%%%%%%%%%%%%%%%%%%%%%%%%%%%%%%%%%%%%%%%%%%%%%%%%%%%%%
 
\noindent {\bf Proof of Theorem \ref{Sec4.Teo.HAS}:} 
Let $\mu=\dist_{\mathcal{H}}(K;\,L)$ with $\mu$ finite since $K$ and $L$ are compact sets. By Definition 
\eqref{Sec2.Def.HausDist} for Hausdorff distance, we have for $x\in \mathbb{R}^n$
\begin{equation}\label{Eq.Sec4.Teo.HAS.1}
	|\dist(x;\,K)-\dist(x;\,L)|\leq \mu\,,
\end{equation}
hence, 
\begin{equation*}
	\begin{split}
		\dist^2(x;\,K)	&\leq \mu^2+2\mu\dist(x;\,L)+\dist^2(x;\,L) \\[1.5ex]
				&\leq \mu^2+\mu(1+\dist^2(x;\,L))+\dist^2(x;\,L)=\mu(1+\mu)+(1+\mu)\dist^2(x;\,L)\,.
	\end{split}
\end{equation*}
After adding $\lambda|x|^2$ to both sides and taking the convex envelope, we find
\begin{equation*}
	\co[\dist^2(\cdot;\,K)+\lambda|\cdot|^2](x)\leq \mu(1+\mu) + 
		(1+\mu)\co[\dist^2(\cdot;\,L)+\frac{\lambda}{1+\mu}|\cdot|^2](x)
\end{equation*}
which yields
\begin{equation}\label{Eq.Sec4.Teo.HAS.2}
	C_{\lambda}^l(\dist^2(\cdot;\,K))(x)\leq 
		\mu(1+\mu) + 
		(1+\mu)C_{\lambda/(1+\mu)}^l(\dist^2(\cdot;\,L))(x)\,.
\end{equation}
Since
\begin{equation*}
	C^l_{\lambda/(1+\mu)}(\dist^2(\cdot;\,L))(x)\leq C^l_\lambda(\dist^2(\cdot;\,L))(x)\text{ and }
	C^l_{\lambda/(1+\mu)}(\dist^2(\cdot;\,L))(x)\leq \dist^2(x,\,L)\,,
\end{equation*}
we obtain, from \eqref{Eq.Sec4.Teo.HAS.2}, after using \eqref{Eq.Sec4.Teo.HAS.1}, that
\begin{equation}\label{Eq.Sec4.Teo.HAS.3}
\begin{split}
	C^l_\lambda(\dist^2(\cdot;\,K))(x)&\leq \mu(1+\mu) +C_{\lambda}^l(\dist^2(\cdot;\,L))(x) + \mu \dist^2(x;\,L) \\[1.5ex]
				   &\leq C_{\lambda}^l(\dist^2(\cdot;\,L))(x) +	 \mu(1+\mu) + \mu (\mu+\dist(x;\,K))^2)\,.
\end{split}
\end{equation}
With a similar argument, we find that 
\begin{equation}\label{Eq.Sec4.Teo.HAS.4}
\begin{split}
	C^l_\lambda(\dist^2(\cdot;\,L))(x)&\leq \mu(1+\mu) +C_{\lambda}^l(\dist^2(\cdot;\,K))(x) + \mu \dist^2(x;\,K) \\[1.5ex]
				   &\leq C_{\lambda}^l(\dist^2(\cdot;\,K))(x) +	 \mu(1+\mu) + \mu (\mu+\dist(x;\,K))^2)\,.
\end{split}
\end{equation}
By comparing \eqref{Eq.Sec4.Teo.HAS.3} and  \eqref{Eq.Sec4.Teo.HAS.4} we therefore conclude that given a compact 
set $K\subset \R^n$, for any compact set $L\subset \R^n$, we have that, for any $x\in \R^n$,
\begin{equation}\label{Eq.Sec4.Teo.HAS.5}
	\left|C^l_\lambda(\dist^2(\cdot;\,\,K))(x) - C_{\lambda}^l(\dist^2(\cdot;\,L))(x)\right|
	\leq \mu\Big((1+\mu) + (\mu+\dist(x;\,K))^2)\Big)\,,
\end{equation}
which proves \eqref{Sec4.Eq.Teo.HASLwtr}.
To show \eqref{Sec4.Eq.Teo.HAS},
observe that after using \eqref{Eq.Sec4.Teo.HAS.1} we have for any $x\in \R^n$
\begin{equation}\label{Eq.Sec4.Teo.HAS.6}
\begin{split}
	|\dist^2(x;\,K)-\dist^2(x;\,L)|&\leq	|\dist(x;\,K)-\dist(x;\,L)||\dist(x;\,K)+\dist(x;\,L)| \\[1.5ex]
				       &\leq	\mu(2\dist(x;\,K)+\mu)\,,
\end{split}
\end{equation}
and from the definition of the multiscale medial axis map and the triangle inequality, we obtain
\begin{equation}\label{Eq.Sec4.Teo.HAS.7}
	\left|M_{\lambda}(x;\,K)-M_{\lambda}(x;\,L)\right|\leq \mu(1+\lambda)
	\Big((\dist(x;\,K)+\mu)^2+2\dist(x;\,K)+2\mu+1\Big)\,,
\end{equation}
where we have taken into account \eqref{Eq.Sec4.Teo.HAS.5} and \eqref{Eq.Sec4.Teo.HAS.6}. This concludes the proof.
\hfill \qed\\

%%%%%%%%%%%%%%%%%%%%%%%%%%%%%%%%%%%%%%%%%%%%%%%%%%%%%%%%%%%%%%%%%%%%%%%%%%%
%%%%%%%%%%%%%%%%%%%%%%%%%%%%%%%%%%%%%%%%%%%%%%%%%%%%%%%%%%%%%%%%%%%%%%%%%%%
%%%%%%%%%%%%%%%%%%%%%%%%%%%%%%%%%%%%%%%%%%%%%%%%%%%%%%%%%%%%%%%%%%%%%%%%%%%

\noindent {\bf Proof of Corollary \ref{Sec4.Cor.HAS}:}
This follows from Theorem \ref{Sec4.Teo.HAS},  since
$
	\dist(x;\,\partial \Omega)\leq \diam(\Omega)
$
if $x\in \bar{\Omega}$.
\hfill \qed\\

%%%%%%%%%%%%%%%%%%%%%%%%%%%%%%%%%%%%%%%%%%%%%%%%%%%%%%%%%%%%%%%%%%%%%%%%%%%%%%%%%%%%%%%
%%%%%%%%%%%%%%%%%%%%%%%%%%%%%%%%%%%%%%%%%%%%%%%%%%%%%%%%%%%%%%%%%%%%%%%%%%%%%%%%%%%%%%%
%%%%%%%%%%%%%%%%%%%%%%%%%%%%%%%%%%%%%%%%%%%%%%%%%%%%%%%%%%%%%%%%%%%%%%%%%%%%%%%%%%%%%%%

{\bf Acknowledgement.}  
We thank the referees for valuable suggestions.
KZ wishes to thank The University of Nottingham for its support,
EC is grateful for the financial support of the College of Science, Swansea University,
and
AO acknowledges the financial support of the Argentinean Agency through the 
Project Prestamo BID PICT PRH 30 No 94, the National University of Tucum\'{a}n through the project PIUNT E527
and the Argentinean Research Council CONICET.

\end{document}